\documentclass[10pt,reqno,a4paper]{article}
\usepackage{amsmath,amsfonts,amscd,amsthm,amssymb}
\usepackage{graphicx}
\usepackage{subcaption}
\usepackage[dvipsnames]{xcolor}
\usepackage{cite}
\usepackage{tikz}
\usepackage[colorlinks=true]{hyperref}\hypersetup{urlcolor=blue, citecolor=blue}
\usepackage{scalerel,stackengine}
\usepackage{bm}
\usepackage{enumerate}
\usepackage{multicol}
\allowdisplaybreaks

\newcommand{\R}{\mathbb{R}}

\newcommand{\pt}{\partial_t}

\newcommand{\bw}{{\bm w}}

\newcommand{\bP}{{\bm P}}
\newcommand{\bv}{{\bm v}}

\newcommand{\bS}{{\bm S}}
\newcommand{\bV}{{\bm V}}

\newcommand{\bH}{{\bm H}}
\newcommand{\bL}{{\bm L}}

\newcommand{\f}{{\bm f}}

\newcommand{\bu}{{\bm u}}

\newcommand{\bn}{{\bm n}}

\newcommand{\bphi}{\mbox{\boldmath $\phi$}}

\newcommand{\bxi}{\mbox{\boldmath $\xi$}}

\newcommand{\btheta}{\mbox{\boldmath $\theta$}}

\newcommand{\trho}{\tilde{\rho}}

\newcommand{\e}{{\bm e}}

\newcommand{\Eh}{\mathcal{E}_h}
\newcommand{\dx}{\mathrm{\,d}x}
\newcommand{\ds}{\mathrm{\,d}s}
\newcommand{\dt}{\mathrm{\,d}t}

\newcommand{\jump}[1]{[\kern-0.85ex[ #1 ]\kern-0.85ex]} %
\newcommand{\avr}[1]{\{\kern-0.85ex\{ #1 \}\kern-0.85ex\}}
	
\DeclareSymbolFont{myletters}{OML}{ztmcm}{m}{it}
\DeclareMathSymbol{\uplambda}{\mathord}{myletters}{"15}

\numberwithin{equation}{section}

\begin{document}

\newtheorem{theorem}{Theorem}[section]
\newtheorem{lemma}{Lemma}[section]
\newtheorem{remark}{Remark}[section]
\newtheorem{example}{Example}[section]
\newtheorem{proposition}{Proposition}[section]

\def\cydot{\leavevmode\raise.4ex\hbox{.}}

\title
{ \large\bf On a Completely Discrete Discontinuous Galerkin Method for Incompressible Chemotaxis-Navier-Stokes Equations}

%[DG FEM for the Chemotaxis-Navier-Stokes]

\author{Bikram Bir\thanks{Department of Mathematics, 
		Indian Institute of Technology Bombay, Powai, Mumbai-400076, India. Email: bikram@math.iitb.ac.in}, ~~	  
	  Harsha Hutridurga\thanks{Department of Mathematics, 
	  	Indian Institute of Technology Bombay, Powai, Mumbai-400076, India. Email: hutri@math.iitb.ac.in}
	  	~~and
	Amiya K. Pani\thanks{Department of Mathematics, 
		BITS-Pilani, KK Birla Goa Campus, NH 17 B, Zuarinagar Goa-403726, India. Email: amiyap@goa.bits-pilani.ac.in,  akp@math.iitb.ac.in}}

\date{}
\maketitle

\begin{abstract}
	This paper deals with a fully discrete numerical scheme for the incompressible Chemotaxis(Keller-Segel)-Navier-Stokes system. Based on a discontinuous Galerkin finite element scheme in the spatial directions, a semi-implicit first-order finite difference method in the temporal direction is applied to derive a completely discrete scheme. With the help of a new projection, optimal error estimates in $L^2$ and $H^1$-norms for the cell density, the concentration of chemical substances and the fluid velocity are derived. Further, optimal error bound in $L^2$-norm for the fluid pressure is obtained. Finally, some numerical simulations are performed, whose results confirm the theoretical findings.
\end{abstract}

\vspace{1em} 

\vspace{1em} 

\section{Introduction}

In this paper, a completely discrete method is applied to the following time-dependent Chemotaxis-Navier-Stokes system in $\Omega\times(0,T]$:
\begin{align}
	&\rho_t - \mu \Delta\rho + \bu\cdot\nabla\rho + \beta \nabla\cdot(\rho\nabla c) = 0,\quad\mbox {in}~ \Omega\times(0,T], \label{eqrho}\\
	&c_t - \kappa\Delta c + \bu\cdot\nabla c + \gamma\rho	c = 0,\quad\mbox {in}~ \Omega\times(0,T], \label{eqc}\\
	&\bu_t - \nu\Delta\bu + \bu\cdot\nabla\bu + \nabla p = \rho\nabla\Phi, \quad\mbox {in}~ \Omega\times(0,T], \label{equ}\\
	&\nabla \cdot \bu=0,\quad\mbox {in}~ \Omega\times(0,T], \label{eqdiv}\\
	&\nabla\rho\cdot n = 0,~~ \nabla c\cdot n = 0, ~~\bu = 0,\quad\mbox {on}~\partial\Omega\times (0,T], \label{eqbd}\\
	&\rho(\cdot,0) = \rho_0,~~ c(\cdot,0) = c_0,~~\bu(\cdot,0) = \bu_0,\quad\mbox {in}~\Omega. \label{eqint}
\end{align}
Here, $\Omega$ is a bounded polygonal domain in $\mathbb{R}^2$ with boundary $\partial \Omega$. $\rho(x,t)$, $c(x,t)$, $\bu(x,t)$ and $p(x,t)$ denote the cell density, the concentration of chemical substances, the fluid velocity and the pressure of the fluid, respectively. The constants $\mu>0, \kappa>0, \beta>0,\gamma>0$ and $\nu>0$ represent the diffusion constant of the cell, the diffusion constant of the chemical substances, the consumption rate of the chemical substances, the Chemotactic coefficient and the kinematic coefficient of viscosity of the fluid, respectively. The term $\nabla\Phi$ represents the effect due to the density variation caused by cell aggregation. Further, the initial data $\rho_0, c_0, \bu_0$ are given functions in their respective domains of definition. 

The above-mentioned model is proposed by Tuval \textit{et. al.} \cite{TCDWKG05} in the context of modeling the movement of cells towards the concentration gradient in the fluid region. 
An experimental observation reveals that the mechanism is a chemotactic movement of the swimming bacteria in the direction of a higher concentration of oxygen that they consume, a gravitational effect on the motion of the fluid by the heavier bacteria, and a convective transport of both cells and oxygen through the water (see \cite{Lor10} and \cite{TCDWKG05}).
Due to the presence of the coupling between incompressible Navier-Stokes equation and the Chemotaxis system, the existence and uniqueness of the solution of \eqref{eqrho}-\eqref{eqint} are quite challenging. However, numerous researchers have dedicated their efforts to demonstrating the existence, uniqueness, regularity, and qualitative properties of the continuous solutions of the above system. For an extensive exploration of these aspects, we refer to \cite{CL16,DLM10,LL11,Win12,Wan23,ZZ14} and the references, therein.

In literature, only a limited number of numerical studies are available pertaining to this model. For instance, Chertock et al. \cite{CFKLM12} have presented a high-resolution, vorticity-based hybrid finite volume finite difference scheme designed for the examination of the nonlinear dynamics inherent in a two-dimensional chemotaxis-fluid system. 
Lee and Kim \cite{LK15} have investigated numerically falling bacterial plumes caused by bioconvection in a three-dimensional chamber by solving the full chemotaxis–fluid coupled system using an operator splitting-type Navier–Stokes solver such as the advection term is addressed using a semi-Lagrangian method, while the diffusion term is tackled using the backward Euler method. A fully discrete mixed finite element method by introducing a new variable by the gradient of the chemical concentration to the nonlinear chemo-attraction term in the cell-density equation has been analyzed by Duarte-Rodr{\'\i}guez \textit{et. al.} \cite{DRRV21}. In the work by Caucao \textit{et. al.}\cite{CCGI23}, a Banach-spaces-based fully-mixed finite element method has been examined for a stationary Chemotaxis-Navier-Stokes system. In \cite{LXLC23}, Li \textit{et. al.} have discussed a continuous Galerkin finite element method with a semi-implicit first-order backward Euler approach for the model described above. A numerical study with an upwind finite element method to investigate the pattern formation and the hydrodynamical stability of the system for our model has been discussed by Deleuze \textit{et. al.} \cite{DCTS16}  and a time scale analysis has been performed. In \cite{TZL23}, a linear, decoupled fully discrete unconditional energy stable finite element method combined with scalar auxiliary variable (SAV) approach, implicit-explicit (IMEX) scheme, and pressure-projection method has been discussed for the chemo-repulsion-Navier-Stokes system with linear production term.  A fully discrete discontinuous Galerkin finite element method with an implicit explicit approach for nonlinear and coupling term and a pressure correction method for the Navier-Stokes equation is applied for the chemo-repulsion-Navier-Stokes system by Wang \textit{et. al.} \cite{WZWZ23}. 

In this paper, we have analyzed a fully discrete scheme based on the discontinuous Galerkin finite element method in the spatial directions along with a first-order semi-implicit finite difference scheme (backward Euler method) in the temporal direction. The discontinuous Galerkin finite element method, initially introduced in the 1970s by Douglas and Dupont \cite{DD75}, Wheeler \cite{Whe78}, Arnold \cite{Arn82}, and Reed and Hill \cite{RH73} for the elliptic, parabolic and hyperbolic problems, has gained significant popularity in the numerical approximation of a diverse range of mathematical problems. Its widespread adoption can be attributed to its versatility and notable characteristics, including local mesh adaptivity, local mass conservation, arbitrary order accuracy, and a relatively straightforward implementation when compared to alternative methods such as the finite volume technique. A substantial body of research has been dedicated to the application of the discontinuous Galerkin method to problems involving the Stokes and the Navier-Stokes equations, with notable references available in the literature \cite{GRW05domain,GRW05spliting,MLR23,BGR23} and references, therein. On the other hand, the use of the discontinuous Galerkin approach in the context of the Keller-Segel system has received a growing interest now-a-days, see, \cite{Eps09, EI09, EK08, LSY17, QLY23}.
Nevertheless, and up to our knowledge, there is no work based on discontinuous Galerkin finite element methods for the coupled chemotaxis-Navier-Stokes equations \eqref{eqrho}-\eqref{eqint}. Although, in \cite{WZWZ23}, a discontinuous Galerkin method has been employed for the chemo-repulsion-Navier-Stokes system the error bounds are sub-optimal in $L^2$-norm. As mentioned earlier, the goal of this paper is to find optimal error estimates in $L^2$-norm for the incompressible coupled time-dependent chemotaxis-Navier-Stokes system. 
Also, for the Keller-Segel model, there are no optimal error bounds in $L^2$-norm in the literature except \cite{LSY17}. In \cite{LSY17}, a local discontinuous Galerkin finite element has been applied to the chemotaxis Keller-Segel model and an optimal error estimate has been derived under the assumption of high regularity of the exact solution and using the special finite element spaces (for $\mathcal{Q}_k$-polynomial spaces). 

Due to the presence of chemotactic term ($\beta \nabla\cdot(\rho\nabla c)$) in the equation \eqref{eqrho}, a use of the standard elliptic projection operator shows a sub-optimal results in the error estimates with respect to space variable (see, \cite[Lemma 6]{LXLC23}). We can derive optimal estimates with the help of a modified version of the elliptic projection in the discontinuous space using a linearized form of the chemotactic term.

The main contributions of this work are as follows:
\begin{enumerate}[(i)]
	\item { Discrete mass conservation of the cell density}
	\item Introduce a modified version of the elliptic projection operator in the discontinuous space which helps to find optimal error estimates
	\item Optimal \textit{a priori} error estimates in $L^2$ and $H^1$-norms for the fully discrete cell density, concentration and fluid velocity, and in $L^2$-norm for the fluid pressure
	\item Numerical experiments with known solutions to verify the rate of convergence. Another example is to show the movement of the swimming bacteria toward the concentration in a fluid medium.
\end{enumerate}

Although optimal $L^2$ error estimates are derived for the discontinuous Galerkin methods, this present analysis can easily extended to classical finite element method for deriving optimal $L^2$ error estimates. 
The rest of the paper is organized as follows: In Section 2, we define some functional spaces, broken Sobolev spaces and associated norms, bilinear and trilinear forms, continuous trace inequalities, and some regularity assumptions. Section 3 contains finite element spaces, discrete formulation, some properties of the bilinear and trilinear forms, discrete trace inequalities, and projection operators and their properties. \textit{A priori} error estimates are discussed in Section 4. In Section 5, we present a few numerical examples for validating our findings. Finally, the conclusion part is given in Section 6.

Throughout the paper, we will use $C$ as a generic constant which does not depend on the discretizing parameter $h$ and $\Delta t$ but may depend on the given data. Let $a\lesssim b$ denotes $a\le Cb$.

\section{Preliminaries}

We begin this section by introducing some standard functional spaces.
For a non-negative integer $m$ and $p$ with $0\le p\le \infty$, we denote by $W^{m,p}(\Omega)$ and $L^p(\Omega)$ the usual Sobolev spaces and the Lebesgue spaces with associated norms $\|\cdot\|_{W^{m,p}}$ and $\|\cdot\|_{L^p}$, respectively \cite{A75}. If $p=2$, $W^{m,2}(\Omega)=H^m(\Omega)$ is the Hilbert spaces equipped with norm $\|\cdot\|_{H^m}$. We define $(\cdot,\cdot)$ as the $L^2$-inner product equipped with norm $\|\cdot\|$.
Let $H^m/\mathbb{R}$ be the quotient space of equivalent classes of functions in $H^m(\Omega)$ differ by constant with norm $\|\phi\|_{H^m/\mathbb{R}}=\inf_{c\in \mathbb{R}}\|\phi+c\|_{H^m}.$ On the other way, one can define $H^m/\mathbb{R}$ as $$H^m/\mathbb{R} = \left\{\psi\in H^m(\Omega): \int_\Omega \psi(s) \rm{d}s = 0\right\}.$$
For $m=0$, it is denoted by $L^2/\mathbb{R}$.
Further, $H_0^1(\Omega)$ is defined as a subspace of $H^1(\Omega)$ whose elements vanish on $\partial\Omega$ in the sense of trace. For simplicity, we denote $H^m(\Omega)$ as $H^1$.

\noindent
For our subsequent analysis, we denote the vector-valued function spaces by boldface letters such as
\begin{align*}
	\bH^m=[H^m(\Omega)]^2, \quad \bH_0^1 = [H_0^1(\Omega)]^2, \quad\mbox{and}\quad\bL^2 = [L^2(\Omega)]^2.
\end{align*}

For dG formulation, we consider a shape regular family of triangulations $\Eh$ of $\bar{\Omega}$ \cite{Cia78}, consisting of triangles or quadrilaterals of maximum diameter $h$. Let $h_E$ denote the diameter of a triangle $E$ and $\pi_E$ the diameter of its inscribed circle.  Then, the shape regular triangulation means there exists a constant $m_1>0$ such that
\begin{equation}
	\frac{h_E}{\pi_E}\leq m_1, \quad \forall~ E\in \Eh.\label{regular}
\end{equation}
We denote by $\Gamma_h^0$ the set of all interior edges of $\Eh$ and set $\Gamma_h = \Gamma_h^0\cup\partial\Omega$. With each edge $e$, we associate a unit normal vector $\bn_e$. If $e$ is on the boundary $\partial\Omega$, then $\bn_e$ is taken to be the unit outward vector normal to $\partial \Omega$. Let $e$ be an edge shared by two elements $E_i$ and $E_j$ of $\Eh$; we associate with $e$, once and for all, a unit normal vector $\bn_e$ directed from $E_i$ to $E_j$. We define the jump $\jump{\psi}$ and the average $\avr{\psi}$ of a function $\psi$ on $e$ by
\[ \jump{\psi} = (\psi|_{E_i})-(\psi|_{E_j}), \quad \avr{\psi}=\frac{1}{2}\left((\psi|_{E_i})+(\psi|_{E_j})\right).\] 
If $e$ is adjacent to $\partial\Omega$, then the jump and the average of $\psi$ on $e$ coincide with the value of $\psi$ on $e$. Similar definitions are considered for vector-valued functions.

\noindent
We consider the following discontinuous spaces \cite{GRW05domain} (broken Sobolev spaces) for our analysis.
\begin{align*}
	X &= \left\{ \psi\in L^2 : \psi|_E \in W^{2, \frac{4}{3}}(E)~~\forall~ E\in\Eh \right\}, \\
	X^0 &= \left\{ \chi\in L^2/\R : \chi|_E \in W^{2, \frac{4}{3}}(E)~~\forall~ E\in\Eh \right\},\\
	\bV &= \left\{ \bv\in \bL^2 : \bv|_E \in \left(W^{2, \frac{4}{3}}(E)\right)^2~~\forall~ E\in\Eh \right\},\\
	M &= \left\{ q\in L^2/\R : q|_E \in W^{1, \frac{4}{3}}(E)~~\forall~ E\in\Eh \right\}.
\end{align*}
We also define the associated norms of the above spaces as follows:
\begin{align*}
	\|\chi\|_{\rho} &= \left(\sum_{E\in \Eh} \|\nabla\chi\|_{L^2(E)}^{2} + \sum_{e\in \Gamma_h} \frac{\sigma_\rho}{h_e} \|[\chi]\|_{L^2(e)}^{2}\right)^\frac{1}{2} \quad \forall ~\chi\in X^0,\\
	\|\psi\|_{c} &= \left(\sum_{E\in \Eh} \|\nabla\psi\|_{L^2(E)}^{2} + \sum_{e\in \Gamma_h} \frac{\sigma_c}{h_e} \|[\psi]\|_{L^2(e)}^{2}\right)^\frac{1}{2} \quad \forall~ \psi\in X,\\
	\|\bv\|_{u} &= \left(\sum_{E\in \Eh} \|\nabla\bv\|_{L^2(E)}^{2} + \sum_{e\in \Gamma_h} \frac{\sigma_u}{h_e} \|[\bv]\|_{L^2(e)}^{2}\right)^\frac{1}{2} \quad \forall~ \bv\in \bV,
\end{align*}
where $h_e = |e|: $ length of edge and $\sigma_\rho, \sigma_c, \sigma_u >0$ are the penalty parameters. 

\noindent
We define the bilinear forms: $a_\alpha: X\times X\to\R$ with $\alpha=\{\rho,c\}$, $a_u:\bV\times\bV\to\R$, and $d:\bV\times M\to\R$ as
\begin{align}
	a_\alpha(\psi,\phi) = & \sum_{E\in \Eh} \int_E \nabla\psi \cdot \nabla\phi \dx - \sum_{e\in\Gamma_h^0} \int_e \avr{\nabla\psi}\cdot\bn_e\jump{\phi}\ds \nonumber\\
	&\qquad\qquad - \sum_{e\in\Gamma_h^0} \int_e \avr{\nabla\phi}\cdot\bn_e \jump{\psi} \ds + \sum_{e\in\Gamma_h^0} \frac{\sigma_\alpha}{h_e}\int_e \jump{\psi}\jump{\phi}\ds \quad\forall~ \psi, \phi \in X, \label{biac}\\
	a_u(\bv,\bw) = & \sum_{E\in \Eh} \int_E \nabla\bv : \nabla\bw \dx - \sum_{e\in\Gamma_h} \int_e \avr{\nabla\bv}\bn_e\cdot\jump{\bw}\ds \nonumber\\
	&\quad - \sum_{e\in\Gamma_h} \int_e \avr{\nabla\bw}\bn_e\cdot\jump{\bv} \ds + \sum_{e\in\Gamma_h} \frac{\sigma_u}{h_e}\int_e \jump{\bv}\cdot\jump{\bw}\ds \quad\forall~ \bv, \bw \in \bV, \label{biau} \\
	d(\bv,q) = & \sum_{E\in\Eh} \int_E q \nabla\cdot \bv \dx - \sum_{e\in\Gamma_h} \int_e \avr{q}\bn_e\cdot\jump{\bv} \ds \nonumber\\
	= & -\sum_{E\in\Eh} \int_E \bv \cdot\nabla q \dx + \sum_{e\in\Gamma_h} \int_e \avr{\bv}\cdot\bn_e \jump{q} \ds \quad \forall~ \bv\in \bV, q\in M.
\end{align}

\noindent
The trilinear forms $b_1:\bV \times X \times X \to \R$ and $b_2:\bV \times \bV \times \bV \to \R$ are defined as follows:
\begin{align}
	b_1(\bv,\psi,\phi) =& \sum_{E\in\Eh} \int_E (\bv\cdot\nabla\psi)\phi \dx + \frac{1}{2}\sum_{E\in\Eh}\int_E (\nabla\cdot\bv)\psi\phi \dx \nonumber\\
	&\quad - \sum_{e\in\Gamma_h^0} \int_{e} \avr{\bv}\cdot\bn_e \jump{\psi}\avr{\phi} \ds  - \frac{1}{2}\sum_{e\in\Gamma_h^0}\int_e \jump{\bv}\cdot\bn_e \avr{\psi\phi} \ds \quad\forall~ \bv\in\bV, \psi,\phi\in X, \label{tribc}\\
	b_2(\bv,\bw,\bphi) =& \sum_{E\in\Eh}  \int_E (\bv\cdot\nabla)\bw\cdot\bphi \dx + \frac{1}{2}\sum_{E\in\Eh}\int_E (\nabla\cdot\bv)\bw\cdot\bphi \dx  \nonumber\\
	&~~ - \sum_{e\in\Gamma_h} \int_{e} \avr{\bv}\cdot\bn_e \jump{\bw}\cdot\avr{\bphi} \ds - \frac{1}{2}\sum_{e\in\Gamma_h}\int_e \jump{\bv}\cdot\bn_e \avr{\bw\cdot\bphi} \ds \quad\forall~\bv, \bw, \bphi \in \bV. \label{tribu}
\end{align}

\noindent 
We also define another trilinear form for the chemotaxic term $g:X\times X\times X\to\R$ as
\begin{align}
	g(\chi,\psi,\phi) &= \sum_{E\in\Eh} \int_E \chi\nabla\psi\cdot\nabla\phi \ds - \sum_{e\in\Gamma_h^0} \int_e \avr{\chi\nabla\psi}\cdot\bn_e \jump{\phi} \ds\nonumber\\
	&\qquad - \sum_{e\in\Gamma_h^0} \int_e \avr{\chi\nabla\phi}\cdot\bn_e \jump{\psi} \ds \quad \forall~ \chi, \psi,\phi \in X. \label{defg}
\end{align} 

Now, the dG formulation of \eqref{eqrho}-\eqref{eqint} reads as: Find $(\trho, c, \bu, p)\in X^0 \times X \times \bV \times M$ such that for all $t>0$
\begin{align}
	&(\trho_t,\chi) + \mu a_\rho(\trho,\chi) + b_1(\bu,\trho,\chi) - \beta g(\trho+m_0, c, \chi)  = 0 \quad\forall~\chi\in X^0, \label{dg1rho}\\
	&(c_t,\psi) + \kappa a_c(c,\psi) + b_1(\bu, c, \psi) + \gamma((\trho+m_0)c,\psi) = 0 \quad\forall~\psi\in X, \label{dg2fc}\\
	&(\bu_t,\bv) + \nu a_u(\bu,\bv) + b_2(\bu,\bu,\bv) - d(\bv,p) = ((\trho+m_0)\nabla\Phi,\bv)  \quad\forall~\bv\in\bV, \label{dg3u}\\
	&d(\bu,q)=0 \quad\forall~ q\in M, \label{dg4div}
\end{align}
where $\trho = \rho-m_0$ is a new variable with $m_0=\frac{1}{|\Omega|}\int_\Omega\rho_0\dx$. 

For our subsequence analysis, we will use the following standard trace inequalities (\cite[Section 2.1.3]{Riv08}).
\begin{lemma}\label{edge2elem}
	For each element $E\in\Eh$ with diameter $h_E$, there exists a constant $C$ independent of $h_E$ such that the followings hold:
	\begin{align*}
		\|\phi\|_{L^2(e)} & \le C\left(h_E^{-\frac{1}{2}}\|\phi\|_{L^2(E)} + h_E^{\frac{1}{2}}\|\nabla\phi\|_{L^2(E)}\right) \quad \forall~\phi\in H^1(E),\\
		\|\nabla\phi\|_{L^2(e)} & \le C\left(h_E^{-\frac{1}{2}}\|\nabla\phi\|_{L^2(E)} + h_E^{\frac{1}{2}}\|\nabla^2\phi\|_{L^2(E)}\right) \quad \forall~\phi\in H^2(E),
	\end{align*}
	where $e$ is an edge of the element $E$. These results are also valid for vector-valued function $\bphi\in \bV$.
\end{lemma}
\noindent
We will subsequently use the Gagliardo-Nirenberg inequality \cite{HS00}
\begin{equation} \label{GNI}
	\|\phi\|_{L^p} \le C \|\phi\|^{2/p}\|\nabla\phi\|^{1-2/p} \quad\forall~\phi\in H^1(E),
\end{equation}
where $2\le p<\infty$ and $C=C(p,E)$. Also, we will consider the Agmon's inequality \cite{HS00}
\begin{equation}\label{AI}
	\|\phi\|_{L^\infty} \le C \|\phi\|^{1/2}\|\nabla^2\phi\|^{1/2} \quad\forall~\phi\in H^2(E),
\end{equation}
where $C=C(E)$.

Let us take the assumption on the exact solution.\\
(\textbf{A1}) Suppose the exact solutions satisfy the following regularity assumptions:
\begin{align*}
	& \bu\in L^\infty(0,T;\bH^2)\cap  L^2(0,T;\bH^{k+1}),~\bu_t\in L^2(0,T;\bH^{k+1}), ~\bu_{tt}\in L^2(0,T;\bL^2), \\
	& \rho\in L^\infty(0,T;H^2)\cap  L^2(0,T;H^{k+1}),~\rho_t\in L^2(0,T;H^{k+1}), ~\rho_{tt}\in L^2(0,T;L^2),\\
	& c\in L^\infty(0,T;H^3)\cap  L^2(0,T;H^{k+1}),~c_t\in L^2(0,T;H^{k+1}), ~c_{tt}\in L^2(0,T;L^2),\\
	& p\in L^2(0,T;H^{k}),~p_t\in L^2(0,T;H^{k}),~\Phi\in L^\infty(0,T;W^{1,4}(\Omega)^2).
\end{align*}

\section{Discrete Formulation}

In this section, we first define discrete finite element spaces, semidiscrete weak formulation, and some properties of bilinear and trilinear operators. Then, we define a few projection operators and prove their properties. Finally, we discuss the fully discrete weak formulation. 

\subsection{Semidiscrete dG Formulation}
For the semidiscrete formulation, we first consider the non-conforming finite element spaces $X^0_h \times X_h \times \bV_h \times M_h \subset X^0 \times X \times \bV \times M$, which are defined as follows: For any integer $k\ge 1$,
\begin{align*}
	X_h &= \left\{ \psi_h\in L^2 : \psi_h|_E \in \mathbb{P}_{k}(E) ~~\forall~ E\in\Eh \right\}, \\
	X^0_h &= \left\{ \chi_h\in L^2/\R : \chi_h|_E \in \mathbb{P}_{k}(E) ~~\forall~ E\in\Eh \right\},\\
	\bV_h &= \left\{ \bv_h\in \bL^2 : \bv_h|_E \in \left(\mathbb{P}_{k}(E)\right)^2 ~~\forall~ E\in\Eh \right\},\\
	M_h &= \left\{ q_h\in L^2/\R : q_h|_E \in \mathbb{P}_{k-1}(E) ~~\forall~ E\in\Eh \right\},
\end{align*} 
where $\mathbb{P}_k(E)$ is the polynomial space of order less or equal to $k$ over $E$.

Further, we will assume that the meshes are quasi-uniform and the following inverse hypothesis holds for all discrete $\phi_h\in X_h ~\mbox{or}~ X_h^0 ~\mbox{or}~ \bV_h$, see \cite[Theorem 3.2.6]{Cia78}
\begin{align}\label{inv.hypo}
	\| \phi_h\|_{W^{m,p}(E)^d} \leq  Ch^{n-m-d(\frac{1}{q}-\frac{1}{p})} \|\phi_h\|_{W^{n,q}(E)^d},
\end{align}
where $0\le n \le m \le 1$, $0\le q \le p \le \infty$, $h$ be the diameter of the mesh cell $E\in\mathcal{E}_h$ and $\| \cdot\|_{W^{m,p}(E)^d}$ is the norm in Sobolev space $W^{m,p}(E)^d$. In our case $d=2$.

The semidiscrete dG formulation of \eqref{dg1rho}-\eqref{dg4div} is given as: Find $(\trho_h, c_h, \bu_h, p_h)\in X_h^0 \times X_h \times \bV_h \times M_h$ such that for all $t>0$
\begin{align}
	&(\trho_{ht},\chi_h) + \mu a_\rho(\trho_h,\chi_h) + b_1(\bu_h,\trho_h,\chi_h) - \beta g(\trho_h+m_0, c_h, \chi_h) = 0 \quad\forall~\chi_h\in X_h^0, \label{smdg1rho}\\
	&(c_{ht},\psi_h) + \kappa a_c(c_h,\psi_h) + b_1(\bu_h, c_h, \psi_h) + \gamma((\trho_h+m_0)c_h,\psi_h) = 0 \quad\forall~\psi_h\in X_h, \label{smdg2fc}\\
	&(\bu_{ht},\bv_h) + \nu a_u(\bu_h,\bv_h) + b_2(\bu_h,\bu_h,\bv_h) - d(\bv_h,p_h) = ((\trho_h+m_0)\nabla\Phi,\bv_h)  \quad\forall~\bv_h\in\bV_h, \label{smdg3u}\\
	&d(\bu_h,q_h)=0 \quad\forall~ q_h\in M_h, \label{smdg4div}
\end{align}
with initial data $\trho_h(0), c_h^0(0), \bu_h^0(0)$ are the appropriate approximation of $\trho_0, c_0, \bu_0$ respectively.

To show the well-posedness of the above discrete system, first, we discuss the following properties:
\begin{lemma}[Discrete coercivity{\cite[Lemma 4.12]{DE11}}]\label{lem:prop_bilinear}
	For some sufficiently large positive constant $\sigma_\alpha^*$ with $\sigma_\alpha >\sigma_\alpha^*$, $\alpha=\rho,c,u$, there exists some $\gamma^*>0$ independent of $h$ such that the following hold
	\begin{align*}
		a_\rho(\chi_h,\chi_h) \ge \gamma^*\|\chi_h\|_{\rho}^2 \quad \chi_h\in X_h^0,\\
		a_c(\psi_h,\psi_h) \ge \gamma^*\|\psi_h\|_{c}^2 \quad \psi_h\in X_h,\\
		a_u(\bv_h,\bv_h) \ge \gamma^*\|\bv_h\|_{u}^2 \quad \bv_h\in \bV_h.
	\end{align*}
	Additionally, for any $C>0$ independent of $h$, the following continuity properties \cite[Lemma 4.16]{DE11} hold:
	\begin{align*}
		a_\rho(\chi_h,\phi_h) &\le C\|\chi_h\|_{\rho}\|\phi_h\|_{\rho} \quad \chi_h,\phi_h\in X_h^0,\\
		a_c(\psi_h,\xi_h) &\le C\|\psi_h\|_{c}\|\xi_h\|_{c} \quad \psi_h,\xi_h\in X_h,\\
		a_u(\bv_h,\bw_h) &\le C\|\bv_h\|_{u}\|\bw_h\|_{u} \quad \bv_h,\bw_h\in \bV_h.
	\end{align*}
\end{lemma}

\noindent
We also state below the discrete inf-sup condition \cite[Theorem 6.8]{Riv08}:
\begin{lemma}\label{lem:infsup}
	There exists a constant $\beta^*>0$, independent of $h$ such that 
	\begin{equation*}
		\inf_{q_h\in M_h}\sup_{\bv_h\in\tilde{\bV}_h} \frac{d(\bv_h,q_h)}{\|\bv_h\|_{u}\|q_h\|} \ge \beta^*,
	\end{equation*}
	where
	\[\tilde{\bV}_h = \left\{\bv_h\in\bV_h: \jump{\bv_h}|_e\cdot \bn_e = 0 \quad\forall~ e \in \Gamma_h\right\}.\]
\end{lemma}

We now state the positivity property and the boundedness properties of trilinear forms in the following lemmas. For proof, we refer to \cite[Lemma 6.39, Lemma 6.40]{DE11}.
\begin{lemma}[Skew-symmetric property] \label{lem:skew}
	Using integration by parts, the  trilinear forms $b_1(\cdot,\cdot,\cdot)$ and $b_2(\cdot,\cdot,\cdot)$ defined in \eqref{tribc}  and \eqref{tribu} satisfy the followings:
	\begin{align}
		b_1(\bv_h,\psi_h,\phi_h) = - b_1(\bv_h,\phi_h,\psi_h) \quad \mbox{and} \quad b_2(\bv_h,\bw_h,\bphi_h) = - b_2(\bv_h,\bphi_h,\bw_h).
	\end{align}
	In particular, if $\psi_h=\phi_h$, then $b_1(\bv_h,\psi_h,\psi_h)= 0$ and if $\bw_h=\bphi_h$, then $b_2(\bv_h,\bw_h,\bw_h)= 0$.
\end{lemma}

\begin{lemma}[Boundedness of trilinear terms] \label{lem:bdd}
	There exists a positive constant $C$ such that the nonlinear terms $b_1(\cdot,\cdot,\cdot)$ and $b_2(\cdot,\cdot,\cdot)$ defined in \eqref{tribc}  and \eqref{tribu} satisfy the followings:
	\begin{align}
		b_1(\bv_h,\psi_h,\phi_h) &\le C \|\bv_h\|_u \|\psi_h\|_\rho \|\phi_h\|_\rho \quad\forall~ \bv_h\in\bV_h, \psi_h,\phi_h\in X_h^0, \\
		b_1(\bv_h,\psi_h,\phi_h) &\le C \|\bv_h\|_u \|\psi_h\|_c \|\phi_h\|_c \quad\forall~ \bv_h\in\bV_h, \psi_h,\phi_h\in X_h, \\
		b_2(\bv_h,\bw_h,\bphi_h) &\le C \|\bv_h\|_u \|\bw_h\|_u \|\bphi_h\|_u \quad\forall~\bv_h, \bw_h, \bphi_h \in \bV_h.
	\end{align}
\end{lemma}

\begin{lemma}[Boundedness of chemotactic term] \label{lem:gbdd}
	There exists a positive constant $C$ such that the chemotactic term $g(\cdot,\cdot,\cdot)$  defined in \eqref{defg} satisfy the following:
	\begin{align}
	g(\chi_h,\psi_h,\phi_h) &\le C \|\chi_h\|_{L^\infty} \|\psi_h\|_c \|\phi_h\|_\rho \quad\forall~ \chi_h,\phi_h\in X^0_h, \psi_h\in X_h.
	\end{align}
\end{lemma}
From the discrete coercivity (Lemma \ref{lem:prop_bilinear}), discrete inf-sup condition (Lemma \ref{lem:infsup}), boundedness of nonlinear and chemotactic terms (Lemmas \ref{lem:bdd}, \ref{lem:gbdd}), one can easily prove the existence and uniqueness of the semidiscrete dG solution of \eqref{smdg1rho}-\eqref{smdg4div} following the similar argument as \cite{KR05}. 

\noindent
Now we state an approximation property of an interpolation operator. For a proof, see \cite[pp. 31, (1.46)]{DE11}. %Pedro
\begin{lemma}[Optimal approximation property] \label{polyapprox}
	For all $E\in\Eh$ and all polynomial degree $r$, there is a linear interpolation operator $\mathcal{I}_h:L^2(E)\to \mathbb{P}_{k}$ such that for all $s\in\{0,\dots,k+1\}$ and for all $\phi\in H^s(E)$ the following relation holds 
	\begin{equation*}
		\|\phi-\mathcal{I}_h\phi\|_{H^m(E)} \le Ch_E^{s-m}\|\phi\|_{H^s(E)} \quad \forall~ m\in \{0,\dots,s\},
	\end{equation*} 
	where $C$ is independent of both $h$ and $E$.
\end{lemma}
\noindent
\textbf{Projection Operator:}
First we introduce the orthogonal $L^2$-projections $\tilde{P}_h:L^2\to X_h^0$, $P_h: L^2\to X_h$ and $\bP_h:\bL^2\to \bV_h$ satisfying for any function $\phi\in L^2$, $\bphi\in\bL^2$
\begin{align*}
	\left(\phi-\tilde{P}_h\phi, \chi_h\right) &= 0 \quad \forall~ \chi_h\in X^0_h,\\
	\left(\phi-P_h\phi, \psi_h\right) &= 0 \quad \forall~ \psi_h\in X_h,\\
	\left(\bphi-\bP_h\bphi,\bv_h\right) &= 0 \quad \forall~ \bv_h\in \bV_h.
\end{align*}
The orthogonal projections satisfy the following approximation properties:
\begin{lemma}\label{lem:ortho}
	There exists a positive constant $C$ independent of $h$ such that the followings hold for all $s\in\{0,1,\dots,k+1\}$:
	\begin{align*}
		\|\phi-\tilde{P}_h\phi\| + h \|\phi-\tilde{P}_h\phi\|_\rho &\le Ch^{s}\|\phi\|_{H^{s}} \quad\forall~\phi\in H^s,\\
		\|\phi-P_h\phi\| + h \|\phi-P_h\phi\|_c &\le Ch^{s}\|\phi\|_{H^{s}} \quad\forall~\phi\in H^s,\\
		\|\bphi-\bP_h\bphi\| + h \|\bphi-\bP_h\bphi\|_u &\le Ch^{s}\|\bphi\|_{\bH^{s}} \quad\forall~\bphi\in \bH^s.
	\end{align*}
\end{lemma}

\noindent
We define broken elliptic projection operator $R_hc\in X_h$ for a given $c$ as
\begin{align}\label{ellipj}
	a_c(c-R_hc, \psi_h) = 0 \quad \forall~ \psi_h \in X_h\quad \text{with}\quad (c-R_hc,1)=0.
\end{align}
One can easily derive the following approximation properties for the above-mentioned operator. For a proof, see \cite[(1.34)]{GRW05spliting}.
\begin{lemma}\label{lem:ritz}
	There exists a positive constant $C$ independent of $h$ such that the followings hold for all $s\in\{0,1,\dots,k+1\}$:
	\begin{align*}
	\|c-R_hc\| + h \|c-R_hc\|_c \le Ch^{s}\|c\|_{H^{s}} \quad\forall~c\in H^s.
	\end{align*}
\end{lemma}
\noindent
We now define a modified elliptic projection operator $Q_h\trho\in X_h^0$ of the solution $\trho$ of \eqref{dg1rho} for a given $c$ and $R_hc$, satisfying
\begin{align}\label{modi:eq}
	\mu a_\rho(\trho-Q_h\trho, \chi_h) -  \beta g(\trho+m_0, c-R_hc,\chi_h) = 0 \quad \forall~ \chi_h \in X_h^0.
\end{align}
The following approximation properties hold for the above-mentioned operator.
\begin{lemma}\label{lem:modi}
	There exists a positive constant $C$ independent of $h$ such that the followings hold for all $s\in\{0,1,\dots,k+1\}$:
	\begin{align*}
		\|\trho-Q_h\trho\| + h \|\trho-Q_h\trho\|_\rho \le Ch^{s}\left(\|\trho\|_{H^{s}}+\|c\|_{H^s}\right) \quad\forall~\trho, c\in H^s.
	\end{align*}
\end{lemma}

\begin{proof}
Let $\xi_\rho=\trho-Q_h\trho$ and $\xi_c=c-R_hc$. Choose $\chi_h=\tilde{P}_h\xi_\rho=\xi_\rho-(\trho-\tilde{P}_h\trho)$ in \eqref{modi:eq} to obtain
\begin{align}
	\mu a_\rho(\xi_\rho,\xi_\rho) = \mu a_\rho(\xi_\rho,\trho-\tilde{P}_h\trho) + \beta g(\trho+m_0, \xi_c,\tilde{P}_h\xi_\rho).
\end{align}
A use of Lemma \ref{lem:prop_bilinear} with $\trho = \rho-m_0$ shows
\begin{align}\label{ag1}
	\mu \|\xi_\rho\|_\rho^2 \le \mu \|\xi_\rho\|_\rho\|\trho-\tilde{P}_h\trho\|_\rho + \beta g(\rho, \xi_c,\tilde{P}_h\xi_\rho).
\end{align}
From the definition of $g(\cdot,\cdot,\cdot)$, it follows that
\begin{align}\label{ag3}
	g(\rho,\xi_c,\tilde{P}_h\xi_\rho) &= \sum_{E\in\Eh}\int_E \rho \nabla\xi_c\cdot\nabla\tilde{P}_h\xi_\rho \dx - \sum_{e\in\Gamma_h^0}\int_e \avr{\rho\nabla\xi_c}\cdot\bn_e \jump{\tilde{P}_h\xi_\rho} \ds - \sum_{e\in\Gamma_h^0}\int_e \avr{\rho\nabla\tilde{P}_h\xi_\rho}\cdot\bn_e \jump{\xi_c} \ds \nonumber\\
	&= A_1+A_2+A_3.
\end{align}
An application of the H{\"o}lder inequality with the Young's inequality yields
\begin{align}\label{ag3a1}
	|A_1| &\le \sum_{E\in\Eh} \|\rho\|_{L^\infty(E)} \|\nabla\xi_c\|_{L^2(E)}\|\nabla\tilde{P}_h\xi_\rho\|_{L^2(E)}  \leq  \|\rho\|_{L^\infty} \|\xi_c\|_c\|\tilde{P}_h\xi_\rho\|_\rho   \le C\|\rho\|_{H^2}^2 \|\xi_c\|_c^2   + \frac{\mu}{8}\|\xi_\rho\|_\rho^2.
\end{align}
We use the H{\"o}lder inequality with the Young's inequality and Lemma \ref{edge2elem} to bound $A_2$ as
\begin{align}\label{ag3a2}
	|A_2| &\le  \sum_{e\in\Gamma_h^0}\|\rho\|_{L^\infty(e)} \|\nabla\xi_c\|_{L^2(e)}\|\jump{\tilde{P}_h\xi_\rho}\|_{L^2(e)}  \nonumber\\
	& \lesssim   \|\rho\|_{L^\infty}\left(\sum_{e\in\Gamma_h^0}|e|\|\nabla\xi_c\|^2_{L^2(e)}\right)^{\frac{1}{2}} \left(\sum_{e\in\Gamma_h^0}\frac{\sigma_\rho}{|e|}\|\jump{\tilde{P}_h\xi_\rho}\|_{L^2(e)}^2\right)^{\frac{1}{2}}  \nonumber\\
	& \le C\|\rho\|_{H^2}^2 \left(\sum_{E\in\Eh}( \|\nabla\xi_c\|^2_{L^2(E)} + h_E^2\|\nabla^2\xi_c\|_{L^2(E)}^2)\right)   + \frac{\mu}{16}\|\xi_\rho\|_\rho^2.
\end{align}
To estimate $\|\nabla^2\xi_c\|_{L^2(E)}^2$, we apply the triangle inequality and the inverse inequality with Lemma \ref{polyapprox} to derive
\begin{align}\label{ag4}
	\|\nabla^2\xi_c\|_{L^2(E)}^2 &= \|\nabla^2(c-R_hc)\|_{L^2(E)}^2 \nonumber\\
	& \le \|\nabla^2(c-\mathcal{I}_hc)\|_{L^2(E)}^2 + \|\nabla^2(\mathcal{I}_hc-R_hc)\|_{L^2(E)}^2 \nonumber\\
	& \le \|\nabla^2(c-\mathcal{I}_hc)\|_{L^2(E)}^2 + Ch_E^{-2}\|\nabla(\mathcal{I}_hc-R_hc)\|_{L^2(E)}^2 \nonumber\\
	& \le \|\nabla^2(c-\mathcal{I}_hc)\|_{L^2(E)}^2 + Ch_E^{-2}\left(\|\nabla(c-\mathcal{I}_hc)\|_{L^2(E)}^2+\|\nabla(c-R_hc)\|_{L^2(E)}^2\right) \nonumber\\
	&\lesssim h_E^{2k-2}\|c\|_{H^{k+1}(E)}^2.
\end{align}
With the help of the H{\"o}lder inequality, the Young's inequality and Lemma \ref{dised2el}, one can bound $A_3$ as follows:
\begin{align}\label{ag3a3}
	|A_3| &\le  \sum_{e\in\Gamma_h^0}\|\rho\|_{L^\infty(e)} \|\nabla\tilde{P}_h\xi_\rho\|_{L^2(e)}\|\jump{\xi_c}\|_{L^2(e)} \nonumber\\
	& \lesssim  \|\rho\|_{L^\infty} \left(\sum_{e\in\Gamma_h^0}\frac{\sigma_\rho}{|e|}\|\jump{\xi_c}\|_{L^2(e)}^2\right)^{\frac{1}{2}}\left(\sum_{e\in\Gamma_h^0}|e|\|\nabla\tilde{P}_h\xi_\rho\|^2_{L^2(e)}\right)^{\frac{1}{2}}\nonumber\\
	& \le  C \|\rho\|_{H^2}^2 \|\xi_c\|^2_c  + \frac{\mu}{16}\|\xi_\rho\|_\rho^2.
\end{align}
A use of \eqref{ag3} with \eqref{ag4} and the Lemma \ref{lem:ortho}, \ref{lem:ritz} in \eqref{ag1} yields
\begin{align} \label{nablarho}
	\mu\|\xi_\rho\|_\rho^2 \le Ch^{2k}\left(\|\trho\|_{H^{k+1}}^2+\|c\|_{H^{k+1}}^2\right).
\end{align}
To find the $L^2$ estimate, we consider the following dual problem: Consider $w\in H^2\cap H^1/\R$  as a solution of 
\begin{equation}\label{dualprob}
	-\mu\Delta w = \xi_\rho \quad \mbox{in}~ \Omega,
\end{equation}
which satisfy the following regularity result
\begin{equation}\label{dualest}
	\mu\|w\|_{H^2} \lesssim \|\xi_\rho\|.
\end{equation}
Now, multiply \eqref{dualprob} by $\xi_\rho$ and take integration over $\Omega$ to obtain
\begin{equation}
	\|\xi_\rho\|^2 = \mu(-\Delta w, \xi_\rho) = \mu a_\rho(w,\xi_\rho).
\end{equation}
Take $\chi_h=\tilde{P}_hw$ in \eqref{modi:eq} and combine with the above equation to find
\begin{equation}\label{dual2}
	\|\xi_\rho\|^2 = \mu a_\rho(w-\tilde{P}_hw,\xi_\rho)+\beta g(\trho+m_0,\xi_c,\tilde{P}_hw). 
\end{equation}
An application of the Lemma \ref{lem:ortho} helps to bound the first term on the right-hand side of \eqref{dual2} as
\begin{equation}\label{dual22}
	a_\rho(w-\tilde{P}_hw,\xi_\rho) 
	\le C\|w-\tilde{P}_hw\|_\rho\|\xi_\rho\|_\rho 
	\le C h\|w\|_{H^2}\|\xi_\rho\|_\rho.
\end{equation}
We rewrite the last term of \eqref{dual2} in two parts and use $\trho+m_0=\rho$ which implies
\begin{equation}\label{dual222}
	g(\rho,\xi_c,\tilde{P}_hw) = g(\rho,\xi_c,\tilde{P}_hw-w) + g(\rho,\xi_c,w).
\end{equation}
First, we proceed similar way of \eqref{ag3}-\eqref{ag3a3} by replacing $\tilde{P}_h\xi_\rho$ by $\tilde{P}_hw-w$ and without using the Young's inequality. Then, we use \eqref{ag4} and Lemma \ref{lem:ritz} to find
\begin{align}\label{dual3}
g(\rho,\xi_c,\tilde{P}_hw-w) 
& \lesssim  \|\rho\|_{L^\infty} \|\xi_c\|_c\|\tilde{P}_hw-w\|_\rho + \|\rho\|_{L^\infty}\left(\sum_{e\in\Gamma_h^0}|e|\|\nabla\xi_c\|^2_{L^2(e)}\right)^{\frac{1}{2}} \left(\sum_{e\in\Gamma_h^0}\frac{\sigma_\rho}{|e|}\|\jump{\tilde{P}_hw-w}\|_{L^2(e)}^2\right)^{\frac{1}{2}} \nonumber\\
& \qquad + \|\rho\|_{L^\infty}\left(\sum_{e\in\Gamma_h^0}|e|\|\nabla(\tilde{P}_hw-w)\|^2_{L^2(e)}\right)^{\frac{1}{2}} \left(\sum_{e\in\Gamma_h^0}\frac{\sigma_c}{|e|}\|\jump{\xi_c}\|_{L^2(e)}^2\right)^{\frac{1}{2}}\nonumber\\
	& \lesssim  h\|\rho\|_{H^2} \|\xi_c\|_c\|w\|_{H^2} + \|\rho\|_{H^2}^2 \left(\sum_{E\in\Eh}( \|\nabla\xi_c\|^2_{L^2(E)} + h_E^2\|\nabla^2\xi_c\|_{L^2(E)}^2)\right)^{\frac{1}{2}}\|\tilde{P}_hw-w\|_\rho\nonumber\\
	&\qquad+ \|\rho\|_{H^2}^2 \left(\sum_{E\in\Eh}( \|\nabla(\tilde{P}_hw-w)\|^2_{L^2(E)} + h_E^2\|\nabla^2(\tilde{P}_hw-w)\|_{L^2(E)}^2)\right)^{\frac{1}{2}}\|\xi_c\|_c\nonumber\\
	& \lesssim h^{k+1}\|\rho\|_{H^2}\|w\|_{H^2}\|c\|_{H^{k+1}}.
\end{align}
For the second term of \eqref{dual222}, a use of integration by parts yields
\begin{align}\label{dual4}
	g(\trho,\xi_c,w) 
	&= \sum_{E\in\Eh}\int_E \rho \nabla\xi_c\cdot\nabla w \dx - \sum_{e\in\Gamma_h^0}\int_e \left(\avr{\rho\nabla\xi_c}\cdot\bn_e\jump{w} + \avr{\rho\nabla w}\cdot\bn_e\jump{\xi_c} \right) \ds \nonumber\\
	& = -\sum_{E\in\Eh}\int_E \left(\xi_c\nabla\rho \cdot\nabla w + \xi_c\rho\Delta w\right) \dx \nonumber\\
	& \lesssim h^{k+1}\|\rho\|_{H^2}\|w\|_{H^2}\|c\|_{H^{k+1}}.
\end{align}
Combining all the estimates from \eqref{dual22}-\eqref{dual4} and using \eqref{nablarho}, we finally conclude that
\begin{equation}
	\|\xi_\rho\|^2\lesssim h^{k+1}\|\rho\|_{H^2}\|w\|_{H^2}\left(\|\trho\|_{H^{k+1}} + \|c\|_{H^{k+1}}\right).
\end{equation}
Now, a use of \eqref{dualest} completes the rest of the proof.
\end{proof}

\noindent
We also define the modified Stokes operator  $(\bS_h\bu,\bS_h p)\in \bV_h \times M_h$ for the weak solution $(\bu(t), p(t))$ of \eqref{dg3u}-\eqref{dg4div}  satisfying

\begin{equation}\label{stokesproj}
	\left.
	\begin{aligned}
		&\nu a_u(\bu-\bS_h\bu, \bv_h) - d(\bv_h, p-\bS_hp) = 0 \quad \forall~ \bv_h \in \bV_h, \\
		& ~d(\bu-\bS_h\bu, q_h) = 0 \quad \forall~ q_h\in M_h.
	\end{aligned}
	\right\}
\end{equation}
One can easily derive the following approximation properties for the above-mentioned operator. For a proof, see \cite{GR79}.
\begin{lemma}\label{lem:stokes}
	There exists a positive constant $C$ independent of $h$ such that the followings hold for all $s\in\{1,2,\dots,k+1\}$:
	\begin{align*}
		&\|\bu-\bS_h\bu\| + h \left(\|\bu-\bS_h\bu\|_u + \|p-\bS_hp\|\right) \le Ch^{s}\left(\|\bu\|_{\bH^{s}} + \|p\|_{H^{s-1}}\right).
	\end{align*}
\end{lemma}

\subsection{Fully Discrete dG Formulation}

For time discretization, we take a uniform partition of the time interval $[0,T]$ as $0=t_0<t_1<\cdots<t_n=T$ with time step $\Delta t = t_i-t_{i-1}, 1\le i\le n$ and $t_n = n\Delta t$. For a finite sequence $\{\phi^n\}_{1\le n <\infty}$, we define the backward difference quotient $\partial_t\phi^n = \frac{1}{\Delta t}(\phi^n - \phi^{n-1})$. For a continuous function $v(t)$, we set $v^n=v(t_n)$. 
The fully discrete semi-implicit dG formulation of \eqref{dg1rho}-\eqref{dg4div} reads as:

\noindent
\textbf{Step I:} Seek $(\trho_h^n, c_h^n, \bu_h^n, p_h^n)\in X_h^0 \times X_h \times \bV_h \times M_h$ such that for all $n>0$
\begin{align}
	&(\partial_t\trho_{h}^n,\chi_h) + \mu a_\rho(\trho_h^n,\chi_h) + b_1(\bu_h^{n},\trho_h^n,\chi_h) - \beta g(\trho_h^{n-1}+m_0, c_h^{n}, \chi_h)= 0 \quad\forall~\chi_h\in X_h^0, \label{fddg1rho}\\
	&(\partial_t c_{h}^{n},\psi_h) + \kappa a_c(c_h^n,\psi_h) + b_1(\bu_h^{n}, c_h^{n}, \psi_h) + \gamma((\trho_h^{n-1}+m_0)c_h^{n},\psi_h) = 0 \quad\forall~\psi_h\in X_h, \label{fddg2fc}\\
	&(\partial_t\bu_{h}^{n},\bv_h) + \nu a_u(\bu_h^{n},\bv_h) + b_2(\bu_h^{n},\bu_h^{n},\bv_h) - d(\bv_h,p_h^{n}) = ((\trho_h^{n-1}+m_0)\nabla\Phi^{n},\bv_h)  \quad\forall~\bv_h\in\bV_h, \label{fddg3u}\\
	&d(\bu_h^{n},q_h)=0 \quad\forall~ q_h\in M_h, \label{fddg4div}
\end{align}
with $(\trho_h^0, c_h^0, \bu_h^0)=(Q_h\trho_0, R_hc_0, S_h\bu_0)$.\\[3pt]
\textbf{Step II:} Set $\rho_h^n=\trho_h^n+m_0$.

\vspace{1mm}

It is clear that, \eqref{fddg1rho}-\eqref{fddg4div} is a finite dimensional algebraic system. Using Lemmas \ref{lem:prop_bilinear}-\ref{lem:gbdd}, one can easily show that there exists a unique solution $(\trho_h^n, c_h^n, \bu_h^n, p_h^n)\in X_h^0 \times X_h \times \bV_h \times M_h$  satisfying \eqref{fddg1rho}-\eqref{fddg4div}.

\noindent 
Now, one can easily prove the following qualitative property. 
\begin{lemma}\label{lem:mass}
	The fully discrete solution $\rho_h^n$ satisfies the following mass conservation property:
	\begin{equation}
		\sum_{E\in \Eh}\int_E \rho_h^n(x) \dx = \sum_{E\in \Eh}\int_E \rho_0(x) \dx = m_0|\Omega|.
	\end{equation}
\end{lemma}

\noindent
For our subsequence analysis, we will use the following standard discrete versions of the trace inequalities (\cite[Lemma 1.46]{DE11}), the Ladyzhenskaya's inequality \cite{Tem84}, and the Agmon's inequality \cite[Lemma 2.2]{KSS09}.
\begin{lemma}\label{dised2el}
	For each element $E\in\Eh$ with diameter $h_E$, there exists a constant $C$ independent of $h_E$ such that the followings hold: For all $\phi_h\in X_h ~ \mbox{or}~ X_h^0 ~\mbox{or}~\bV_h$
	\begin{align*}
		\|\phi_h\|_{L^2(e)}  \le C h_E^{-\frac{1}{2}}\|\phi_h\|_{L^2(E)}, \quad \text{and}\quad 
		\|\nabla\phi_h\|_{L^2(e)}  \le C h_E^{-\frac{1}{2}}\|\nabla\phi_h\|_{L^2(E)},
	\end{align*}
	where $e$ is an edge of the element $E$. The non-Hilbertian version of the discrete trace inequality is as follows \cite[Lemma 1.52]{DE11}: For $1\le p \le \infty$,
	\begin{align*}
		\|\phi_h\|_{L^p(e)} & \le C h_E^{-\frac{1}{p}}\|\phi_h\|_{L^p(E)}.
	\end{align*}
\end{lemma}
\begin{lemma}\label{disl42l2}
	For each edge $E\in\Eh$, there exists a constant $C>0$, independent of $h_E$, such that the followings hold
	\begin{align*}%\label{lady}
		\|\phi_h\|_{L^4(E)} &\le C \|\phi_h\|_{L^2(E)}^{\frac{1}{2}}\|\nabla\phi_h\|_{L^2(E)}^{\frac{1}{2}},\\
		\|\nabla\phi_h\|_{L^4(E)} &\le C \|\nabla\phi_h\|_{L^2(E)}^{\frac{1}{2}}\|\nabla^2\phi_h\|_{L^2(E)}^{\frac{1}{2}},\\
		\|\phi_h\|_{L^\infty(E)} &\le C \|\phi_h\|_{L^2(E)}^{\frac{1}{2}}\|\nabla^2\phi_h\|_{L^2(E)}^{\frac{1}{2}},
	\end{align*}
	for all $\phi_h\in X_h ~ \mbox{or}~ X_h^0 ~\mbox{or}~\bV_h$.
\end{lemma}

\section{A Priori Error Analysis}

In this section, we will analyze the error due to space and temporal discretization.
We define any continuous in time function $\phi(t)$ at a time level $t_n$ as $\phi^n = \phi(t_n)$. We denoted the errors as
\begin{align}\label{errorsplit}
	\e_u^n = \bu^n - \bu_h^n, \quad
	e_\rho^n = \trho^n - \trho_h^n, \quad
	e_c^n = c^n - c_h^n.
\end{align}

\noindent
From \eqref{dg1rho}-\eqref{dg4div} and \eqref{fddg1rho}-\eqref{fddg4div}, the error equations are written  as:
\begin{align}
	(\pt e_\rho^n,\chi_h) + \mu a_\rho(e_\rho^n,\chi_h) & = E_\rho(\chi_h) + \Lambda_\rho(\chi_h) + \beta F_\rho(\chi_h) \quad \forall~\chi_h\in X_h^0, \label{erreqrho} \\
	(\pt e_c^n,\psi_h) + \kappa~a_c(e_c^n,\psi_h) + \gamma m_0(e_c^n,\psi_h) & = E_c(\psi_h) + \Lambda_c(\psi_h) + \gamma F_c(\psi_h) \quad \forall~\psi_h\in X_h, 	\label{erreqc} \\
	(\pt\e_u^n, \bv_h) + \nu~a_u(\e_u^n,\bv_h) - d(\bv_h,e_p^n) + d(\e_u^n,q_h) & = E_u(\bv_h) + \Lambda_u(\bv_h) + F_u(\bv_h) \quad \forall~\bv_h\in \bV_h, q_h\in M_h, \label{erreqvel}
\end{align}
where 
\begin{align}
	E_\rho(\chi_h) &= \left( \pt \trho^n - \trho_t^n, \chi_h \right) = -\frac{1}{k}\int_{t_{n-1}}^{t_n} (s-t_{n-1})\left(\trho_{tt}(t), \chi_h\right)\dt, \label{Erho}\\
	\Lambda_\rho(\chi_h)& = b_1(\bu_h^{n},\trho_h^n,\chi_h) - b_1(\bu^n,\trho^n,\chi_h) = b_1(\e_u^n,e_\rho^n,\chi_h)- b_1(\bu^n,e_\rho^n,\chi_h) - b_1(\e_u^n,\trho^n,\chi_h), \label{Lrho}\\
	F_\rho(\chi_h) &=  g(\trho^{n},c^n,\chi_h)-g(\trho_{h}^{n-1},c_h^n,\chi_h) +   g(m_0,e_c^n,\chi_h), \label{Frho}
\end{align}
and
\begin{align}
	E_c(\psi_h) &= \left( \pt c^n - c_t^n, \psi_h \right) = -\frac{1}{k}\int_{t_{n-1}}^{t_n} (s-t_{n-1})\left(c_{tt}(t), \psi_h\right)\dt, \label{Ec}\\
	\Lambda_c(\psi_h)& = b_1(\bu_h^{n},c_h^n,\psi_h) - b_1(\bu^n,c^n,\psi_h) = b_1(\e_u^n,e_c^n,\psi_h)- b_1(\bu^n,e_c^n,\psi_h) - b_1(\e_u^n,c^n,\psi_h), \label{Lc}\\
	F_c(\psi_h) &=  (\trho_{h}^{n-1}c_h^n,\psi_h)-(\trho^{n}c^n,\psi_h) \nonumber\\
	& = (e_\rho^{n-1}e_c^n,\psi_h) - (\trho^{n-1}e_c^n,\psi_h) - (e_\rho^{n-1}c^n,\psi_h)- ((\trho^n-\trho^{n-1})c^n,\psi_h), \label{Fc}
\end{align}
and
\begin{align}
	E_u(\bv_h) &= \left( \pt\bu^n - \bu_t^n, \bv_h \right) = -\frac{1}{k}\int_{t_{n-1}}^{t_n} (s-t_{n-1})\left(\bu_{tt}(t), \bv_h\right)\dt, \label{Eu}\\
	\Lambda_u(\bv_h)& = b_2(\bu_h^{n},\bu_h^n,\bv_h) - b_2(\bu^n,\bu^n,\bv_h) = b_2(\e_u^n,\e_u^n,\bv_h)- b_2(\bu^n,\e_u^n,\bv_h) - b_2(\e_u^n,\bu^n,\bv_h), \label{Lu}\\
	F_u(\bv_h) &= ((\trho_h^{n-1}+m_0)\nabla\Phi^{n},\bv_h) - ((\trho^{n}+m_0)\nabla\Phi^{n},\bv_h) \nonumber\\ &=-(e_{\rho}^{n-1}\nabla\Phi^n,\bv_h) - \left(\left(\trho^n-\trho^{n-1}\right)\nabla\Phi^n,\bv_h\right).  \label{Fu}
\end{align}
Using the projection operators, we split the errors as
\begin{equation}\label{usplit}
	\begin{aligned}
		\e_u^n = \bu^n - \bu_h^n = \underbrace{\bu^n - \bS_h\bu^n}_{\bxi_u^n} + \underbrace{\bS_h\bu^n - \bu_h^n}_{\btheta_u^n},
		\quad  \quad
		e_p^n = p^n - p_h^n = \underbrace{p^n - \bS_hp^n}_{\xi_p^n} + \underbrace{\bS_hp^n - p_h^n}_{\theta_p^n},\\
		e_c^n = c^n - c_h^n = \underbrace{c^n - R_hc^n}_{\xi_c^n} + \underbrace{R_hc^n - c_h^n}_{\theta_c^n},
		\quad  \quad
		e_\rho^n = \trho^n - \trho_h^n = \underbrace{\trho^n - Q_h\trho^n}_{\xi_\rho^n} + \underbrace{Q_h\trho^n - \trho_h^n}_{\theta_\rho^n}.
	\end{aligned}
\end{equation}
Since the estimates of $\xi_\rho, \xi_c, \bxi_u^{n}$ and $\xi_p$ are known from the Lemmas \ref{lem:ritz}, \ref{lem:modi} and \ref{lem:stokes}, it is enough to estimate $\theta_\rho, \theta_c, \btheta_u^{n}$ and $\theta_p$.
A use of the elliptic projection \eqref{ellipj}, the modified elliptic projection \eqref{modi:eq} and the Stokes projection \eqref{stokesproj} with \eqref{usplit} help us to find the equations for $\theta_\rho, \theta_c, \btheta_u$ and $\theta_p$ as follow: For all $(\chi_h, \psi_h, \bv_h, q_h)\in X_h^0\times X_h \times \bV_h\times M_h$
\begin{align}
	\label{erreqthetarho}
	&(\pt\theta_\rho^n, \chi_h) + \mu a_\rho(\theta_\rho^n,\chi_h) = -(\pt\xi_\rho^n, \chi_h) + E_\rho(\chi_h) + \Lambda_\rho(\chi_h) + \beta G_\rho(\chi_h),\\	\label{erreqthetac}
	&(\pt\theta_c^n, \psi_h) + \kappa~a_c(\theta_c^n,\psi_h) + \gamma m_0(\theta_c^n,\psi_h) = -(\pt\xi_c^n, \psi_h) - \gamma m_0(\xi_c^n,\psi_h) + E_c(\psi_h) + \Lambda_c(\psi_h) + F_c(\psi_h), \\
	\label{erreqtheta}
	&(\pt\btheta_u^n, \bv_h) + \nu~a_u(\btheta_u^n,\bv_h) - d(\bv_h,\theta_p^n) + d(\btheta_u^n,q_h) = -(\pt\bxi_u^n, \bv_h)  + E_u(\bv_h) + \Lambda_u(\bv_h) + F_u(\bv_h),
\end{align}
where 
\begin{align}
	G_\rho(\chi_h) := F_\rho(\chi_h) - g(\trho^{n}+m_0,\xi_c^n,\chi_h) = &   g(e_\rho^{n-1},R_hc^n,\chi_h) + g(\trho^{n-1},\theta_c^n,\chi_h) - g(e_\rho^{n-1},\theta_c^n,\chi_h) \nonumber \\ 
	&  + g(\trho^n-\trho^{n-1},R_hc^n,\chi_h) +   g(m_0,\theta_c^n,\chi_h).
\end{align}\label{Grho}

%%%%%%%%%%%%%%%%%%%%%%%%%%%%%%%%%%%%%%%%%%%%%%%%%%%% L2 Error %%%%%%%%%%%%%%%%%%%%%%%%%%%%%%

\subsection{$L^2$-Error Estimates}
In this section, we obtain optimal $L^2$-error estimates for the fluid velocity, the concentration and the cell density. First, we prove a couple of lemmas to conclude the final result.

\begin{lemma}\label{lem:thetau}
	Under the regularity assumption (\textbf{A1}), the following holds true
	\begin{equation}\label{utheta16}
		\|\btheta_u^N\|^2 + \nu\Delta t\sum_{n=1}^N\|\btheta_u^n\|_u^2 
		\lesssim \left(h^{2k+2} + (\Delta t)^2\right)  + \Delta t \sum_{n=1}^N \left(\|\theta_\rho^{n-1}\|^2+\|\btheta_u^n\|^2\right).
	\end{equation}
\end{lemma}
\begin{proof}
	Choose $\bv_h = \btheta_u^n$ and $q_h = \theta_p^n$ in \eqref{erreqtheta} and use the fact $(\pt\btheta_u^n,\btheta_u^n)\ge \frac{1}{2}\pt\|\btheta_u^n\|^2$ to obtain
	\begin{equation}\label{utheta1}
		\frac{1}{2}\pt\|\btheta_u^n\|^2 + \nu\|\btheta_u^n\|_u^2 \le  -(\pt\bxi_u^n, \btheta_u^n)  + E_u(\btheta_u^n) + \Lambda_u(\btheta_u^n) + F_u(\btheta_u^n).
	\end{equation}
	First term on the right hand side of \eqref{utheta1} can be bounded as
	\begin{align}\label{utheta2}
		|(\pt\bxi_u^n, \btheta_u^n)| 
		\le \frac{1}{\Delta t}\int_{t_{n-1}}^{t_n} |(\bxi_{ut}(t),\btheta_u^n)| \dt
		\le \frac{C}{\Delta t}\int_{t_{n-1}}^{t_n} \|\bxi_{ut}(t)\|^2 \dt + \frac{\nu}{16}\|\btheta_u^n\|_u^2.
	\end{align}
	A use of the Cauchy-Schwarz inequality and the Young's inequality with \eqref{Eu} helps to bound the second term on the right hand side of \eqref{utheta1} as
	\begin{align}\label{utheta3}
		|E_u(\btheta_u^n)| \le \left(\frac{1}{\Delta t} \int_{t_{n-1}}^{t_n} (s-t_{n-1})\|\bu_{tt}(t)\|\dt\right) \|\btheta_u^n\| 
		\le C \Delta t \int_{t_{n-1}}^{t_n} \|\bu_{tt}(t)\|^2\dt + \frac{\nu}{16}\|\btheta_u^n\|_u^2.
	\end{align}
	Using \eqref{usplit} in \eqref{Lu}, we rewrite the third term on the right hand side of \eqref{utheta1} as
	\begin{align}
		\Lambda_u(\btheta_u^n) = & b_2(\bxi_u^n,\bxi_u^n,\btheta_u^n) + b_2(\bxi_u^n,\btheta_u^n,\btheta_u^n) + b_2(\btheta_u^n,\bxi_u^n,\btheta_u^n) + b_2(\btheta_u^n,\btheta_u^n,\btheta_u^n) \nonumber\\
		& - 
		b_2(\bu^n,\bxi_u^n,\btheta_u^n) - b_2(\bu^n,\btheta_u^n,\btheta_u^n) -
		b_2(\bxi_u^n,\bu^n,\btheta_u^n) - b_2(\btheta_u^n,\bu^n,\btheta_u^n). \label{Lutheta}
	\end{align}
	The second, fourth, and sixth terms of \eqref{Lutheta} are zero due to Lemma \ref{lem:skew}. For the first term, we use Lemma \ref{lem:bdd} to bound
	\begin{align}\label{utheta4}
		|b_2(\bxi_u^n,\bxi_u^n,\btheta_u^n)| \le C\|\bxi_u^n\|_u^2\|\btheta_u^n\|_u \le  C\|\bxi_u^n\|_u^4 + \frac{\nu}{16}\|\btheta_u^n\|_u^2.
	\end{align}
	From the definition of the trilinear form \eqref{tribu}, an application of the H{\"o}lder inequality yields
	\begin{align}\label{utheta5}
		|b_2(\btheta_u^n,\bxi_u^n,\btheta_u^n)| 
		& \le  \sum_{E\in\Eh} \|\nabla\bxi_u^n\|_{\bL^2(E)}  \|\btheta_u^n\|_{\bL^4(E)}^2 + \frac{1}{2}\sum_{E\in\Eh}\|\nabla\cdot\btheta_u^n\|_{\bL^2(E)}\|\bxi_u^n\|_{\bL^4(E)}\|\btheta_u^n\|_{\bL^4(E)} \nonumber\\
		&\quad+ \sum_{e\in\Gamma_h}\|\btheta_u^n\|_{\bL^4(e)}^2 \|\jump{\bxi_u^n}\|_{\bL^2(e)} + \frac{1}{2} \sum_{e\in\Gamma_h}\|\jump{\btheta_u^n}\|_{\bL^2(e)} \|\bxi_u^n\|_{\bL^2(e)}\|\btheta_u^n\|_{\bL^\infty(e)} \ds \nonumber\\
		& = I_1 + I_2 + I_3 + I_4.
	\end{align}\label{utheta6}
	A use of the Gagliardo-Nirenberg inequality \eqref{GNI} with Lemma \ref{disl42l2} and the Young's inequality shows
	\begin{align}
		I_1+I_2 &\lesssim  \sum_{E\in\Eh}\left( \|\nabla\bxi_u^n\|_{\bL^2(E)}  \|\btheta_u^n\|_{\bL^2(E)} \|\nabla\btheta_u^n\|_{\bL^2(E)} + \|\bxi_u^n\|_{\bL^2(E)}^{\frac{1}{2}}\|\nabla\bxi_u^n\|_{\bL^2(E)}^{\frac{1}{2}} \|\btheta_u^n\|_{\bL^2(E)}^{\frac{1}{2}}\|\nabla\btheta_u^n\|_{\bL^2(E)}^{\frac{3}{2}}\right) \nonumber\\
		& \le C \sum_{E\in\Eh}\left( \|\nabla\bxi_u^n\|_{\bL^2(E)}^2  + \|\bxi_u^n\|_{\bL^2(E)}^2\|\nabla\bxi_u^n\|_{\bL^2(E)}^2 \right)\|\btheta_u^n\|_{\bL^2(E)}^2 + \frac{\nu}{32}\sum_{E\in\Eh}\|\nabla\btheta_u^n\|_{\bL^2(E)}^2\nonumber\\
		&\le C \left(1 + \|\bxi_u^n\|^2\right)\|\bxi_u^n\|_u^2\|\btheta_u^n\|^2 + \frac{\nu}{32}  \|\btheta_u^n\|_u^2.
	\end{align}
	We apply the Young's inequality and Lemmas \ref{edge2elem}, \ref{dised2el} and \ref{disl42l2} with \eqref{inv.hypo} to derive
	\begin{align}\label{utheta7}
		I_3+I_4 &\lesssim \left(\sum_{e\in\Gamma_h} |e|\|\btheta_u^n\|_{\bL^4(e)}^4\right)^{\frac{1}{2}}   \left(\sum_{e\in\Gamma_h} \frac{\sigma_u}{|e|}\|\jump{\bxi_u^n}\|_{\bL^2(e)}^2\right)^{\frac{1}{2}} \nonumber\\
		&\qquad+ \left(\sum_{e\in\Gamma_h} |e|\|\bxi_u^n\|_{\bL^2(e)}^2\|\btheta_u^n\|_{\bL^\infty(e)}^2\right)^{\frac{1}{2}}   \left(\sum_{e\in\Gamma_h} \frac{\sigma_u}{|e|}\|\jump{\btheta_u^n}\|_{\bL^2(e)}^2\right)^{\frac{1}{2}} \nonumber\\
		& \lesssim 
		\left(\sum_{E\in\Eh}  h_E^{-1}\|\btheta_u^n\|_{\bL^2(E)}^2\|\nabla\btheta_u^n\|_{\bL^2(E)}^2\right)^{\frac{1}{2}}   \|\bxi_u^n\|_u \nonumber\\
		&\qquad+ \left(\sum_{E\in\Eh} \left( \|\bxi_u^n\|_{\bL^2(E)}^2+h_E^{2}\|\nabla\bxi_u^n\|_{\bL^2(E)}^2\right) h_E^{-3}\|\btheta_u^n\|_{\bL^2(E)}^2\right)^{\frac{1}{2}}   \|\btheta_u^n\|_u \nonumber\\
		& \le C\left(h^{-3} \|\bxi_u^n\|^2 + h^{-1}\|\bxi_u^n\|_u^2\right)\|\btheta_u^n\|^2 + \frac{\nu}{32} \|\btheta_u^n\|_u^2.
	\end{align}
	Using Lemma \ref{lem:skew} with \eqref{tribu}, we rewrite the fifth term of \eqref{Lutheta} as
	\begin{align}
		b_2(\bu^n,\bxi_u^n,\btheta_u^n) =  -b_2(\bu^n,\btheta_u^n,\bxi_u^n)   = - \sum_{E\in\Eh}\int_E (\bu^n\cdot\nabla)\btheta_u^n \cdot \bxi_u^n \dx - \frac{1}{2}\sum_{E\in\Eh}\int_E (\nabla\cdot\bu^n)\btheta_u^n\cdot\bxi_u^n \dx \nonumber\\
		+ \sum_{e\in\Gamma_h}\int_e \avr{\bu^n} \cdot\bn_e \jump{\btheta_u^n}\cdot\avr{\bxi_u^n}\ds + \frac{1}{2} \sum_{e\in\Gamma_h}\int_e \jump{\bu^n}\cdot\bn_e \avr{\bxi_u^n\cdot\btheta_u^n} \ds.
	\end{align}
	The last term on the above equality becomes zero due to $\jump{\bu^n}\cdot\bn_e=0$. The other terms can be bounded by applying the H{\"o}lder inequality, Young's inequality with the Gagliardo-Nirenberg inequality \eqref{GNI}, the Agmon's inequality \eqref{AI} and Lemmas \ref{disl42l2}, \ref{edge2elem}, \ref{dised2el} by
	\begin{equation}
		|b_2(\bu^n,\bxi_u^n,\btheta_u^n)| \le C \left(
		\|\bxi_u^n\|^2 + h^2\|\bxi_u^n\|_u^2\right)\|\bu^n\|_{H^2}^2 + \frac{\nu}{16}\|\btheta_u^n\|_u^2.
	\end{equation}
	From \eqref{tribu}, we write the seventh term of \eqref{Lutheta} as below:
	\begin{align}\label{utheta10}
		b_2(\bxi_u^n,\bu^n,\btheta_u^n) & = \sum_{E\in\Eh}\int_E (\bxi_u^n\cdot\nabla)\bu^n \cdot \btheta_u^n \dx + \frac{1}{2}\sum_{E\in\Eh}\int_E (\nabla\cdot\bxi_u^n)\bu^n\cdot\btheta_u^n \dx  \nonumber\\
		&\quad- \sum_{e\in\Gamma_h}\int_e \avr{\bxi_u^n} \cdot\bn_e \jump{\bu^n}\cdot\avr{\btheta_u^n}\ds - \frac{1}{2} \sum_{e\in\Gamma_h}\int_e \jump{\bxi_u^n}\cdot\bn_e \avr{\bu^n\cdot\btheta_u^n} \ds.
	\end{align}
	An application of the integration by parts in the second term of the above equation shows
	\begin{align}\label{utheta11}
		\frac{1}{2}\sum_{E\in\Eh}\int_E (\nabla\cdot\bxi_u^n)\bu^n\cdot\btheta_u^n \dx &= - \frac{1}{2}\sum_{E\in\Eh}\int_E (\bxi_u^n\cdot\nabla)\bu^n\cdot\btheta_u^n \dx - \frac{1}{2}\sum_{E\in\Eh}\int_E (\bxi_u^n\cdot\nabla)\btheta_u^n\cdot\bu^n \dx \nonumber\\
		&\quad + \frac{1}{2}\sum_{E\in\Eh}\int_E (\bxi_u^n\cdot\bn_E)\bu^n\cdot\btheta_u^n \dx.
	\end{align}
	One can transfer the last sum of \eqref{utheta11} from each element to each edge. Let $e$ be the common edge of two elements $E_l$ and $E_r$ with normal $\bn_l$ and $\bn_r$, respectively. Then, thanks to $[ab]=[a]\{b\}+\{a\}[b]$
	\begin{align} \label{utheta12}
		\frac{1}{2}\sum_{E\in\Eh}\int_E  (\bxi_u^n\cdot\bn_E)&\bu^n\cdot\btheta_u^n \dx = \frac{1}{2}\sum_{e\in\Gamma_h} \left(\int_e (\bxi_u^n|_{E_l}\cdot\bn_l)\bu^n|_{E_l}\cdot\btheta_u^n|_{E_l} \dx + \int_e (\bxi_u^n|_{E_r}\cdot\bn_r)\bu^n|_{E_r}\cdot\btheta_u^n|_{E_r} \dx \right)\nonumber\\
		=& \frac{1}{2}\sum_{e\in\Gamma_h}\int_e [(\bxi_u^n\cdot\bn_e)\bu^n\cdot\btheta_u^n] \dx\nonumber\\
		=& \frac{1}{2}\sum_{e\in\Gamma_h}\int_e \left( \jump{\bxi_u^n}\cdot\bn_e\avr{\bu^n\cdot\btheta_u^n}  +  \avr{\bxi_u^n}\cdot\bn_e\jump{\bu^n}\cdot\avr{\btheta_u^n}   +  \avr{\bxi_u^n}\cdot\bn_e\avr{\bu^n}\cdot\jump{\btheta_u^n}\right) \dx.
	\end{align}
	Combining the estimates \eqref{utheta10}-\eqref{utheta12} follows
	\begin{align}\label{utheta13}
		b_2(\bxi_u^n,\bu^n,\btheta_u^n) & = \frac{1}{2}\sum_{E\in\Eh}\int_E (\bxi_u^n\cdot\nabla)\bu^n \cdot \btheta_u^n \dx - \frac{1}{2}\sum_{E\in\Eh}\int_E (\bxi_u^n\cdot\nabla)\btheta_u^n\cdot\bu^n \dx  \nonumber\\
		&\quad- \frac{1}{2}\sum_{e\in\Gamma_h}\int_e \avr{\bxi_u^n} \cdot\bn_e \jump{\bu^n}\cdot\avr{\btheta_u^n}\ds + \frac{1}{2} \sum_{e\in\Gamma_h}\int_e\avr{ \bxi_u^n}\cdot\bn_e \avr{\bu^n}\cdot\jump{\btheta_u^n} \ds.
	\end{align}
	Applying the H{\"o}lder inequality, Young's inequality with the Gagliardo-Nirenberg inequality \eqref{GNI}, the Agmon's inequality \eqref{AI} and Lemmas \ref{disl42l2}, \ref{edge2elem}, \ref{dised2el} and the fact $[\bu^n]\cdot\bn_e=0$, we arrive at
	\begin{equation}
		|b_2(\bxi_u^n,\bu^n,\btheta_u^n)| \le C \left(
		\|\bxi_u^n\|^2 + h^2\|\bxi_u^n\|_u^2\right)\|\bu^n\|_{H^2}^2 + \frac{\nu}{16}\|\btheta_u^n\|_u^2.
	\end{equation}
	For the last term of \eqref{Lutheta}, one proceed similar to \eqref{utheta5}-\eqref{utheta7} to find 
	\begin{align}\label{utheta131}
		|b_2(\btheta_u^n,\bu^n,\btheta_u^n)| 
		& \le C\|\bu^n\|_{H^2}^2\|\btheta_u^n\|^2 + \frac{\nu}{16}\|\btheta_u^n\|_u^2.
	\end{align}
	From \eqref{Fu} and \eqref{errorsplit}, a use of the Cauchy-Schwarz inequality and the Young's inequality yields
	\begin{align}\label{utheta14}
		|F_u(\btheta_u^n)| &\le \left(\|e_\rho^{n-1}\|+\|\trho^n-\trho^{n-1}\|\right)\|\nabla\Phi\|_{\bL^4}\|\btheta_u^n\|_{\bL^4}\nonumber\\
		& \le C \left(\|\xi_\rho^{n-1}\|^2+\|\theta_\rho^{n-1}\|^2+\|\trho^n-\trho^{n-1}\|^2\right)\|\nabla\Phi\|_{\bL^4}^2 + \frac{\nu}{16}\|\btheta_u^n\|_u^2\nonumber\\
		& \le C \left(\|\xi_\rho^{n-1}\|^2+\|\theta_\rho^{n-1}\|^2+\Delta t\int_{t_{n-1}}^{t_n}\|\trho_t(t)\|^2\dt\right)\|\nabla\Phi\|_{\bL^4}^2 + \frac{\nu}{16}\|\btheta_u^n\|_u^2. 
	\end{align}
	Inserting all the estimates \eqref{utheta2}-\eqref{utheta14} in \eqref{utheta1} to obtain
	\begin{align}\label{utheta15}
		\pt\|\btheta_u^n\|^2 + \nu\|\btheta_u^n\|_u^2 
		&\lesssim \frac{1}{\Delta t}\int_{t_{n-1}}^{t_n} \|\bxi_{ut}(t)\|^2\dt + \Delta t \int_{t_{n-1}}^{t_n} \left(\|\bu_{tt}(t)\|^2+\|\trho_t(t)\|^2\|\nabla\Phi\|_{\bL^4}^2\right)\dt +  \|\bxi_u^n\|^4  \nonumber\\
		&\quad + \|\bu^n\|_{H^2}^2 \|\bxi_u^n\|^2+ h^2\|\bu^n\|_{H^2}^2\|\bxi_u^n\|_u^2 +  \|\nabla\Phi\|_{\bL^4}^2\|\xi_\rho^{n-1}\|^2+\|\nabla\Phi\|_{\bL^4}^2\|\theta_\rho^{n-1}\|^2 \nonumber\\
		& \quad + \left(\|\bxi_u^n\|_u^2 + \|\bxi_u^n\|^2\|\bxi_u^n\|_u^2 + h^{-3} \|\bxi_u^n\|^2 + h^{-1}\|\bxi_u^n\|_u^2  + \|\bu^n\|_{H^2}^2\right)\|\btheta_u^n\|^2 .
	\end{align}
	Now, we first use Lemma \ref{lem:ritz} and \ref{lem:stokes} with the regularity assumptions (\textbf{A1}). Then, multiply the resulting inequality by $\Delta t$ and take summation from $n=1$ to $N$ with the fact $\Delta t \sum_{n=1}^N \pt\|\btheta_u^n\|^2 = \|\btheta_u^N\|^2$ and $\|\btheta_u^0\|=0$ to  derive
	\begin{align}
		\|\btheta_u^N\|^2 + \nu\Delta t\sum_{n=1}^N\|\btheta_u^n\|_u^2 
		&\lesssim h^{2k+2}\int_{0}^{t_N} \left(\|\bu(t)\|_{\bH^{k+1}}^2 + \|p(t)\|_{H^k}+\|\bu_t(t)\|_{\bH^{k+1}}^2 + \|p_t(t)\|_{H^k} + \|\trho(t)\|_{H^{k+1}}\right)\dt  \nonumber\\
		& \qquad + (\Delta t)^2 \int_{0}^{t_N} \left(\|\bu_{tt}(t)\|^2+\|\trho_t(t)\|^2\right)\dt + \Delta t \sum_{n=1}^N \left(\|\theta_\rho^{n-1}\|^2+\|\btheta_u^n\|^2\right).
	\end{align}
	A use of regularity assumptions (\textbf{A1}) completes the rest of the proof.
\end{proof}

\begin{lemma}\label{lem:thetac}
	Under the assumption of Lemma \ref{lem:thetau}, the following inequality holds true
	\begin{align}
		\|\theta_c^N\|^2 + \kappa\Delta t\sum_{n=1}^N\|\theta_c^n\|_c^2  +\gamma m_0\Delta t\sum_{n=1}^N\|\theta_c^n\|^2 
		\lesssim \left(h^{2k+2} + (\Delta t)^2\right)   &+ \Delta t \sum_{n=1}^N \left(\|\theta_\rho^{n-1}\|^2+\|\theta_c^n\|^2+ \|\btheta_u^n\|^2\right) \nonumber\\
		&+  \Delta t \sum_{n=1}^N \|\theta_\rho^{n-1}\|^2\|\theta_c^n\|^2.
	\end{align}
\end{lemma}
\begin{proof}
	We choose $\psi_h = \theta_c^n$ in \eqref{erreqthetac} to obtain
	\begin{align}\label{ctheta1}
		\frac{1}{2}\pt\|\theta_c^n\|^2 + \kappa \|\theta_c^n\|_c^2 + \gamma m_0\|\theta_c^n\|^2 \le -(\pt\xi_c^n,\theta_c^n) - \gamma m_0(\xi_c^n,\theta_c^n) + E_c(\theta_c^n) + \Lambda_c(\theta_c^n) + F_c(\theta_c^n).
	\end{align}
	Processed similarly to \eqref{utheta2} and \eqref{utheta3}, one can bound the first and third term on the right hand side of \eqref{ctheta1} as
	\begin{align}\label{ctheta2}
		|(\pt\xi_c^n, \theta_c^n)| +|E_c(\theta_c^n)|
		\le  \frac{C}{\Delta t}\int_{t_{n-1}}^{t_n} \|\xi_{ct}(t)\|^2 \dt+ C\Delta t \int_{t_{n-1}}^{t_n} \|c_{tt}(t)\|^2\dt + \frac{\kappa}{8}\|\theta_c^n\|_c^2
	\end{align}
	An application of the Cauchy-Schwarz inequality and the Young's inequality in the second term of \eqref{ctheta1} yields
	\begin{equation}
		|\gamma m_0(\xi_c^n,\theta_c^n)|\le C \|\xi_c^n\|^2 + \frac{\gamma m_0}{2}\|\theta_c^n\|^2.
	\end{equation}
	We first use \eqref{usplit} to rewrite the nonlinear terms and then apply Lemma \ref{lem:skew} which gives 
	\begin{align}
		\Lambda_c(\theta_c^n) 
		= & b_1(\bxi_u^n,\xi_c^n,\theta_c^n)  + b_1(\btheta_u^n,\xi_c^n,\theta_c^n) - 
		b_1(\bu^n,\xi_c^n,\theta_c^n) -
		b_1(\bxi_u^n,c^n,\theta_c^n) - b_1(\btheta_u^n,c^n,\theta_c^n). \label{Lctheta}
	\end{align}
	One can consider the similar argument \eqref{utheta4}-\eqref{utheta131} to bound the above nonlinear terms as
	\begin{align}\label{ctheta5}
		|\Lambda_c(\theta_c^n)| &\le C\left( \left(\|\xi_c^n\|_c^2 + \|c^n\|_{H^2}^2\right)\|\bxi_u^n\|^2 + h^2\|\bxi_u^n\|_u^2\|c^n\|_{H^2}^2 \right.  \nonumber\\
		&\quad+ \left.\left(h^{-3} \|\xi_c^n\|^2 + h^{-1}\|\xi_c^n\|_u^2+ \|c^n\|_{H^2}^2\right)\|\btheta_u^n\|^2\right)  + \frac{\kappa}{8}\|\theta_c^n\|_c^2.
	\end{align}
	From \eqref{Fc} and \eqref{usplit}, we rewrite the last term on the right hand side of \eqref{ctheta2} as
	\begin{align}
		F_c(\theta_c^n) &= 
		(\xi_\rho^{n-1}\xi_c^n,\theta_c^n) + (\xi_\rho^{n-1}\theta_c^n,\theta_c^n) + (\theta_\rho^{n-1}\xi_c^n,\theta_c^n) + (\theta_\rho^{n-1}\theta_c^n,\theta_c^n) - (\trho^{n-1}\xi_c^n,\theta_c^n) \nonumber\\
		&\quad  - (\trho^{n-1}\theta_c^n,\theta_c^n) - (\xi_\rho^{n-1}c^n,\theta_c^n) - (\theta_\rho^{n-1}c^n,\theta_c^n) - ((\trho^n-\trho^{n-1})c^n,\theta_c^n).
	\end{align}
	A use of the H{\"o}lder inequality and the Young's inequality yields
	\begin{align}\label{ctheta7}
		|F_c(\theta_c^n)| \le & C\left(\|\xi_\rho^{n-1}\|^2\|\xi_c^n\|\|\xi_c^n\|_c + \|\nabla\trho^{n-1}\|^2\|\xi_c^n\|^2 + \|\nabla c^n\|^2\|\xi_\rho^{n-1}\|^2 + \left(\|\xi_\rho^{n-1}\|^2+\|\trho^{n-1}\|^2\right)\|\theta_c^n\|^2 \right.\nonumber\\
		& \quad \left.+ \left(\|\xi_c^n\|\|\xi_c^n\|_c+\|\nabla c^n\|^2\right)\|\theta_\rho^{n-1}\|^2 + \|\theta_\rho^{n-1}\|^2\|\theta_c^n\|^2 \right)  + \frac{\kappa}{8}\|\theta_c^n\|_c^2\nonumber\\
		&\quad + C\left(\Delta t \int_{t_{n-1}}^{t_n}\|\trho_t(t)\|^2\dt\right)\|\nabla c^n\|^2.
	\end{align}
	Combining all the estimates from \eqref{ctheta2}-\eqref{ctheta7} and put it in \eqref{ctheta1} to obtain
	\begin{align}\label{ctheta8}
		\pt\|\theta_c^n\|^2 + \kappa\|\theta_c^n\|_c^2 
		&\lesssim \frac{1}{\Delta t}\int_{t_{n-1}}^{t_n} \|\xi_{ct}(t)\|^2\dt + \Delta t \int_{t_{n-1}}^{t_n} \|c_{tt}(t)\|^2\dt + \left(\Delta t \int_{t_{n-1}}^{t_n}\|\trho_t(t)\|^2\dt\right)\|\nabla c^n\|^2   \nonumber\\
		& + \left(\|\xi_c^n\|_c^2 + \|c^n\|_{H^2}^2\right)\|\bxi_u^n\|^2 + h^2\|\bxi_u^n\|_u^2\|c^n\|_{H^2}^2 + \|\xi_\rho^{n-1}\|^2\|\xi_c^n\|\|\xi_c^n\|_c   \nonumber\\
		& + \left(1+\|\nabla\trho^{n-1}\|^2\right)\|\xi_c^n\|^2 + \|\nabla c^n\|^2\|\xi_\rho^{n-1}\|^2 + \left(h^{-3} \|\xi_c^n\|^2 + h^{-1}\|\xi_c^n\|_u^2+ \|c^n\|_{H^2}^2\right)\|\btheta_u^n\|^2  \nonumber\\
		& + \left(\|\xi_\rho^{n-1}\|^2+\|\trho^{n-1}\|^2\right)\|\theta_c^n\|^2 + \left(\|\xi_c^n\|\|\xi_c^n\|_c+\|\nabla c^n\|^2\right)\|\theta_\rho^{n-1}\|^2 + \|\theta_\rho^{n-1}\|^2\|\theta_c^n\|^2.
	\end{align}
	Now, we first use Lemma \ref{lem:ritz} and \ref{lem:stokes} with the regularity assumptions (\textbf{A1}). Then, multiply the resulting inequality by $\Delta t$ and take summation from $n=1$ to $N$ with the fact $\Delta t \sum_{n=1}^N \pt\|\theta_c^n\|^2 = \|\theta_c^N\|^2$ to  derive
	\begin{align}
		\|\theta_c^N\|^2 +  \Delta t\sum_{n=1}^N(\kappa\|\theta_c^n\|_c^2  + \gamma m_0 \|\theta_c^n\|^2)
		&  
		\lesssim h^{2k+2}\int_{0}^{t_N} \left(\|\bu(t)\|_{\bH^{k+1}}^2 + \|p(t)\|_{H^k} + \|\rho(t)\|_{H^{k+1}}\right)\dt  \nonumber\\
		& + h^{2k+2}\int_{0}^{t_N} \left(\|c(t)\|_{H^{k+1}}^2 + \|c_t(t)\|_{H^{k+1}}^2 \right)\dt \nonumber\\
		&+  (\Delta t)^2  \int_{0}^{t_N} \|c_{tt}(t)\|^2\dt +  (\Delta t)^2  \int_{0}^{t_N} \|\trho_{t}(t)\|^2\dt  \nonumber\\		
		& + \Delta t \sum_{n=1}^N \left(\|\theta_\rho^{n-1}\|^2+\|\theta_c^n\|^2+ \|\btheta_u^n\|^2 +  \|\theta_\rho^{n-1}\|^2\|\theta_c^n\|^2\right).
	\end{align}
	A use of regularity assumptions (\textbf{A1}) completes the rest of the proof.
\end{proof}

%%%%%%%%%%%%%%%%%%%%%%%%%%%%%%%%%%%%%%%%%%%%%%%%%%%%%%%%%%%%%%%%%%%%%%%%%%%%%%%%%%

\begin{lemma}\label{lem:thetarho}
	Under the assumption of Lemma \ref{lem:thetau}, the following inequality holds true
\begin{align*}
	\|\theta_\rho^N\|^2 +  \mu \Delta t\sum_{n=1}^N \|\theta_\rho^n\|_\rho^2  & 
	\lesssim \left(h^{2k+2} +  (\Delta t)^2\right)  + \Delta t \sum_{n=1}^N \left(\|\theta_\rho^{n-1}\|^2+\|\theta_c^n\|^2+ \|\btheta_u^n\|^2 +  h^{-4}\|\theta_\rho^{n-1}\|^2\|\theta_c^n\|^2\right).
\end{align*}
\end{lemma}
\begin{proof}
	Choose $\chi_h=\theta_\rho^n$ in \eqref{erreqthetarho} to obtain
	\begin{equation}\label{rhotheta1}
	\frac{1}{2}\pt\|\theta_\rho^n\|^2 + \mu\|\theta_\rho^n\|_\rho^2 \le -(\pt\xi_\rho^n, \theta_\rho^n) + E_\rho(\theta_\rho^n) + \Lambda_\rho(\theta_\rho^n) + \beta G_\rho(\theta_\rho^n) + \beta m_0 a_c(\theta_c^n,\theta_\rho^n).
	\end{equation}
	First three terms on the right hand side of \eqref{rhotheta1} can be bounded exactly similar to  \eqref{ctheta2}-\eqref{ctheta5} as
\begin{align}\label{rhotheta2}
		|(\pt\xi_\rho^n, \theta_\rho^n)| + |E_\rho(\theta_\rho^n)|  +  |\Lambda_\rho(\theta_\rho^n)|
		\le  C\bigg( \frac{1}{\Delta t}\int_{t_{n-1}}^{t_n} \|\xi_{\rho t}(t)\|^2 \dt + \Delta t \int_{t_{n-1}}^{t_n} \|\rho_{tt}(t)\|^2\dt +  \|\xi_\rho^n\|_\rho^2  \|\bxi_u^n\|^2 \nonumber\\
		+ \|\trho^n\|_{H^2}^2 \|\bxi_u^n\|^2 + h^2\|\bxi_u^n\|_u^2\|\trho^n\|_{H^2}^2 + \left(h^{-3} \|\xi_\rho^n\|^2 + h^{-1}\|\xi_\rho^n\|_u^2+ \|\trho^n\|_{H^2}^2\right)\|\btheta_u^n\|^2\bigg) + \frac{\mu}{8}\|\theta_\rho^n\|_\rho^2.
\end{align}
We rewrite the fourth term of the right hand of \eqref{rhotheta1} using \eqref{Frho} and \eqref{usplit} as
\begin{align}\label{grho}
	G_\rho(\theta_\rho^n) &= 
	g(\xi_\rho^{n-1},R_hc^n,\theta_\rho^n) +  g(\theta_\rho^{n-1},R_hc^n,\theta_\rho^n) + g(\trho^{n-1},\theta_c^n,\theta_\rho^n) - g(\xi_\rho^{n-1},\theta_c^n,\theta_\rho^n) \nonumber\\
	&\quad   - g(\theta_\rho^{n-1},\theta_c^n,\theta_\rho^n) + g(\trho^n-\trho^{n-1},R_hc^n,\theta_\rho^n) +   g(m_0,\theta_c^n,\theta_\rho^n).
\end{align}
From the definition of $g(\cdot,\cdot,\cdot)$, it follows that
\begin{align}\label{groh1}
	g(\xi_\rho^{n-1}, R_hc^n,\theta_\rho^n) 
	& = \sum_{E\in\Eh} \int_E \xi_\rho^{n-1}\nabla R_hc^n\cdot\nabla\theta_\rho^n \ds - \sum_{e\in\Gamma_h^0} \int_e \avr{\xi_\rho^{n-1}\nabla R_hc^n}\cdot\bn_e \jump{\theta_\rho^n} \ds \nonumber\\
	& \qquad- \sum_{e\in\Gamma_h^0} \int_e  \avr{\xi_\rho^{n-1}\nabla \theta_\rho^n}\cdot\bn_e \jump{R_hc^n} \ds   = G_{11} + G_{12} + G_{13}.
\end{align}
One can bound $G_{11}$ by using the H{\"o}lder inequality and the Young's inequality as	
\begin{align}\label{groh1g11}
	|G_{11}| 
	\lesssim  \|\xi_\rho^{n-1}\| \|\nabla R_hc^n\|_{L^\infty} \|\theta_\rho^n\|_\rho 
	\le  C\|\nabla R_hc^n\|_{L^\infty}^2  \|\xi_\rho^{n-1}\|^2  + \frac{\mu}{32} \|\theta_\rho^n\|_\rho^2.
\end{align}
A use of the H{\"o}lder inequality with the Young's inequality, the inverse inequality and Lemma \ref{edge2elem} yields	
\begin{align}\label{groh1g12}
|G_{12}| & \le   \sum_{e\in\Gamma_h^0} \|\xi_\rho^{n-1}\|_{L^2(e)}\|\nabla R_hc^n\|_{L^\infty(e)}\|\jump{\theta_\rho^n}\|_{L^2(e)} \nonumber\\
& \lesssim   \|\nabla R_hc^n\|_{L^\infty}\left(\sum_{e\in\Gamma_h^0} \frac{|e|}{\sigma_\rho}\|\xi_\rho^{n-1}\|_{L^2(e)}^2\right)^{\frac{1}{	2}}\left(\sum_{e\in\Gamma_h^0}\frac{\sigma_\rho}{|e|}\|\jump{\theta_\rho^n}\|_{L^2(e)}^2\right)^{\frac{1}{2}} \nonumber\\
& \le C \|\nabla R_hc^n\|_{L^\infty}^2  \sum_{E\in\Eh} h_E(h_E^{-1}\|\xi_\rho^{n-1}\|_{L^2(E)}^2+h_E\|\nabla\xi_\rho^n\|^2)  + \frac{\mu}{64} \|\theta_\rho^n\|_\rho^2 \nonumber\\
& \le C\|\nabla R_hc^n\|_{L^\infty}^2 \left(\|\xi_\rho^{n-1}\|^2
+ h^2\|\xi_\rho^{n-1}\|_\rho^2\right) + \frac{\mu}{64} \|\theta_\rho^n\|_\rho^2,
\end{align}
and	using the fact $\jump{R_hc^n} = \jump{R_hc^n-c^n} = - \jump{\xi^n}$, we find that
\begin{align}\label{groh1g13}
	|G_{13}| & \le  \sum_{e\in\Gamma_h^0} \|\xi_\rho^{n-1}\|_{L^2(e)} \|\nabla \theta_\rho^n\|_{L^\infty(e)} \|\jump{R_hc^n-c^n}\|_{L^2(e)} \nonumber\\
	& \lesssim  \left(\sum_{e\in\Gamma_h^0} \frac{|e|}{\sigma_\rho}\|\xi_\rho^{n-1}\|_{L^2(e)}^2\|\nabla\theta_\rho^n\|_{L^\infty(e)}^2\right)^{\frac{1}{2}}\left(\sum_{e\in\Gamma_h^0}\frac{\sigma_c}{|e|}\|\jump{\xi_c^n}\|_{L^2(e)}^2\right)^{\frac{1}{2}} \nonumber\\
	& \lesssim   \|\xi_c^n\|_c  \left(\sum_{E\in\Eh} h_E(h_E^{-1}\|\xi_\rho^{n-1}\|_{L^2(E)}^2+h_E\|\nabla\xi_\rho^n\|^2)h_E^{-2}\|\nabla\theta_\rho^n\|_{L^2(E)}^2 \right)^{\frac{1}{2}} \nonumber\\
	& \le C h^{-2}\|\xi_c^n\|_c^2 \left(\|\xi_\rho^{n-1}\|^2
	+ h^2\|\xi_\rho^{n-1}\|_\rho^2\right) + \frac{\mu}{64} \|\theta_\rho^n\|_\rho^2.
\end{align}
The term $\|\nabla R_hc^n\|_{L^\infty}$ located in \eqref{groh1g11} and \eqref{groh1g12} can be bounded as below
\begin{align}
	\|\nabla R_hc^n\|_{L^\infty(E)}^2 
	&\lesssim \|\nabla c^n\|_{L^\infty(E)}^2 + \|\nabla (c^n-R_hc^n)\|^2_{L^\infty(E)}\nonumber\\
	&\lesssim  \|\nabla c^n\|_{H^2(E)}^2 + \|\nabla (c^n-R_hc^n)\|_{L^2(E)}\|\nabla (c^n-R_hc^n)\|_{H^2(E)} \nonumber \\
	&\lesssim \|c^n\|_{H^3(E)}^2 + h^{2}\|c^n\|_{H^3(E)}\left( \|\nabla(c^n-\mathcal{I}_hc^n)\|_{H^2(E)} + \|\nabla(\mathcal{I}_hc^n-R_hc^n)\|_{H^2(E)}\right)\nonumber\\
	& \lesssim \|c^n\|_{H^3(E)}^2,
\end{align}	
where $\mathcal{I}_h:L^2(E)\to \mathbb{P}_{r+1}(E)$ is the optimal polygonal approximation defined in Lemma \ref{polyapprox}. 
In a similar way to \eqref{groh1} with the inverse inequality, we can bound the second term on the right hand side of \eqref{grho} as
\begin{align}\label{groh2}
	g(\theta_\rho^{n-1}, R_hc^n,\theta_\rho^n) 
 \le C\|c^n\|_{H^3}^2  \|\theta_\rho^{n-1}\|^2
	  + \frac{\mu}{16} \|\theta_\rho^n\|_\rho^2. 
\end{align}
For the third term of \eqref{grho}, a use of the H{\"o}lder inequality with the Young's inequality and the Agmon's inequality \eqref{AI} shows
\begin{align}\label{groh3}
	g(\trho^{n-1}, \theta_c^n,\theta_\rho^n) 
	& \lesssim  \|\trho^{n-1}\|_{L^\infty} \|\nabla\theta_c^n\|\|\nabla\theta_\rho^n\| 
	+  \|\trho^{n-1}\|_{L^\infty}\left(\sum_{e\in\Gamma_h^0} \frac{|e|}{\sigma_\rho}\|\nabla\theta_c^n\|_{L^2(e)}^2\right)^{\frac{1}{	2}}\left(\sum_{e\in\Gamma_h^0}\frac{\sigma_\rho}{|e|}\|\jump{\theta_\rho^n}\|_{L^2(e)}^2\right)^{\frac{1}{2}} \nonumber\\
	&\qquad+  \|\trho^{n-1}\|_{L^\infty}\left(\sum_{e\in\Gamma_h^0} \frac{|e|}{\sigma_c}\|\nabla\theta_\rho^n\|_{L^2(e)}^2\right)^{\frac{1}{	2}}\left(\sum_{e\in\Gamma_h^0}\frac{\sigma_c}{|e|}\|\jump{\theta_c^n}\|_{L^2(e)}^2\right)^{\frac{1}{2}} \nonumber\\
	& \le C \|\trho^{n-1}\|_{H^2}^2 \|\theta_c^n\|_c^2
	+ \frac{\mu}{16} \|\theta_\rho^n\|_\rho^2.  
\end{align}
From the definition of $g(\cdot,\cdot,\cdot)$, the fourth term of \eqref{grho} can be write as
\begin{align}\label{groh4}
	g(\xi_\rho^{n-1}, \theta_c^n,\theta_\rho^n) 
	& = \sum_{E\in\Eh} \int_E \xi_\rho^{n-1}\nabla \theta_c^n\cdot\nabla\theta_\rho^n \ds - \sum_{e\in\Gamma_h^0} \int_e \avr{\xi_\rho^{n-1}\nabla \theta_c^n}\cdot\bn_e \jump{\theta_\rho^n} \ds \nonumber\\
	&\qquad - \sum_{e\in\Gamma_h^0} \int_e \avr{\xi_\rho^{n-1}\nabla \theta_\rho^n}\cdot\bn_e \jump{\theta_c^n}\ds = G_{41}+G_{42}+G_{43}.
\end{align}
An application of the H{\"o}lder inequality with the Young's inequality and the inverse hypothesis helps to bound $G_{41}$ as
\begin{align}\label{groh4g1}
	|G_{41}| \le  \|\xi_\rho^{n-1}\| \|\nabla \theta_c^n\|_{L^\infty} \|\theta_\rho^n\|_\rho 
	\le  Ch^{-2} \|\xi_\rho^{n-1}\|^2
	 \|\theta_c^n\|_c^2 + \frac{\mu}{64} \|\theta_\rho^n\|_\rho^2.  
\end{align}
To bound the $G_{42}$ term, we use the H{\"o}lder inequality, the Young's inequality, the inverse hypothesis \eqref{inv.hypo}, Lemma \ref{edge2elem} and the Agmon's inequality \eqref{AI} to derive
\begin{align}\label{groh4g2}
	|G_{42}| & \le \sum_{e\in\Gamma_h^0} \|\xi_\rho^{n-1}\|_{L^2(e)}\|\nabla \theta_c^n\|_{L^\infty(e)}\|\jump{\theta_\rho^n}\|_{L^2(e)} \nonumber\\
	& \lesssim  \|\nabla \theta_c^n\|_{L^\infty} \left(\sum_{e\in\Gamma_h^0} \frac{|e|}{\sigma_\rho}\|\xi_\rho^{n-1}\|_{L^2(e)}^2\right)^{\frac{1}{	2}}\left(\sum_{e\in\Gamma_h^0}\frac{\sigma_\rho}{|e|}\|\jump{\theta_\rho^n}\|_{L^2(e)}^2\right)^{\frac{1}{2}} \nonumber\\
	& \lesssim   \|\nabla \theta_c^n\|_{L^\infty} \left( \sum_{E\in\Eh} h_E(h_E^{-1}\|\xi_\rho^{n-1}\|_{L^2(E)}^2+h_E\|\nabla\xi_\rho^n\|^2)\right)^{\frac{1}{2}} \|\theta_\rho^n\|_\rho \nonumber\\
	& \le Ch^{-2}\left(\|\xi_\rho^{n-1}\|^2
	+ h^2\|\xi_\rho^{n-1}\|_\rho^2\right)\|\theta_c^n\|_c^2 + \frac{\mu}{64} \|\theta_\rho^n\|_\rho^2.  
\end{align}
Considering a similar argument as above, one can bound $G_{43}$ as follows
\begin{align}\label{groh4g3}
	|G_{43}| &\le  \sum_{e\in\Gamma_h^0} \|\xi_\rho^{n-1}\|_{L^2(e)} \|\nabla \theta_\rho^n\|_{L^\infty(e)}\|\jump{\theta_c^n}\|_{L^2(e)} \nonumber\\
	& \lesssim  \|\nabla \theta_\rho^n\|_{L^\infty} \left(\sum_{e\in\Gamma_h^0} \frac{|e|}{\sigma_c} \|\xi_\rho^{n-1}\|_{L^2(e)}^2\right)^{\frac{1}{2}} \left(\sum_{e\in\Gamma_h^0}\frac{\sigma_c}{|e|}\|\jump{\theta_c^n}\|_{L^2(e)}^2\right)^{\frac{1}{2}} \nonumber\\
	& \le  Ch^{-2}\left(\|\xi_\rho^{n-1}\|^2
	+ h^2\|\xi_\rho^{n-1}\|_\rho^2\right)\|\theta_c^n\|_c^2 + \frac{\mu}{32} \|\theta_\rho^n\|_\rho^2.
\end{align}
Arguing similar as \eqref{groh4} (use Lemma \ref{dised2el} instead of Lemma \ref{edge2elem}) with the application of inverse inequality \eqref{inv.hypo}, fifth term can be estimates as
\begin{align}\label{groh5}
	g(\theta_\rho^{n-1}, \theta_c^n,\theta_\rho^n) 
	 \le Ch^{-2}\|\theta_\rho^{n-1}\|^2
	 \|\theta_c^n\|_c^2 + \frac{\mu}{16} \|\theta_\rho^n\|_\rho^2 \le Ch^{-4}\|\theta_\rho^{n-1}\|^2
	 \|\theta_c^n\|^2 + \frac{\mu}{16} \|\theta_\rho^n\|_\rho^2.
\end{align}
Sixth term of \eqref{grho} can be bounded similar to the first term of \eqref{grho} by replacing $\xi_\rho^{n-1}$ by $\trho^n-\trho^{n-1}$ as
\begin{align}
	|g(\trho^n-\trho^{n-1},R_hc^n,\theta_\rho^n)| &\le C\left(\|c^n\|_{H^3}^2+h^{-2}\|\xi_c^n\|_c^2\right)\left(\|\trho^n-\trho^{n-1}\|^2+h^2\|\nabla(\trho^n-\trho^{n-1})\|^2\right) + \frac{\mu}{32}\|\theta_\rho^n\|_\rho^2 \nonumber\\
	&\le C \Delta t \int_{t_{n-1}}^{t_n} \left(\|\trho_t(t)\|^2+h^2\|\nabla\trho_t(t)\|^2\right) \dt + \frac{\mu}{32}\|\theta_\rho^n\|_\rho^2.
\end{align}
Using \eqref{defg}, the last term of \eqref{grho}, can be bounded as
\begin{align}\label{grho6}
	g(m_0,\theta_c^n,\theta_\rho^n) \le C \|\theta_c^n\|^2_c + \frac{\mu}{32}\|\theta_\rho^n\|_\rho^2.
\end{align}
Inserting all the estimates from \eqref{rhotheta2}-\eqref{grho6} in \eqref{rhotheta1} to obtain
\begin{align}
	\partial_t\|\theta_\rho^n\|^2 + \mu \|\theta_\rho^n\|^2_\rho &\lesssim \frac{1}{\Delta t}\int_{t_{n-1}}^{t_n} \|\xi_{\rho t}(t)\|^2 \dt + \Delta t \int_{t_{n-1}}^{t_n} \left(\|\rho_{tt}(t)\|^2 + \|\trho_t(t)\|^2+h^2\|\nabla\trho_t(t)\|^2\right)\dt   \nonumber\\
	&\quad + \left(\|\xi_\rho^n\|_\rho^2 + \|\trho^n\|_{H^2}^2\right)\|\bxi_u^n\|^2  + \left(h^{-3} \|\xi_\rho^n\|^2 + h^{-1}\|\xi_\rho^n\|_u^2+ \|\trho^n\|_{H^2}^2\right)\|\btheta_u^n\|^2 \nonumber\\
	&\quad + h^2\|\bxi_u^n\|_u^2\|\trho^n\|_{H^2}^2 +\|c^n\|_{H^3}^2 \left(\|\xi_\rho^{n-1}\|^2
	+ h^2\|\xi_\rho^{n-1}\|_\rho^2\right) + \|c^n\|_{H^3}^2  \|\theta_\rho^{n-1}\|^2  \nonumber \\
	&\quad + h^{-4}\|\theta_\rho^{n-1}\|^2
	\|\theta_c^n\|^2 + \left(\|\trho^{n-1}\|_{H^2}^2 + h^{-2} \|\xi_\rho^{n-1}\|^2 
	+ \|\xi_\rho^{n-1}\|_\rho^2+1\right)\|\theta_c^n\|_c^2.
\end{align}
Now, we first use Lemma \ref{lem:ritz} and \ref{lem:stokes} with the regularity assumptions (\textbf{A1}). Then, multiply the resulting inequality by $\Delta t$ and take summation from $n=1$ to $N$ with the fact $\Delta t \sum_{n=1}^N \pt\|\theta_c^n\|^2 = \|\theta_c^N\|^2$ to  derive	
\begin{align}
	\|\theta_\rho^N\|^2 +  \mu \Delta t\sum_{n=1}^N \|\theta_\rho^n\|_\rho^2   
	\lesssim & h^{2k+2}\int_{0}^{t_N} \left(\|\bu(t)\|_{\bH^{k+1}}^2 + \|p(t)\|_{H^k}+\|c(t)\|_{H^{k+1}}^2 + \|\rho(t)\|_{H^{k+1}} + \|\rho_t(t)\|_{H^{k+1}}\right)\dt  \nonumber\\
	&+  (\Delta t)^2  \int_{0}^{t_N}\left(\|\rho_{tt}(t)\|^2 + \|\trho_t(t)\|^2+h^2\|\nabla\trho_t(t)\|^2\right)\dt \nonumber\\
	& + \Delta t \sum_{n=1}^N \left(\|\theta_\rho^{n-1}\|^2+\|\theta_c^n\|_c^2+ \|\btheta_u^n\|^2 +  h^{-4}\|\theta_\rho^{n-1}\|^2\|\theta_c^n\|^2\right).
\end{align}
A use of Lemma \ref{lem:thetac} with regularity estimates (\textbf{A1}) completes the rest of the proof.
\end{proof}

Now, we are ready to prove the final result of this section.

\begin{theorem}\label{thm:l2u}
	Let the assumption (\textbf{A1}) holds true. Suppose $(\bu,c,\trho)$ and $(\bu_h^n,c_h^n,\trho_h^n)$ be the solutions of the continuous system \eqref{dg1rho}-\eqref{dg4div} and the fully discrete system \eqref{fddg1rho}-\eqref{fddg4div}, respectively. Then, with $\e_u^n = \bu^n - \bu_h^n, ~
	e_\rho^n = \trho^n - \trho_h^n, ~ e_c^n = c^n - c_h^n$, the following holds
	\begin{align}
	\|\e_u^n\| + \|e_\rho^n\| + \|e_c^n\| \lesssim \left(h^{k+1}+\Delta t\right).
 	\end{align}
\end{theorem}
\begin{proof}
	Since $\e_u^n=\bxi_u^n+\btheta_u^n$, $e_c^n=\xi_c^n+\theta_c^n$ and $e_\rho^n = \xi_\rho^n+\theta_\rho^n$. And we have the estimates of $\bxi_u^n$, $\xi_c^n$ and $\xi_\rho$. It is enough to find the estimates for $\btheta_u^n$, $\theta_c^n$ and $\theta_\rho^n$.
	From Lemmas \ref{lem:thetau}, \ref{lem:thetac} and \ref{lem:thetarho}, it follows that
	\begin{align}\label{esth1}
	\|\btheta_u^N\|^2 + \|\theta_c^N\|^2+ \|\theta_\rho^N\|^2  +  \Delta t\sum_{n=1}^N  \left(\nu\|\btheta_u^n\|_u^2 + \kappa\|\theta_c^n\|_c^2 + m_0\gamma\|\theta_c^n\|^2 + \mu \|\theta_\rho^n\|_\rho^2\right)   
	\lesssim \left(h^{2k+2} +  (\Delta t)^2\right)\nonumber\\
	+ \Delta t \sum_{n=1}^N \left(\|\theta_\rho^{n-1}\|^2+\|\theta_c^n\|^2+ \|\btheta_u^n\|^2 +  h^{-4}\|\theta_\rho^{n-1}\|^2\|\theta_c^n\|^2\right).
	\end{align}
	Due to the presence of a multiplication term on the right hand side of the above inequality, we can not use the discrete Gronwall's lemma directly. But, with the help of the method of induction, one can proceed as follows:\\
	\textbf{Claim:} For any finite $m$, \eqref{esth1} implies
	\begin{align}\label{esth}
		\|\btheta_u^m\|^2 + \|\theta_c^m\|^2+ \|\theta_\rho^m\|^2  \lesssim \left(h^{2k+2} +  (\Delta t)^2\right).
	\end{align}
	For $m=1$ with $\theta_\rho^0=0$ and dropping the positive term from the left-hand side of \eqref{esth1}, we deduce 
	\begin{align}
		\|\btheta_u^1\|^2 + \|\theta_c^1\|^2+ \|\theta_\rho^1\|^2    
		\le C\left(h^{2k+2} +  (\Delta t)^2\right) 
		+ C\Delta t \left( \|\theta_c^1\|^2+ \|\btheta_u^1\|^2\right).
	\end{align}
	One can choose $\Delta t_1$ in such a way that $0<\Delta t < \Delta t_1 :=\frac{1}{2C}<\frac{1}{C}$, which implies $1\le \dfrac{1}{1-C\Delta t}\le 2$. Hence, \eqref{esth} is true for $m=1$. Again for $m=2$, we obtain
	\begin{align}
		\|\btheta_u^2\|^2 + \|\theta_c^2\|^2+ \|\theta_\rho^2\|^2    
		\le C\left(h^{2k+2} +  (\Delta t)^2\right) 
		+ C\Delta t \left( \left(1+Ch^{-4}\left(h^{2k+2}+(\Delta t)^2\right)\right)\|\theta_c^2\|^2+ \|\btheta_u^2\|^2\right).
	\end{align}
	With the assumption on $\Delta t < \min\{\Delta t_1, h^{2}\}$, \eqref{esth} is true for $m=2$. In a similar fashion, since results are known for $m\le N-1$, it is enough to show $m=N$ to complete the induction. Hence,
	\begin{align}
		\|\btheta_u^N\|^2 + \|\theta_c^N\|^2+ \|\theta_\rho^N\|^2  
		& \le C\left(h^{2k+2} +  (\Delta t)^2\right)
		+ C\Delta t \left(\|\theta_\rho^{N-1}\|^2+\|\theta_c^N\|^2+ \|\btheta_u^N\|^2 +  h^{-4}\|\theta_\rho^{N-1}\|^2\|\theta_c^N\|^2\right) \nonumber\\
		& \le C\left(h^{2k+2} +  (\Delta t)^2\right)
		+ C\Delta t \left(\left(1+Ch^{-4}\left(h^{2k+2}+(\Delta t)^2\right)\right)\|\theta_c^N\|^2+ \|\btheta_u^N\|^2\right). \nonumber
	\end{align}
	Considering the same assumption on $\Delta t$,  we obtain
	\begin{align}\label{esth2}
		\|\btheta_u^N\|^2 + \|\theta_c^N\|^2+ \|\theta_\rho^N\|^2  \lesssim \left(h^{2k+2} +  (\Delta t)^2\right).
	\end{align}
	A use of Lemmas \ref{lem:ritz}, \ref{lem:modi} and \ref{lem:stokes} with \eqref{esth2} completes the rest of the proof.
\end{proof}

\begin{remark}\label{rem:l2l2u}
	A use of  Lemmas \ref{lem:ritz}, \ref{lem:modi} and \ref{lem:stokes} and Theorem \ref{thm:l2u} with \eqref{esth1} yields
	\begin{align}
		\Delta t\sum_{n=1}^N  \left(\nu\|\e_u^n\|_u^2 + \kappa\|e_c^n\|_c^2 + m_0\gamma\|e_c^n\|^2 + \mu \|e_\rho^n\|_\rho^2\right)   
		\lesssim \left(h^{2k+2} +  (\Delta t)^2\right).
	\end{align}
\end{remark}

%%%%%%%%%%%%%%%%%%%%%%%%%%%%%%%%%%%% H1 Error %%%%%%%%%%%%%%%%%%

\subsection{$H^1$-Error Estimates}

In this section, we find $H^1$-errors for the fluid velocity, the cell density and the concentration.

\begin{lemma}\label{lem:h1u}
	Under the assumption of Lemma \ref{lem:thetau}, the following inequality holds true
	\begin{align} \label{h1error}
		\nu \|\e_u^N\|_u^2 +  \Delta t\sum_{n=1}^N\|\pt\e_u^n\|^2 
		\lesssim \left(h^{2k} + (\Delta t)^2\right).
	\end{align}
\end{lemma}

\begin{proof}
	To prove \eqref{h1error}, we first rewrite \eqref{erreqtheta} 
	\begin{equation}\label{h1erreqtheta}
		(\pt\btheta_u^n, \bv_h) + \nu~a_u(\btheta_u^n,\bv_h) - d(\bv_h,\theta_p^n) + d(\pt\btheta_u^n,q_h) = -(\pt\bxi_u^n, \bv_h)  + E_u(\bv_h) + \Lambda_u(\bv_h) + F_u(\bv_h).
	\end{equation}
	Here, we change the form of the continuity equation in terms of $\pt$  by taking the difference of time level $t_n$ and $t_{n-1}$ and dividing by $\Delta t$. 
	Now, choose $\bv_h = \pt\btheta_u^n$ and $q_h = \theta_p^n$ in \eqref{h1erreqtheta} and use the fact $a_u(\btheta_u^n,\pt\btheta_u^n)\ge \frac{1}{2}\pt\|\btheta_u^n\|_u^2$ to obtain
	\begin{equation}\label{h1utheta1}
		\|\pt\btheta_u^n\|^2 + 	\frac{\nu}{2}\pt\|\btheta_u^n\|_u^2 \le  -(\pt\bxi_u^n, \pt\btheta_u^n)  + E_u(\pt\btheta_u^n) + \Lambda_u(\pt\btheta_u^n) + F_u(\pt\btheta_u^n).
	\end{equation}
	A use of the Cauchy-Schwarz inequality and the Young's inequality with \eqref{Eu} helps to bound the following terms on the right hand side of \eqref{h1utheta1}:
	\begin{align}\label{h1utheta2}
		|(\pt\bxi_u^n, \pt\btheta_u^n)| + |E_u(\pt\btheta_u^n)|
		\le \frac{C}{\Delta t}\int_{t_{n-1}}^{t_n} \|\bxi_{ut}(t)\|^2 \dt + C \Delta t \int_{t_{n-1}}^{t_n} \|\bu_{tt}(t)\|^2\dt + \frac{1}{16}\|\pt\btheta_u^n\|^2.
	\end{align}
	From \eqref{Lu}, we rewrite the nonlinear terms on the right hand side of \eqref{h1utheta1} as
	\begin{align}
		\Lambda_u(\pt\btheta_u^n) = & b_2(\e_u^n,\e_u^n,\pt\btheta_u^n) - b_2(\e_u^n,\bu^n,\pt\btheta_u^n) - b_2(\bu^n,\e_u^n,\pt\btheta_u^n). \label{h1Lutheta}
	\end{align}
	For the first term of \eqref{h1Lutheta}, we use the H{\"o}lder inequality with the Lemmas \ref{edge2elem} and \ref{disl42l2}, inverse inequality \eqref{inv.hypo} and the Young's inequality to bound
	\begin{align}\label{h1utheta4}
		|b_2(\e_u^n,\e_u^n,\pt\btheta_u^n)| 
		&\le \sum_{E\in\Eh}\|\e_u^n\|_{\bL^2(E)} \|\nabla\e_u^n\|_{\bL^2(E)} \|\pt\btheta_u^n\|_{\bL^\infty(E)} + \frac{1}{2}\sum_{E\in\Eh}\|\nabla\cdot\e_u^n\|_{\bL^2(E)} \|\e_u^n\|_{\bL^2(E)}\|\pt\btheta_u^n\|_{\bL^\infty(E)} \nonumber\\
		& \quad+ \sum_{e\in\Gamma_h} \|\e_u^n\|_{\bL^2(e)} \|\jump{\e_u^n}\|_{\bL^2(e)} \|\pt\btheta_u^n\|_{\bL^\infty(e)} + \frac{1}{2} \sum_{e\in\Gamma_h} \|\jump{\e_u^n}\|_{\bL^2(e)} \|\e_u^n\|_{\bL^2(e)}\|\pt\btheta_u^n\|_{\bL^\infty(e)}  \nonumber\\
		& \le Ch^{-1}\|\e_u^n\|\|\e_u^n\|_u\|\pt\btheta_u^n\| \le  Ch^{-2}\|\e_u^n\|^2 \|\e_u^n\|_u^2 + \frac{1}{16}\|\pt\btheta_u^n\|^2.
	\end{align}
	A use of the H{\"older} inequality with the Gagliardo-Nirenberg inequality \eqref{GNI}, the Agmon's inequality \eqref{AI}, Lemma \ref{disl42l2} and the Young's inequality yields
	\begin{align}\label{h1utheta5}
		|b_2(\e_u^n,\bu^n,\pt\btheta_u^n)| 
		& \lesssim  \sum_{E\in\Eh} \|\e_u^n\|_{\bL^4(E)}\|\nabla\bu^n\|_{\bL^4(E)}  \|\pt\btheta_u^n\|_{\bL^2(E)} +  \sum_{E\in\Eh}\|\nabla\cdot\e_u^n\|_{\bL^2(E)}\|\bu^n\|_{\bL^\infty(E)}\|\pt\btheta_u^n\|_{\bL^2(E)} \nonumber\\
		&\quad+  \sum_{e\in\Gamma_h}\|\jump{\e_u^n}\|_{\bL^2(e)} \|\bu^n\|_{\bL^\infty(e)}\|\pt\btheta_u^n\|_{\bL^2(e)} \nonumber\\
		& \lesssim \|\e_u^n\|_u \|\bu^n\|_{\bH^2} \|\pt\btheta_u^n\| +  \|\bu^n\|_{\bL^\infty}\left(\sum_{e\in\Gamma_h}\frac{\sigma_u}{|e|}\|\jump{\e_u^n}\|_{\bL^2(e)}^2\right)^{\frac{1}{2}} \left(\sum_{e\in\Gamma_h} |e|\|\pt\btheta_u^n\|_{\bL^2(e)}^2\right)^{\frac{1}{2}} \nonumber\\
		& \lesssim \|\e_u^n\|_u \|\bu^n\|_{\bH^2} \|\pt\btheta_u^n\|  \le C\|\e_u^n\|_u^2 \|\bu^n\|_{\bH^2}^2 + \frac{1}{8} \|\pt\btheta_u^n\|^2.
	\end{align}
	Arguing similar to \eqref{h1utheta5}, the third term can be estimated as
	\begin{align}\label{h1utheta6}
		|b_2(\bu^n,\e_u^n,\pt\btheta_u^n)| 
		& \lesssim \|\e_u^n\|_u \|\bu^n\|_{\bH^2} \|\pt\btheta_u^n\|  \le C\|\e_u^n\|_u^2 \|\bu^n\|_{\bH^2}^2 + \frac{1}{8} \|\pt\btheta_u^n\|^2.
	\end{align}
	From \eqref{Fu} and \eqref{errorsplit}, a use of the Cauchy-Schwarz inequality with the Young's inequality yields
	\begin{align}\label{h1utheta14}
		|F_u(\pt\btheta_u^n)| &\le \left(\|e_\rho^{n-1}\|_{L^4}+\|\trho^n-\trho^{n-1}\|_{L^4}\right)\|\nabla\Phi\|_{\bL^4}\|\pt\btheta_u^n\| \nonumber\\
		&\le C \left(\|\e_\rho^{n-1}\|_\rho^2 + \Delta t \int_{t_{n-1}}^{t_n}\|\nabla\trho_t(t)\|^2\dt\right) \|\nabla\Phi\|_{\bL^4}^2 + \frac{1}{8}\|\pt\btheta_u^n\|^2. 
	\end{align}
	Inserting all the estimates \eqref{h1utheta2}-\eqref{h1utheta14} in \eqref{h1utheta1} to obtain
	\begin{align}\label{h1utheta15}
		\nu\pt\|\btheta_u^n\|_u^2 + \|\pt\btheta_u^n\|^2 
		&\lesssim \frac{1}{\Delta t}\int_{t_{n-1}}^{t_n} \|\bxi_{ut}(t)\|^2\dt + \Delta t \int_{t_{n-1}}^{t_n} \|\bu_{tt}(t)\|^2\dt + \left( \Delta t \int_{t_{n-1}}^{t_n}\|\nabla\trho_t(t)\|^2\dt\right) \|\nabla\Phi\|_{\bL^4}^2 \nonumber\\
		& \qquad+ \left(h^{-2}\|\e_u^n\|^2 + \|\bu^n\|_{H^2}^2\right)
		\|\e_u^n\|_u^2  + \|\e_\rho^{n-1}\|_\rho^2 \|\nabla\Phi\|_{\bL^4}^2.
	\end{align}
	Now, multiply the above inequality by $\Delta t$ and take summation from $n=1$ to $N$ with $\|\btheta_u^0\|_u=0$ and the fact $\Delta t \sum_{n=1}^N \pt\|\btheta_u^n\|_u^2 = \|\btheta_u^N\|_u^2$, we  derive
	\begin{align}
		\nu\|\btheta_u^N\|_u^2 +  \Delta t\sum_{n=1}^N\|\pt\btheta_u^n\|^2 
		&\lesssim h^{2k+2}\int_{0}^{t_N} \left(\|\bu_t(t)\|_{\bH^{k+1}}^2 + \|p_t(t)\|_{H^k}\right)\dt  + (\Delta t)^2 \int_{0}^{t_N} \left(\|\bu_{tt}(t)\|^2+\|\nabla\trho_t(t)\|^2\right)\dt  \nonumber\\
		& \qquad + \Delta t \sum_{n=1}^N \left(h^{-2}\|\e_u^n\|^2 + \|\bu^n\|_{H^2}^2\right)
		\|\e_u^n\|_u^2  + \|\e_\rho^{n-1}\|_\rho^2 \|\nabla\Phi\|_{\bL^4}^2.
	\end{align}
	A use of Theorem \ref{thm:l2u}, Remark \ref{rem:l2l2u} with the regularity assumptions (\textbf{A1}) and Lemma \ref{lem:stokes} completes the rest of the proof.
\end{proof}

\begin{lemma}\label{lem:h1c}
	Under the assumption of Lemma \ref{lem:thetau}, the following inequality holds true
	\begin{align*}
		\kappa\|\theta_c^N\|_c^2 + \gamma m_0\|\theta_c^N\|^2 +  \Delta t\sum_{n=1}^N \|\pt\theta_c^n\|^2
		&  
		\lesssim \left(h^{2k} + (\Delta t)^2\right).
	\end{align*}
\end{lemma}
\begin{proof}
	Choose $\psi_h = \pt\theta_c^n$ in \eqref{erreqthetac} to find
	\begin{align}\label{h1ctheta1}
		\|\pt\theta_c^n\|^2 + \frac{\kappa}{2} \pt\|\theta_c^n\|_c^2 + \frac{\gamma m_0}{2} \pt\|\theta_c^n\|^2 \le -(\pt\xi_c^n,\pt\theta_c^n) - \gamma m_0(\xi_c^n,\pt\theta_c^n) + E_c(\pt\theta_c^n) + \Lambda_c(\pt\theta_c^n) + F_c(\pt\theta_c^n).
	\end{align}
	Processed similarly to \eqref{h1utheta2}, it is easily bound the following terms as
	\begin{align}\label{h1ctheta2}
		|(\pt\xi_c^n, \pt\theta_c^n)| +  |E_c(\pt\theta_c^n)|
		\le  \frac{C}{\Delta t}\int_{t_{n-1}}^{t_n} \|\xi_{ct}(t)\|^2 \dt + C\Delta t \int_{t_{n-1}}^{t_n} \|c_{tt}(t)\|^2\dt +  \frac{1}{8}\|\pt\theta_c^n\|^2
	\end{align}
	An application of the Cauchy-Schwarz inequality and the Young's inequality in the second term of \eqref{h1ctheta1} yields
	\begin{equation}
		|\gamma m_0(\xi_c^n,\pt\theta_c^n)|\le C \|\xi_c^n\|^2 + \frac{1}{8}\|\pt\theta_c^n\|^2.
	\end{equation}
	One can consider the similar argument \eqref{h1utheta4}-\eqref{h1utheta6} to bound the nonlinear terms as
	\begin{align}\label{h1ctheta5}
		|\Lambda_c(\pt\theta_c^n)| &\le C\left( \left(h^{-2}\|\e_u^n\|_u^2 + \|\bu^n\|_{H^2}^2\right)\|e_c^n\|_c^2 +  \|\e_u^n\|_u^2\|c^n\|_{H^2}^2 \right)  + \frac{1}{8}\|\pt\theta_c^n\|^2.
	\end{align}
	From \eqref{Fc}, a use of the H{\"o}lder inequality and the Young's inequality yields
	\begin{align}\label{h1ctheta7}
		|F_c(\pt\theta_c^n)| \le & C\left(\|e_\rho^{n-1}\|_\rho^2\|e_c^n\|_c^2 + \|\trho^{n-1}\|_{H^2}^2\|e_c^n\|^2 + \| c^n\|_{H^2}^2\left(\|e_\rho^{n-1}\|^2+\|\trho^n-\trho^{n-1}\|^2\right) \right)  + \frac{1}{8}\|\pt\theta_c^n\|^2.
	\end{align}
	Inserting all the estimates from \eqref{h1ctheta2}-\eqref{h1ctheta7} and put it in \eqref{h1ctheta1} to obtain
	\begin{align}\label{h1ctheta8}
		\|\pt\theta_c^n\|^2 + \kappa\pt\|\theta_c^n\|_c^2 +  \gamma m_0 \pt\|\theta_c^n\|^2
		&\lesssim \frac{1}{\Delta t}\int_{t_{n-1}}^{t_n} \|\xi_{ct}(t)\|^2\dt + \Delta t \int_{t_{n-1}}^{t_n} \left(\|c_{tt}(t)\|^2+\|c^n\|_{H^2}^2\|\trho_t(t)\|^2\right)\dt  \nonumber\\ 
		& +  \|\xi_c^n\|^2 + \left(h^{-2}\|\e_u^n\|_u^2 + \|\bu^n\|_{H^2}^2\right)\|e_c^n\|_c^2 +  \|\e_u^n\|_u^2\|c^n\|_{H^2}^2 \nonumber\\
		& + \| c^n\|_{H^2}^2\|e_\rho^{n-1}\|^2 + \left(\|e_\rho^{n-1}\|_\rho^2 + \|\trho^{n-1}\|_{H^2}^2\right)\|e_c^n\|^2_c.
	\end{align}
	We multiply the above inequality by $\Delta t$ and take summation from $n=1$ to $N$. 
	Then, a use of Theorem \ref{thm:l2u} and remark \ref{rem:l2l2u} with the regularity assumption (\textbf{A1})  completes the rest of the proof.
\end{proof}

\begin{lemma}\label{lem:h1rho}
	Under the assumption of Lemma \ref{lem:thetau}, the following inequality holds true
	\begin{align*}
		\mu\|\theta_\rho^N\|_\rho^2  +  \Delta t\sum_{n=1}^N \|\pt\theta_\rho^n\|^2
		&  
		\lesssim \left(h^{2k} + (\Delta t)^2\right).
	\end{align*}
\end{lemma}

\begin{proof}
	Choose $\chi_h=\pt\theta_\rho^n$ in \eqref{erreqthetarho} to obtain
	\begin{equation}\label{h1rhotheta1}
		\|\pt\theta_\rho^n\|^2 + \frac{\mu}{2}\pt\|\theta_\rho^n\|_\rho^2 \le -(\pt\xi_\rho^n, \pt\theta_\rho^n) + E_\rho(\pt\theta_\rho^n) + \Lambda_\rho(\pt\theta_\rho^n) + \beta G_\rho(\pt\theta_\rho^n).
	\end{equation}
	First three terms on the right hand side of \eqref{h1rhotheta1} can be bounded exactly similar to  \eqref{h1ctheta2} and \eqref{h1ctheta5} as
	\begin{align}\label{h1rhotheta2}
		|(\pt\xi_\rho^n, \pt\theta_\rho^n)| + |E_\rho(\pt\theta_\rho^n)|  +  |\Lambda_\rho(\pt\theta_\rho^n)|
		\le  C\bigg( \frac{1}{\Delta t}\int_{t_{n-1}}^{t_n} \|\xi_{\rho t}(t)\|^2 \dt + \Delta t \int_{t_{n-1}}^{t_n} \|\rho_{tt}(t)\|^2\dt \nonumber\\
		+  (h^{-2}\|\e_u^n\|_u^2+\|\bu^n\|_{H^2})\|e_\rho^n\|^2 + \|\e_u^n\|^2\|\trho^n\|_{H^2}^2\bigg) + \frac{1}{8}\|\pt\theta_\rho^n\|^2.
	\end{align}
	Proceed in a similar fashion of Lemma \ref{lem:thetarho}, one can use of the Holder inequality, the Young's inequality with the inverse inequality, Lemmas \ref{dised2el}, \ref{disl42l2} and \ref{edge2elem} to bound all the terms of $G_\rho(\pt\theta_\rho^n)$ as
	\begin{align}\label{h1grho}
		G_\rho(\pt\theta_\rho^n) 
		& \le C \left((h^{-2}\|^2\|c^n\|_{H^3}^2 + h^{-4}\|\xi_c^n\|_c^2) \|e_\rho^{n-1}\|^2+ (h^{-2}\|\trho^{n-1}\|_{H^2}^2 + h^{-4}\|e_\rho^{n-1}\|^2) \|\theta_c^n\|_c^2\right) \nonumber\\
		& \qquad + C\left(\Delta t \int_{t_{n-1}}^{t_n}\|\nabla\trho_t(t)\|^2 \dt\right)\|c^n\|_{H^3}^2 +  \frac{1}{8}\|\pt\theta_\rho^n\|^2.
	\end{align}
	Inserting all the estimates from \eqref{h1rhotheta2}-\eqref{h1grho} in \eqref{h1rhotheta1} and multiply the resulting inequality by $\Delta t$ and take summation from $n=1$ to $N$. A use of Lemmas \ref{lem:thetac}, \ref{lem:ritz} and \ref{lem:stokes} with regularity assumptions (\textbf{A1}) complete the rest of the proof.
\end{proof}

Combining the results from Lemmas \ref{lem:h1u}, \ref{lem:h1c} and \ref{lem:h1rho}, we can prove the following theorem.
\begin{theorem}\label{thm:h1u}
	Let the assumption (\textbf{A1}) hold true. Suppose $(\bu,c,\trho)$ and $(\bu_h^n,c_h^n,\trho_h^n)$ be the solutions of the continuous system \eqref{dg1rho}-\eqref{dg4div} and the fully discrete system \eqref{fddg1rho}-\eqref{fddg4div}, respectively. Then, with $\e_u^n = \bu^n - \bu_h^n, ~
	e_\rho^n = \trho^n - \trho_h^n, ~ e_c^n = c^n - c_h^n$, the following holds
	\begin{align}
		\|\e_u^n\|_u + \|e_\rho^n\|_\rho + \|e_c^n\|_c \lesssim \left(h^{k}+\Delta t\right).
	\end{align}
\end{theorem}

%%%%%%%%%%%%%%%%%%%%%%% Pressure %%%%%%%%%%%%%%%%%%%%%%

\subsection{Error estimate for fluid pressure}
In this section, we derive optimal $L^2$-error bound for the fluid pressure.

\begin{theorem}\label{thm:l2p}
	Let the assumption (\textbf{A1}) hold true. Suppose $(\bu,c,\trho)$ and $(\bu_h^n,c_h^n,\trho_h^n)$ be the solutions of the continuous system \eqref{dg1rho}-\eqref{dg4div} and the fully discrete system \eqref{fddg1rho}-\eqref{fddg4div}, respectively. Then, with $ e_p^n = p^n - p_h^n$, the following holds
	\begin{align}
		\left(\sum_{n=1}^N\|e_p^n\|^2\right)^{\frac{1}{2}} \lesssim \left(h^{k}+\Delta t\right).
	\end{align}
\end{theorem}
\begin{proof}
	From \eqref{dg3u} and \eqref{fddg3u}, one can find the following error equation:
	\begin{equation}\label{eep1}
		d(\bv_h,e_p^n) =(\pt\e_u^n, \bv_h) + \nu a_u(\e_u^n,\bv_h)   + E_u(\bv_h) + \Lambda_u(\bv_h) + F_u(\bv_h), \quad \forall~\bv_h\in \bV_h,
	\end{equation}
	where $E_u, \Lambda_u$ and $F_u$ are defined on \eqref{Eu}-\eqref{Fu}. A use of the Cauchy-Schwarz inequality and the continuity property of the bilinear form $a_u(\cdot,\cdot)$ yields
	\begin{align}\label{eep2}
		|(\pt\e_u^n, \bv_h) + \nu a_u(\e_u^n,\bv_h)| \lesssim \left(\|\pt\e_u^n\| + \|\e_u^n\|_u\right) \|\bv_h\|_u.
	\end{align}
	Arguing similar to Lemma \ref{lem:thetau}, the following terms estimates as
	\begin{align}\label{eep3}
		|E_u(\bv_h)|+ |F_u(\bv_h)| \lesssim \left(\Delta t\int_{t_{n-1}}^{t_n}\|\bu_{tt}(t)\|^2 \dt\right)^{\frac{1}{2}} \|\bv_h\|_u +  \left(\|e_\rho^{n-1}\|+\left(\Delta t\int_{t_{n-1}}^{t_n}\|\trho_{t}(t)\|^2 \dt\right)^{\frac{1}{2}}\right)\|\nabla \Phi\|_{\bL^4}\|\bv_h\|_u.
	\end{align}
	An application of the boundedness property of the trilinear term $b_2(\cdot,\cdot,\cdot)$ help to bound the following:
	\begin{align}\label{eep5}
		|\Lambda_u(\bv_h)| \lesssim \left(\|\e_u^n\|_u^2+\|\e_u^n\|_u\right)\|\bv_h\|_u.
	\end{align}
	Inserting \eqref{eep2}-\eqref{eep5} in \eqref{eep1} with $e_p^n = \xi_p^n+\theta_p^n$, we obtain
	 \begin{align}\label{eep6}
	 	d(\bv_h,\theta_p^n) \lesssim &\left(\|\xi_p^n\| + \|\pt\e_u^n\| + \|\e_u^n\|_u + \left(\|e_\rho^{n-1}\|+\left(\Delta t\int_{t_{n-1}}^{t_n}\|\trho_{t}(t)\|^2 \dt\right)^{\frac{1}{2}}\right)\|\nabla \Phi\|_{\bL^4}\right.\nonumber\\
	 	& \left.  + \|\e_u^n\|_u^2 + \left(\Delta t\int_{t_{n-1}}^{t_n}\|\bu_{tt}(t)\|^2 \dt\right)^{\frac{1}{2}}\right) \|\bv_h\|_u,
	 \end{align}
	 From the inf-sup condition (Lemma \ref{lem:infsup}) and \eqref{eep6}, we have
	 \begin{align}
	 	\|e_p^n\| &\le \|\xi_p^n\| + \|\theta_p^n\| \le 
	 	\|\xi_p^n\| +  \sup_{\bv_h\in\tilde{\bV}_h} \frac{d(\bv_h,\theta_p^n)}{\|\bv_h\|_{u}} \nonumber\\
	 	& \lesssim  \|\xi_p^n\| + \|\pt\e_u^n\| + \|\e_u^n\|_u + \left(\|e_\rho^{n-1}\|+\left(\Delta t\int_{t_{n-1}}^{t_n}\|\trho_{t}(t)\|^2 \dt\right)^{\frac{1}{2}}\right)\|\nabla \Phi\|_{\bL^4} \nonumber\\
	 	&\qquad  + \|\e_u^n\|_u^2 + \left(\Delta t\int_{t_{n-1}}^{t_n}\|\bu_{tt}(t)\|^2 \dt\right)^{\frac{1}{2}}.
	 \end{align}
	 After taking summation from $n=1$ to $N$ with Theorem \ref{thm:l2u} and \ref{thm:h1u} and Remark \ref{rem:l2l2u}, we complete the rest of the proof.  
\end{proof}

\begin{remark}
	The entire analysis can be extended to the 3D case. To prove the existence and uniqueness of the discrete solution, one needs the smallness assumptions on the given data. The erorr analysis varies due to the application of the trace inequality, inverse hypothesis, Sobolev inequality, etc. The main differences in the 3D case will arise due to the presence of chemotaxic term and the nonlinear (trilinear) terms, which can be handled by using some existing $L^p$ estimate in the 3D case. For example, one can use $L^3$ and $L^6$ estimates instead of $L^4$ estimate to bound the nonlinear terms. 
\end{remark}

\section{Numerical Experiments}

In this section, we demonstrate a few numerical examples for the fully discrete discontinuous Galerkin formulation \eqref{fddg1rho}-\eqref{fddg4div} of \eqref{eqrho}-\eqref{eqint}. We present the first example to verify the convergence rate in both the space and time directions. We consider the second example to show the behavior of the solutions of \eqref{fddg1rho}-\eqref{fddg4div}.

\subsection{Convergence test}
Due to the unavailability of the exact solutions of the system \eqref{eqrho}-\eqref{eqint}, we consider the following modified system with given source terms (in such a way that we can find the exact solution) to verify the convergence rate with respect to space and time variables:\\
\textbf{Step I:} Seek $(\trho_h^n, c_h^n, \bu_h^n, p_h^n)\in X_h^0 \times X_h \times \bV_h \times M_h$ such that for all $n>0$
\begin{align}
	&(\partial_t\trho_{h}^n,\chi_h) + \mu a_\rho(\trho_h^n,\chi_h) + b_1(\bu_h^{n},\trho_h^n,\chi_h) - \beta g(\trho_h^{n-1}+m_0, c_h^{n}, \chi_h) = (f_\rho,\chi_h), \label{ffddg1rho}\\
	&(\partial_t c_{h}^{n},\psi_h) + \kappa a_c(c_h^n,\psi_h) + b_1(\bu_h^{n}, c_h^{n}, \psi_h) + \gamma((\trho_h^{n-1}+m_0)c_h^{n},\psi_h) = (f_c,\psi_h), \label{ffddg2fc}\\
	&(\partial_t\bu_{h}^{n},\bv_h) + \nu a_u(\bu_h^{n},\bv_h) + b_2(\bu_h^{n},\bu_h^{n},\bv_h) - d(\bv_h,p_h^{n}) = ((\trho_h^{n-1}+m_0)\nabla\Phi^{n},\bv_h) + (\bm{f}_u,\bv_h), \label{ffddg3u}\\
	&d(\bu_h^{n},q_h)= 0,\quad\forall~ (\chi_h,\psi_h,\bv_h,q_h)\in X_h\times X_h \times \bH_h \times M_h. \label{ffddg4div}
\end{align}
\textbf{Step II:} Set $\rho_h^n=\trho_h^n+m_0$.

Now, we consider the following two examples with a known solution.

\begin{example}\label{exm1}
For the experiment in 2D, we choose $\f_u, f_c, f_\rho$ in such a way that the solution of \eqref{ffddg1rho}-\eqref{ffddg4div} becomes as follow:
\begin{align*}
u_1(x,y,t) &= -e^{-t}\left( \cos(2\pi x)\sin(2\pi y) - \sin(2\pi y)\right),\\
u_2(x,y,t) &=  e^{-t}\left( \sin(2\pi x)\cos(2\pi y) - \sin(2\pi x)\right),\\
p(x,y,t)   &= e^{-t}\left( \cos(2\pi x) + \sin(2\pi y)\right), \\
\rho(x,y,t) &= e^{-t}\left( \cos(2\pi x)+ \cos(2\pi y) + 3\right),\\
c(x,y,t) &= e^{-t}\left( \cos(2\pi x) + \sin(2\pi y) - 2\pi y + 9 \right).
\end{align*}
\end{example}

For this example, we take the computational domain $\Omega=[0,1]\times[0,1]$. The discrete velocity space $\bH_h$ and the discrete pressure space $L_h$ are constructed using stable discontinuous $(\bm{P}_r,P_{r-1}),r=1,2,3$ pairs. On the other hand, the discrete cell density space $X^0_h$ and the discrete concentration space $X_h$ are comprised of discontinuous linear, quadratic, and cubic elements. For the simulation, we set the parameters $\nu=1, \mu=1, \kappa=1, \gamma=1, \beta=1$ and $\Phi(x,y)=x+y$. 

\begin{figure}[h!]
\centering
\begin{subfigure}[b]{0.48\textwidth}
    \includegraphics[scale=0.50]{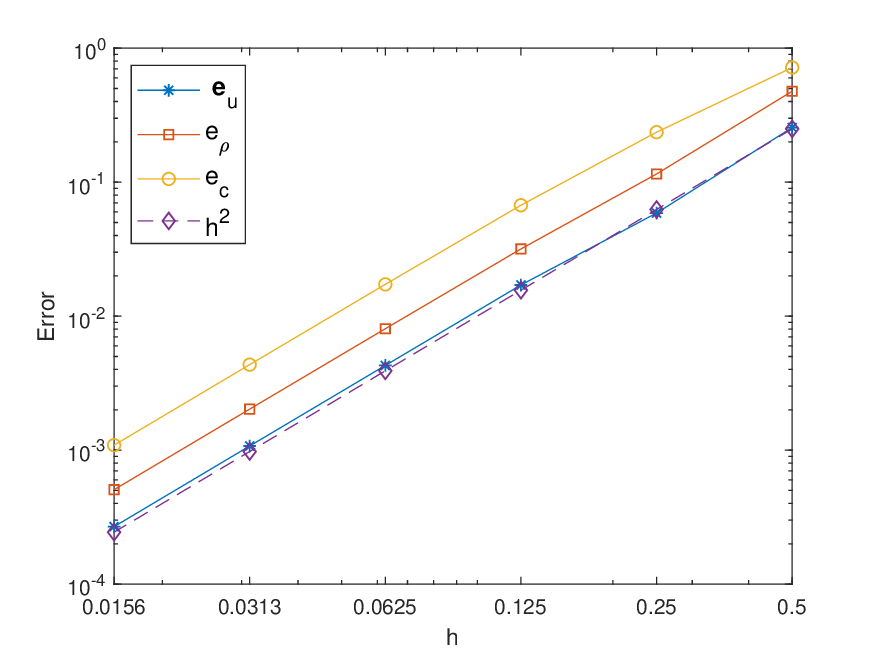}
%    \caption{$L^2$-errors}
%    \label{fig:P1L2ex1}
\end{subfigure}
\hfill
\begin{subfigure}[b]{0.48\textwidth}
    \includegraphics[scale=0.50]{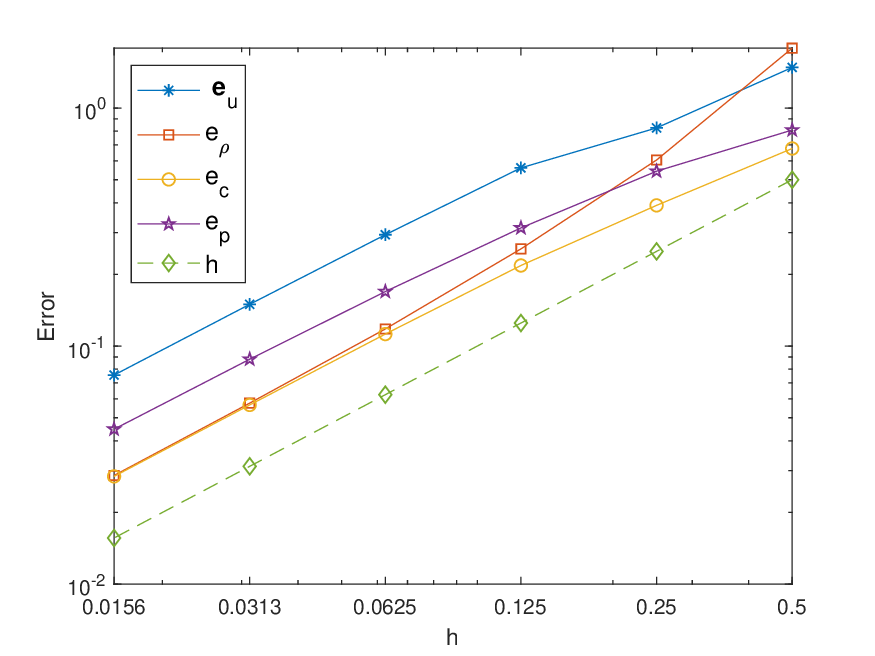}
%    \caption{$H^1$-errors}
%    \label{fig:P1H1ex1}
\end{subfigure}
\caption{Numerical errors with $(\bm{P}_1,P_0,P_1,P_1)$ pairs for Example \ref{exm1}.}
\label{fig:P1ex1}
\end{figure}

\begin{figure}[h!]
\centering
\begin{subfigure}[b]{0.48\textwidth}
    \includegraphics[scale=0.50]{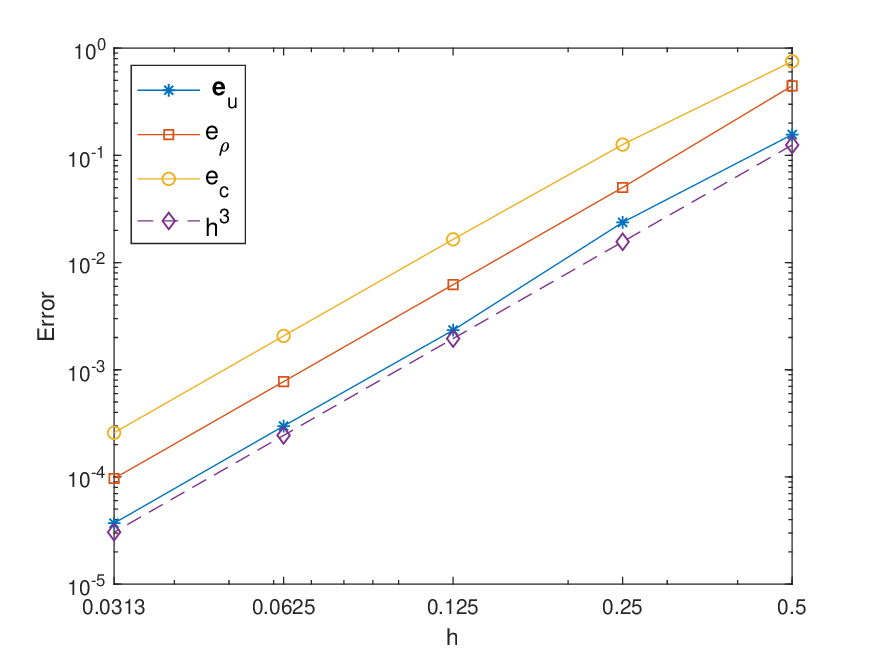}
%    \caption{$L^2$-errors}
%    \label{fig:P2L2ex1}
\end{subfigure}
\hfill
\begin{subfigure}[b]{0.48\textwidth}
    \includegraphics[scale=0.50]{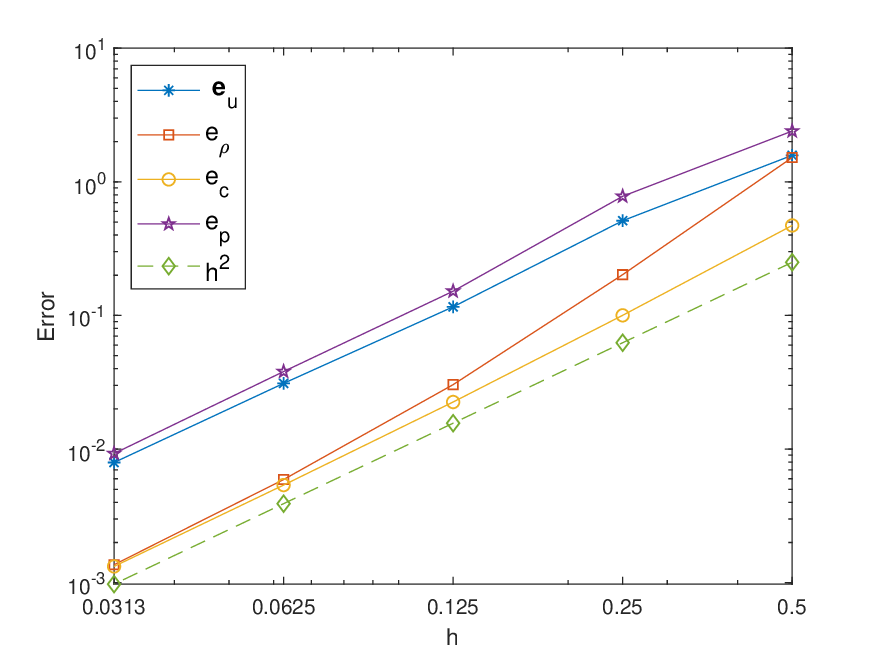}
%    \caption{$H^1$-errors}
%    \label{fig:P2H1ex1}
\end{subfigure}
\caption{Numerical errors with $(\bm{P}_2,P_1,P_2,P_2)$ pairs for Example \ref{exm1}.}
\label{fig:P2ex1}
\end{figure}

\begin{figure}[h!]
\centering
\begin{subfigure}[b]{0.48\textwidth}
    \includegraphics[scale=0.50]{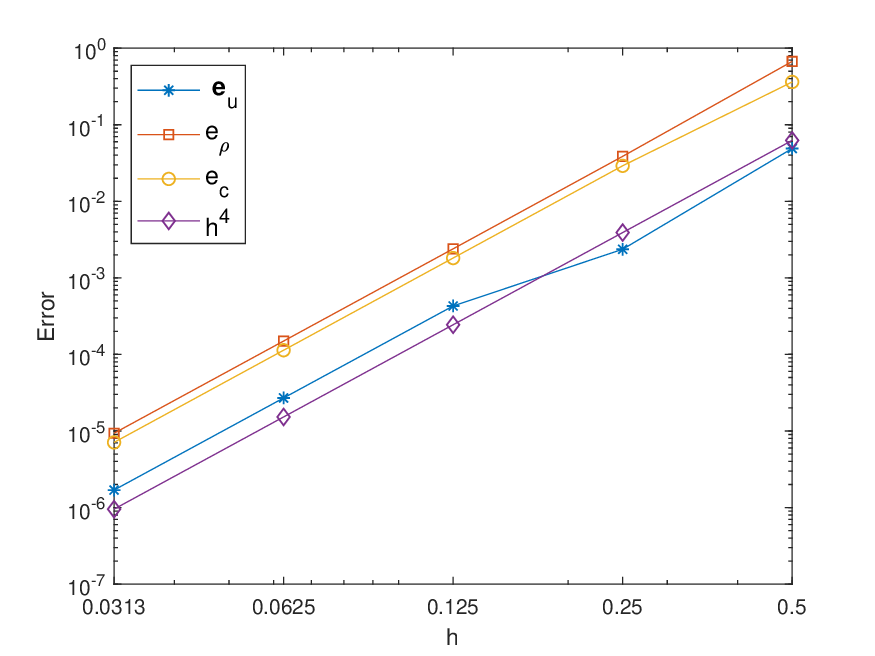}
%    \caption{$L^2$-errors}
%    \label{fig:P2L2ex1}
\end{subfigure}
\hfill
\begin{subfigure}[b]{0.48\textwidth}
    \includegraphics[scale=0.50]{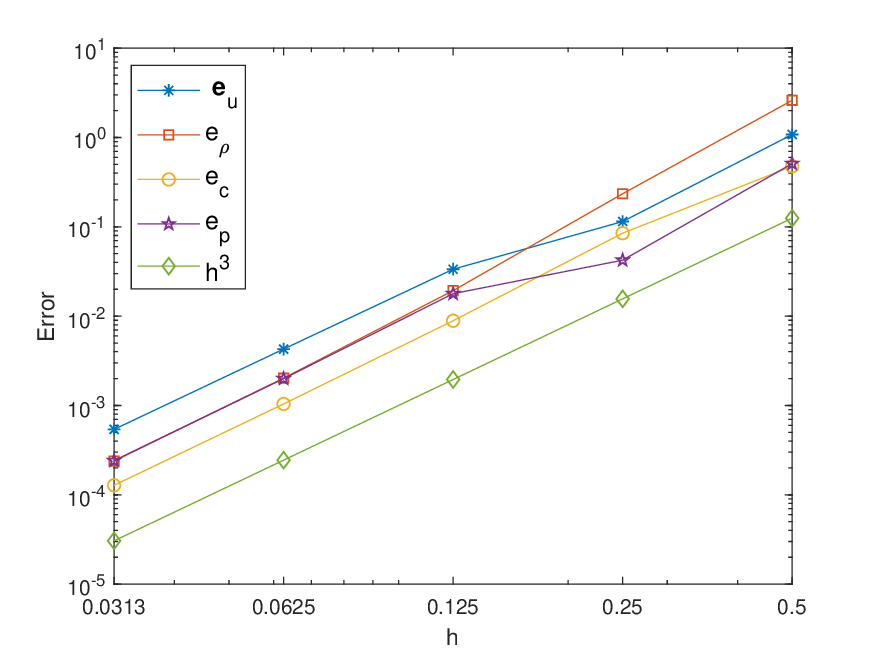}
%    \caption{$H^1$-errors}
%    \label{fig:P2H1ex1}
\end{subfigure}
\caption{Numerical errors with $(\bm{P}_3,P_2,P_3,P_3)$ pairs for Example \ref{exm1}.}
\label{fig:P3ex1}
\end{figure}

First, we present the spacial convergence rates where the domain is discretized based on a regular uniform triangulation with mesh size $h=2^{-i},~i=1,2,\dots,6$. The time discretization parameter $\Delta t$ is chosen as $\mathcal{O}(h^{r+1})$ and the final time is to be taken $T=1$. The Figures \ref{fig:P1ex1}, \ref{fig:P2ex1} and \ref{fig:P3ex1} represent the numerical errors of all the variables with respect to spacial discretizing parameter $h$ for the $(\bm{P}_r,P_{r-1},P_r,P_r)$ pairs $r=1,2,3$, respectively. From these figures, it can be seen that the spatial rates of convergence for the velocity, the cell density, and the concentration are $\mathcal{O}(h^{r+1})$ in $L^2$-norm and  $\mathcal{O}(h^{r})$ in $H^1$-norm. For the pressure, the convergence rate is $\mathcal{O}(h^{r})$ in $L^2$-norm. All these results coincide with the theoretical findings.

% P1dcP0

%P2P1

%P3P2

%%%%%%%%%%%%%%%%%%%%%%%%%%%%%%%%%%%%%%%%%%%%%%%%%%%%%%%%%%%%%%%%%%%%%%%%%%%%%%%%%%%%%%%%%%%%%%%%%%%%%%%%%%%%%

\begin{figure}
\centering
\begin{subfigure}[b]{0.30\textwidth}
    \includegraphics[scale=0.40]{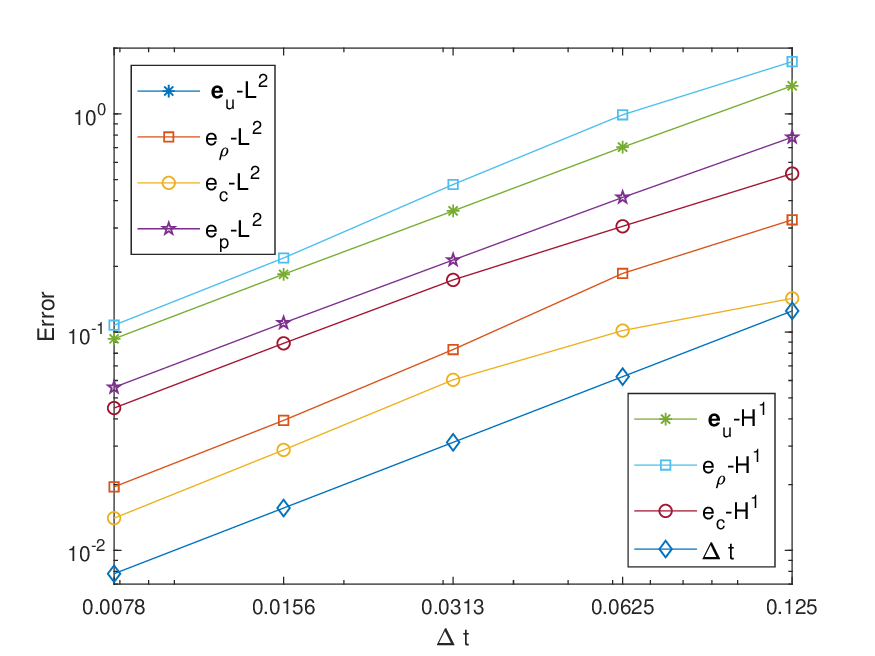}
    \caption{$(\bm{P}_1,P_0,P_1,P_1)$ pairs}
    \label{fig:P1L2ex1time}
\end{subfigure}
\hfill
\begin{subfigure}[b]{0.30\textwidth}
    \includegraphics[scale=0.40]{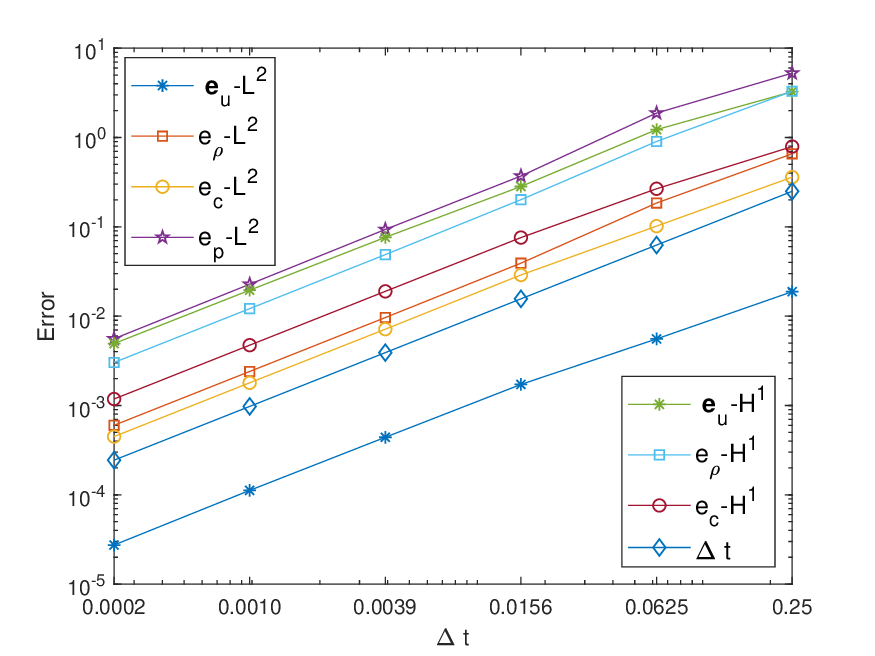}
    \caption{$(\bm{P}_2,P_1,P_2,P_2)$ pairs}
    \label{fig:P1H1ex1time}
\end{subfigure}
\hfill
\begin{subfigure}[b]{0.30\textwidth}
    \includegraphics[scale=0.40]{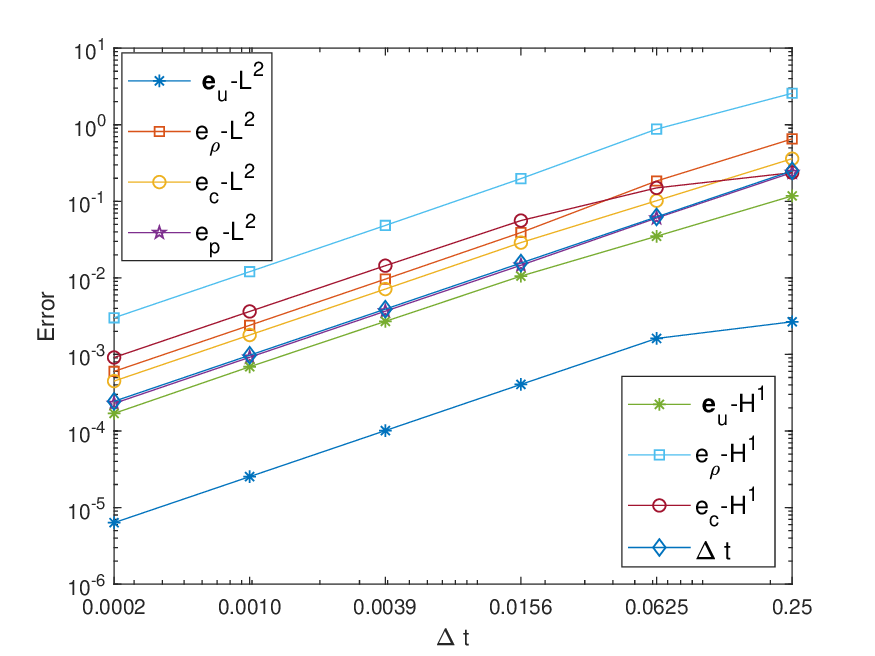}
    \caption{$(\bm{P}_3,P_2,P_3,P_3)$ pairs}
    \label{fig:P2L2ex1time}
\end{subfigure}
\caption{Numerical errors with respect to time variable for Example \ref{exm1}.}
\label{fig:time_ex1}
\end{figure}

To show the convergence rates in temporal direction, we take an uniform partition of the time interval $[0,T]$ with $\Delta t = 2^i,~i=3,4,\dots,7$ for the $(\bm{P}_1,P_{0},P_1,P_1)$ pair and $\Delta t = 2^{2i},~i=1,2,\dots,6$ for the $(\bm{P}_r,P_{r-1},P_r,P_r)$ pairs $r=2,3$. The mesh parameter $h$ has to be fixed as $\mathcal{O}(\Delta t^{\frac{1}{r+1}})$ for $L^2$-norm and $\mathcal{O}(\Delta t^{\frac{1}{r}})$ for $H^1$-norm. In Figure \ref{fig:time_ex1}, we present the numerical errors for the velocity, the cell density, and the concentration in $L^2$ as well as $H^1$-norms and $L^2$-norm for the pressure. From this figure, we conclude that the rate of convergence is linear for all the variables in the temporal direction for the discontinuous finite element pairs $(\bm{P}_r,P_{r-1},P_r,P_r), r=1,2,3$, which are consistent with our theoretical results (since we have applied a first-order scheme in time direction).

% Conv time

\begin{example}\label{exm3d}
For the experiment in 3D, we choose $\f_u, f_c, f_\rho$ in such a way that the solution of \eqref{ffddg1rho}-\eqref{ffddg4div} becomes as follow:

\begin{align*}
u_1(x,y,z,t) &= 10\cos(t)w(x)\left(w_y(y)w(z)-w(y)w_z(z)\right),\\
u_2(x,y,z,t) &=  10\cos(t)w(y)\left(w_z(z)w(x)-w(z)w_x(x)\right),\\
u_3(x,y,z,t) &= 10\cos(t)w(x)\left(w(x)w_y(y)-w_x(x)w_y(y)\right), \\ 
\rho(x,y,z,t) &= 10\cos(t)w(x)w(y)w(z),\\
c(x,y,z,t) &= 10\cos(t)w_x(x)w_y(y)w_z(z),\\
p(x,y,z,t)   &= 10\cos(t)(2x - 1)(2y - 1)(2z - 1),
\end{align*}
where $w(x)=x^2(1-x)^2$ and $w_s(s)$ is the partial derivatives of $w(s)$ with respect to $s$.
\end{example}

\begin{figure}[h!]
\centering
\begin{subfigure}[b]{0.48\textwidth}
    \includegraphics[scale=0.450]{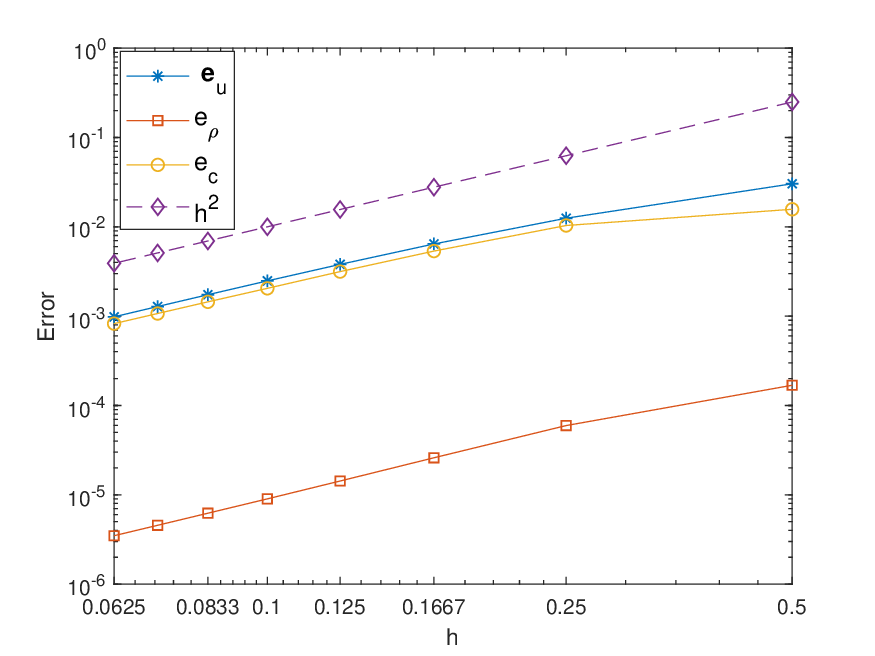}
%    \caption{$L^2$-errors}
%    \label{fig:P1L2ex1}
\end{subfigure}
\hfill
\begin{subfigure}[b]{0.48\textwidth}
    \includegraphics[scale=0.450]{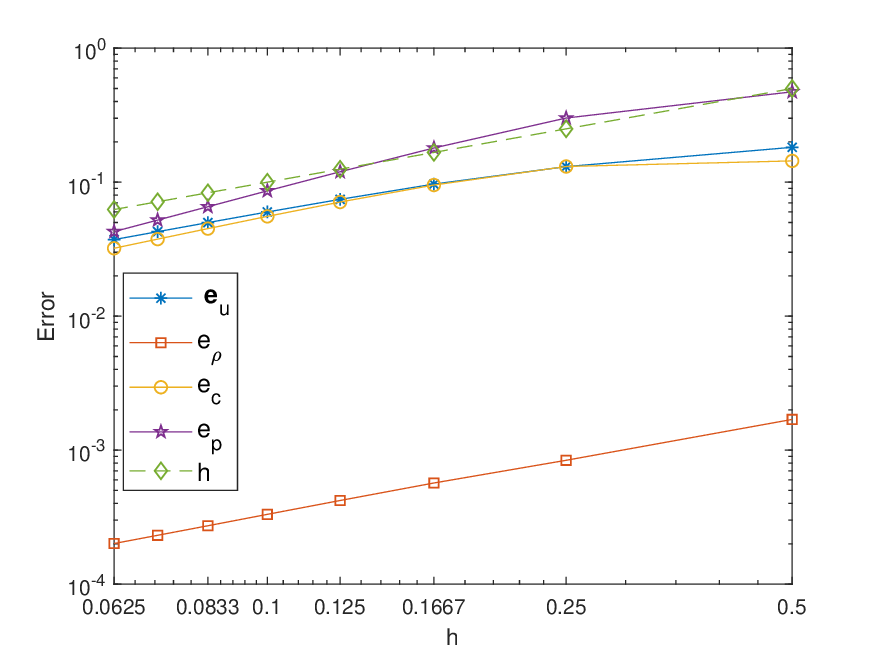}
%    \caption{$H^1$-errors}
%    \label{fig:P1H1ex1}
\end{subfigure}
\caption{Numerical errors with $(\bm{P}_1,P_0,P_1,P_1)$ pairs for Example \ref{exm3d}.}
\label{fig:P1ex3d}
\end{figure}

\begin{figure}[h!]
\centering
\begin{subfigure}[b]{0.48\textwidth}
    \includegraphics[scale=0.450]{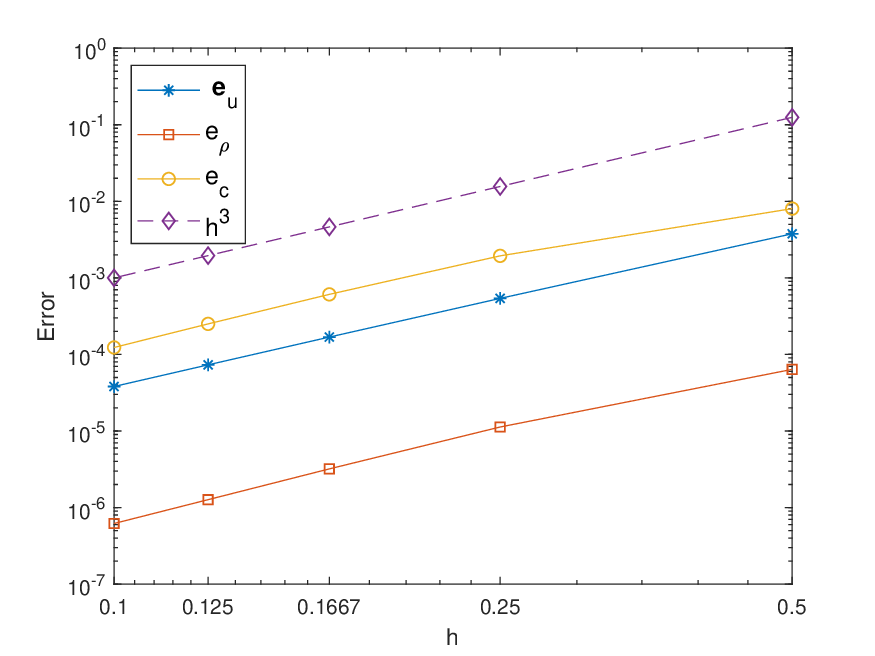}
%    \caption{$L^2$-errors}
%    \label{fig:P2L2ex1}
\end{subfigure}
\hfill
\begin{subfigure}[b]{0.48\textwidth}
    \includegraphics[scale=0.450]{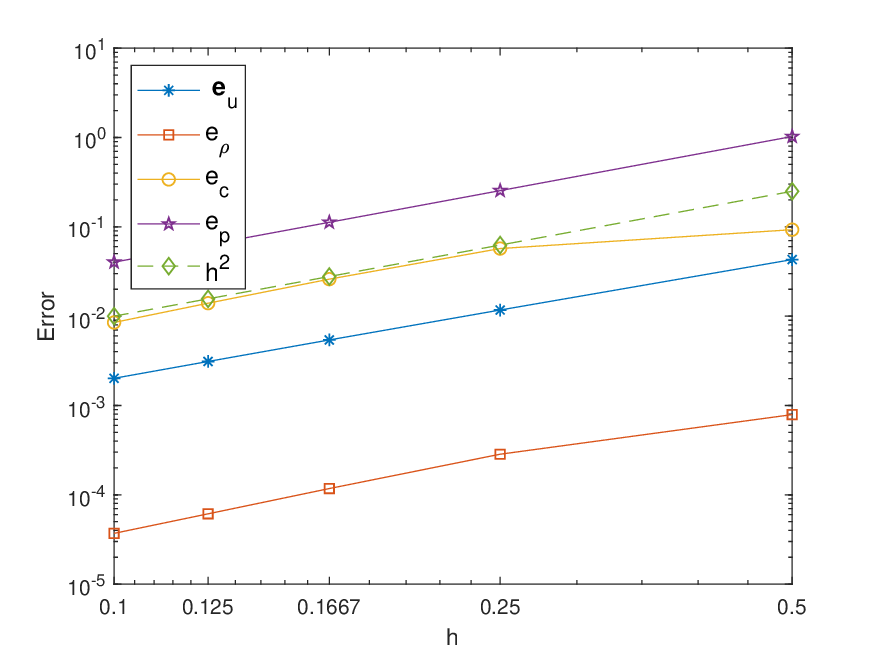}
%    \caption{$H^1$-errors}
%    \label{fig:P2H1ex1}
\end{subfigure}
\caption{Numerical errors with $(\bm{P}_2,P_1,P_2,P_2)$ pairs for Example \ref{exm3d}.}
\label{fig:P2ex3d}
\end{figure}

\begin{figure}
\centering
\begin{subfigure}[b]{0.48\textwidth}
    \includegraphics[scale=0.450]{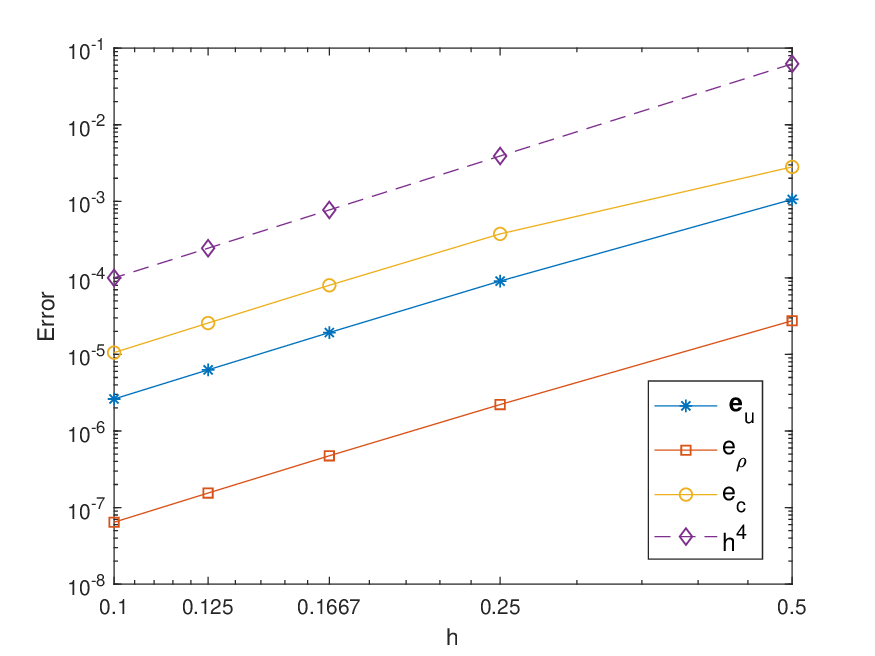}
%    \caption{$L^2$-errors}
%    \label{fig:P2L2ex1}
\end{subfigure}
\hfill
\begin{subfigure}[b]{0.48\textwidth}
    \includegraphics[scale=0.450]{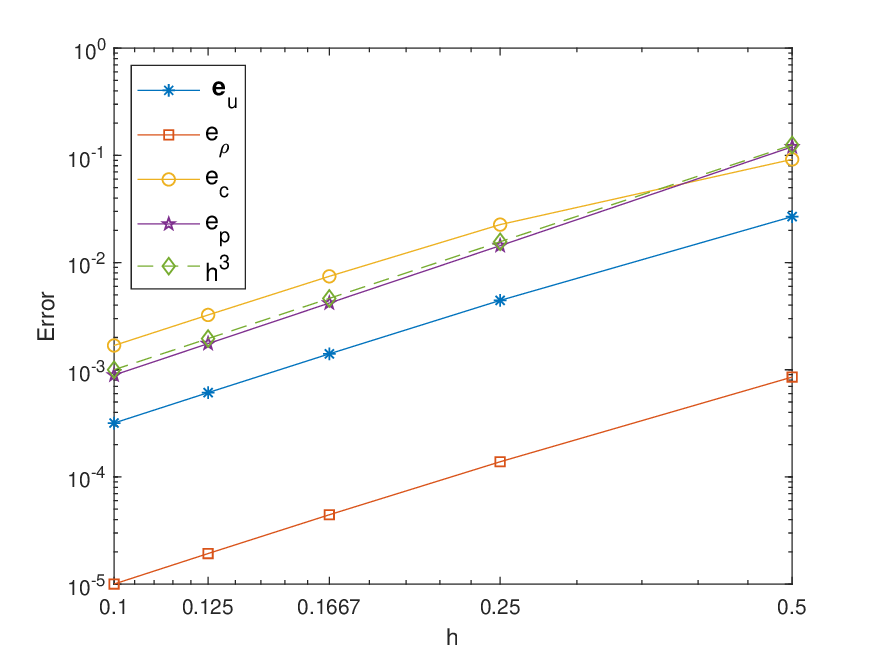}
%    \caption{$H^1$-errors}
%    \label{fig:P2H1ex1}
\end{subfigure}
\caption{Numerical errors with $(\bm{P}_3,P_2,P_3,P_3)$ pairs for Example \ref{exm3d}.}
\label{fig:P3ex3d}
\end{figure}

For this example, we take the computational domain $\Omega=[0,1]^3$. The discrete spaces $(\bH_h, L_h, X^0_h, X_h)$ are constructed using stable discontinuous $(\bm{P}_r,P_{r-1},P_r,P_r),r=1,2,3$ pairs. For the simulation, we set the parameters $\nu=1, \mu=1, \kappa=1, \gamma=1, \beta=1$ and $\Phi(x,y)=x+y+z$. 
The domain is discretized based on a regular uniform triangulation with mesh size $h=1/(2i),~i=1,2,\dots,8$ for $r=1$ and $h=1/(2i),~i=1,2,\dots,5$ for $r=2,3$. The time discretization parameter $\Delta t$ is chosen as $\mathcal{O}(h^{r+1})$ and the final time is to be taken $T=0.1$. The Figures \ref{fig:P1ex3d}, \ref{fig:P2ex3d} and \ref{fig:P3ex3d} represent the numerical errors of all the variables with respect to the spacial discretizing parameter $h$ for the $(\bm{P}_r,P_{r-1},P_r,P_r)$ pairs $r=1,2,3$, respectively. From these figures, it can be seen that the spatial rates of convergence for the velocity, the cell density, and the concentration are $\mathcal{O}(h^{r+1})$ in $L^2$-norm and  $\mathcal{O}(h^{r})$ in $H^1$-norm. For the pressure, the convergence rate is $\mathcal{O}(h^{r})$ in $L^2$-norm. All these results coincide with the theoretical findings.

% P1dcP0 3d

%P2P1 3d

%P3P2 3d

\subsection{Movement of cells guided by the concentration}
We take an example where the numerical experiment is carried out for the system \eqref{fddg1rho}-\eqref{fddg4div} (that is, for no source terms). This example is taken from \cite{DRRV21}.
\begin{example} \label{bench}
For this simulation, we consider a rectangular domain $\Omega = [0,2]\times[0,1]$ and the initial data as:
\begin{align*}
\rho_0(x,y) & = \sum_{i=1}^3 70\exp\left(-8(x-s_i)^2-10(y-1)^2\right),\\
c_0(x,y)   & = 30\exp\left( -5(x-1^2-5(y-0.5)^2\right),\\
\bu_0(x,y) & = \bm{0},
\end{align*}
where $s_1=0.2$, $s_2=0.5$ and $s_3=1.2$.
\end{example}

\begin{figure}
\centering
\begin{multicols}{2}  
	\includegraphics[scale=0.50]{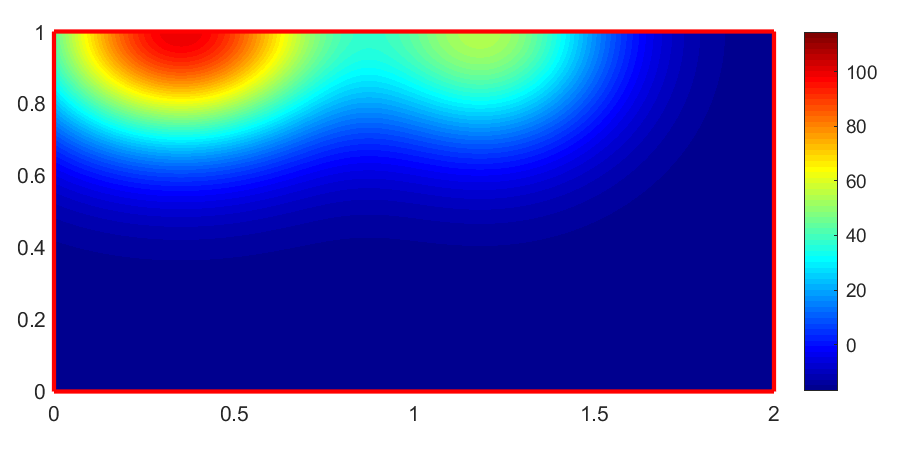}
	\includegraphics[scale=0.50]{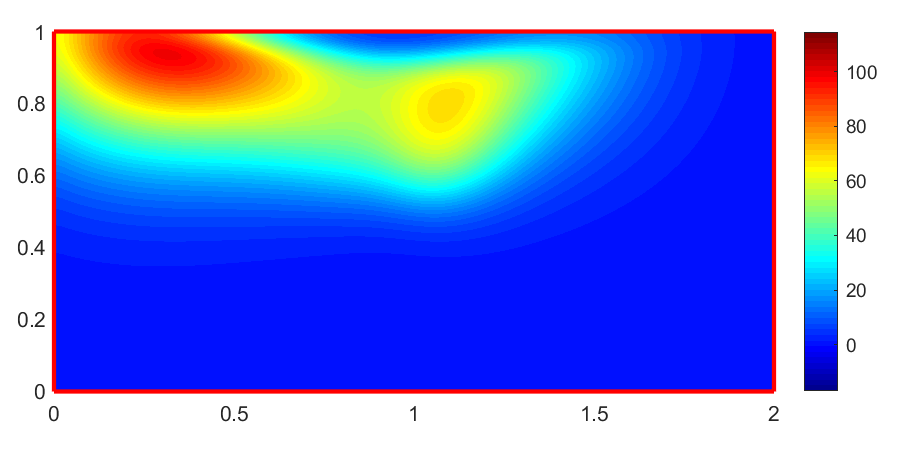}
	\includegraphics[scale=0.50]{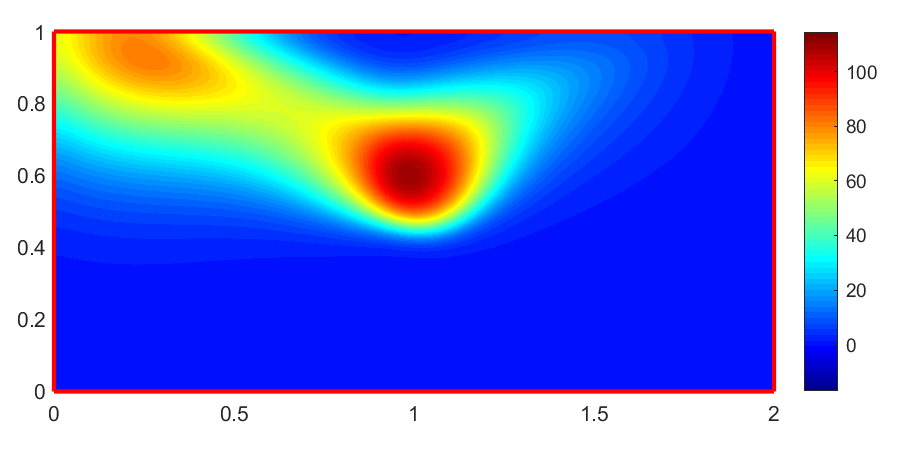}
	\includegraphics[scale=0.50]{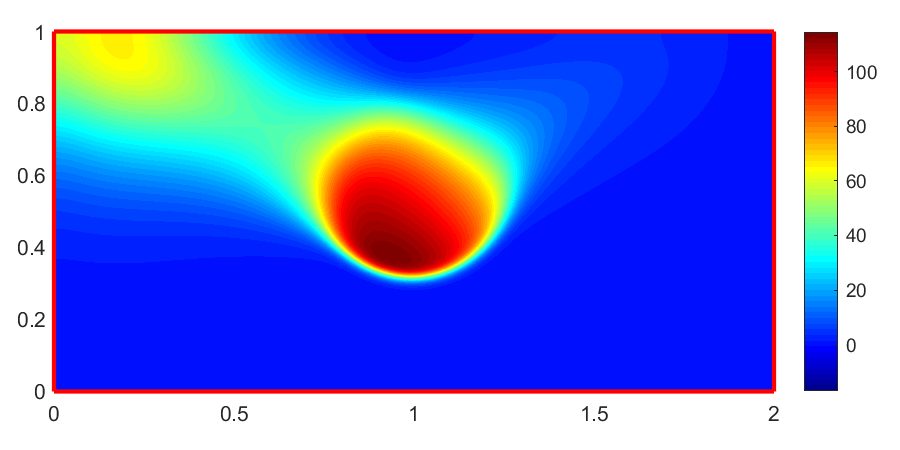}
	\includegraphics[scale=0.50]{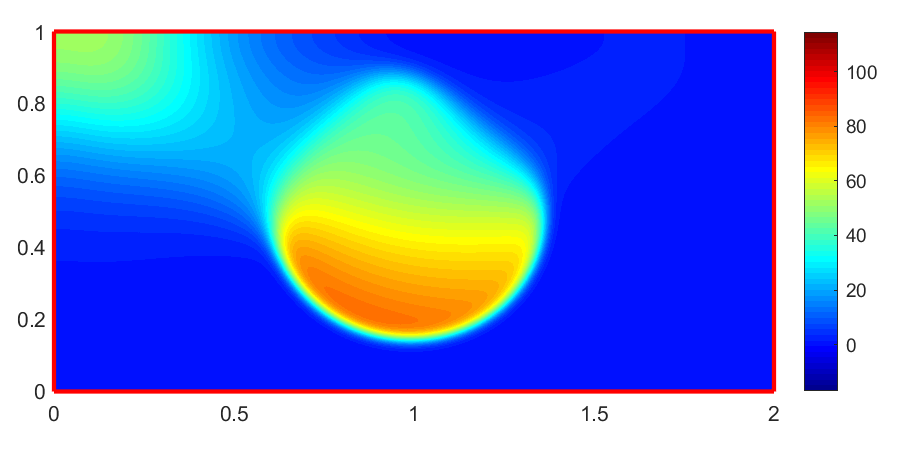}
	\includegraphics[scale=0.50]{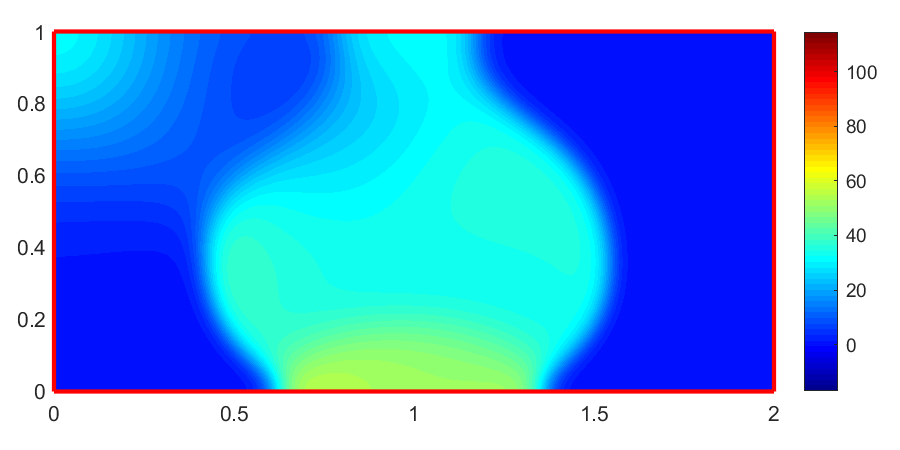}

\columnbreak

	\includegraphics[scale=0.50]{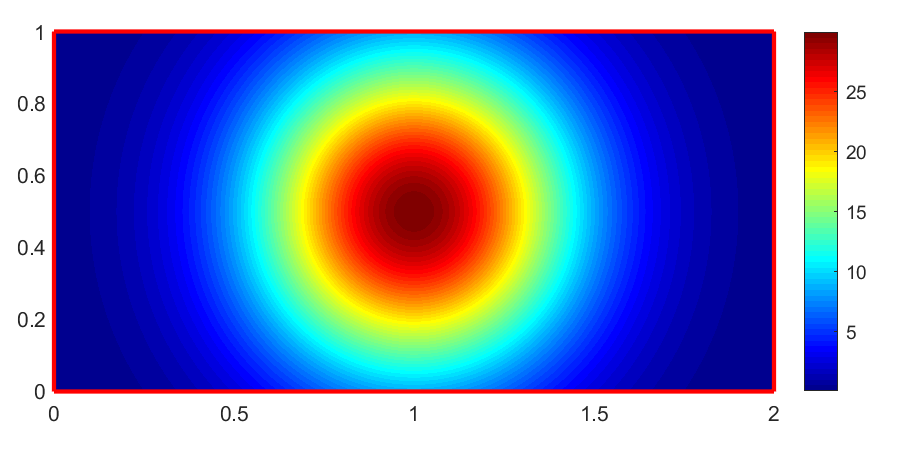}
	\includegraphics[scale=0.50]{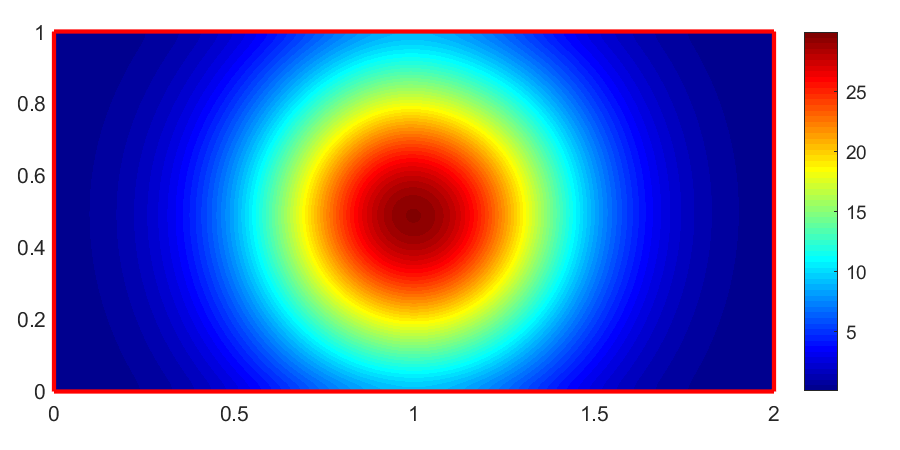}
	\includegraphics[scale=0.50]{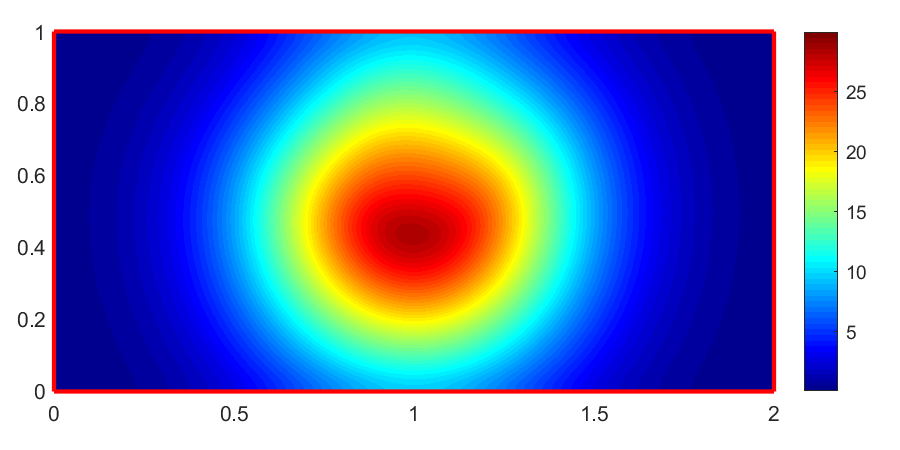}
	\includegraphics[scale=0.50]{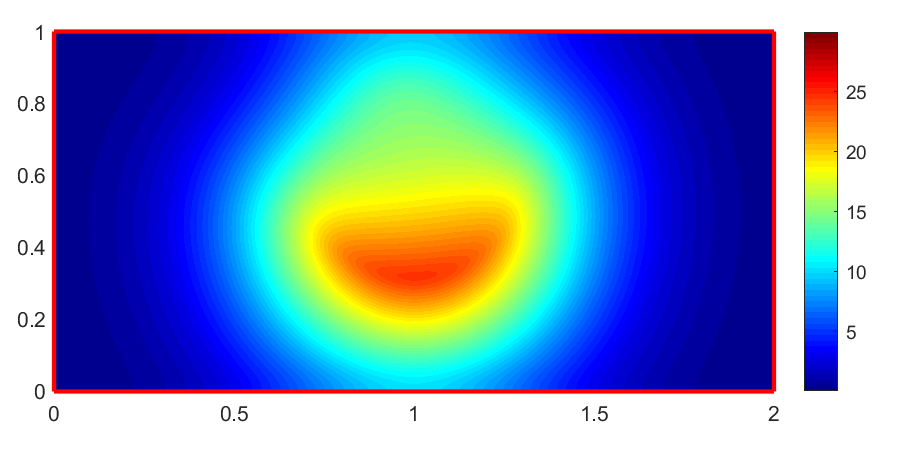}
	\includegraphics[scale=0.50]{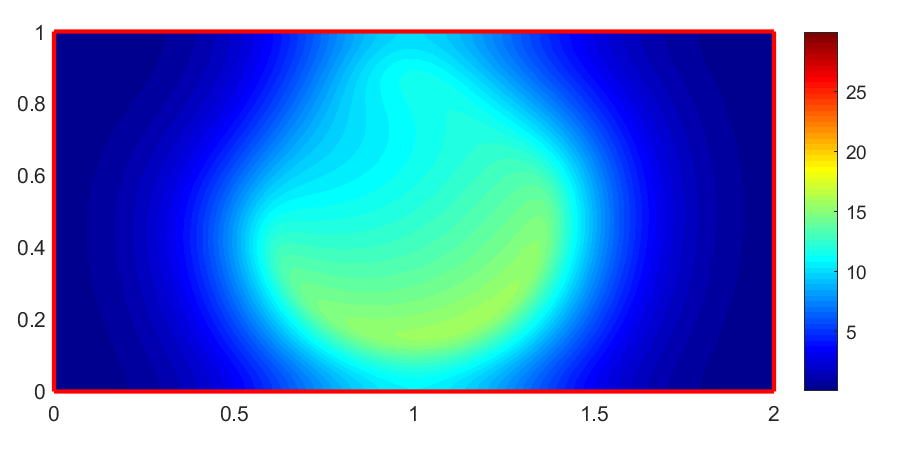}
	\includegraphics[scale=0.50]{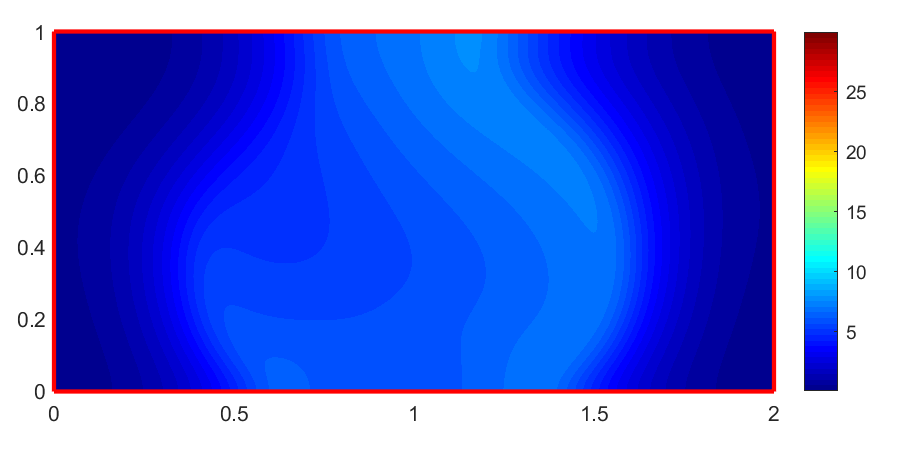}

\end{multicols}
\caption{The evolution of the cell density (left) and the concentration of chemical substances (right) at different time levels $t = 0, 0.0002, 0.0005, 0.001, 0.002$ and $0.005$ (top to down). }
\label{fig:etac}
\end{figure}

\begin{figure}
\centering
\begin{multicols}{2}  
	\includegraphics[scale=0.50]{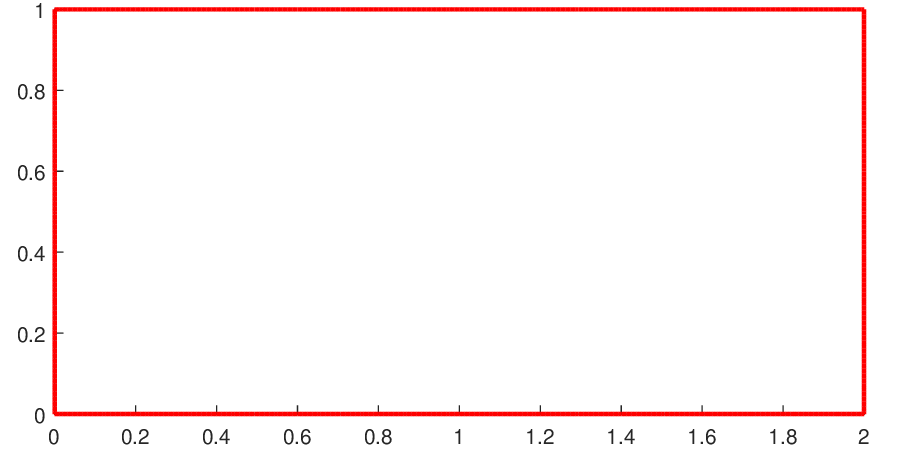}
	\includegraphics[scale=0.50]{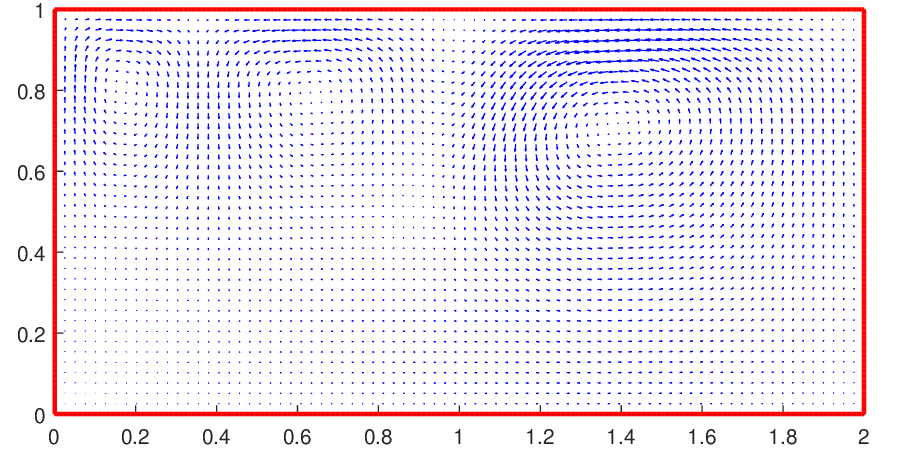}
	\includegraphics[scale=0.50]{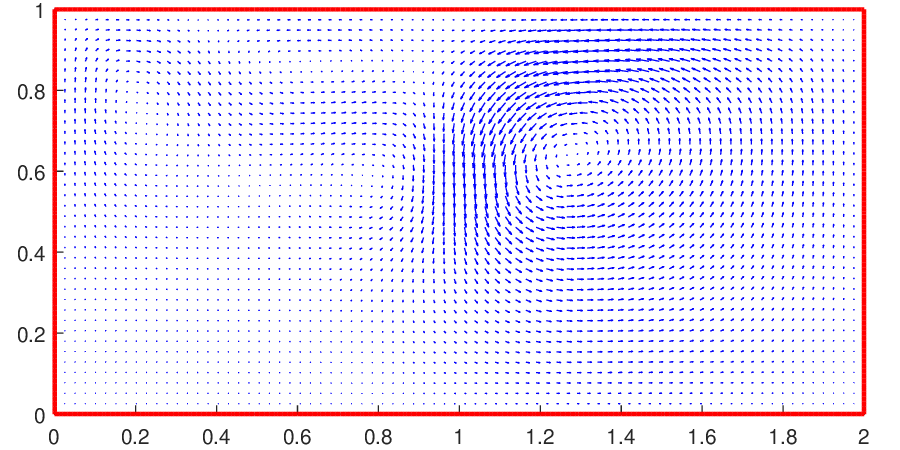}
	\includegraphics[scale=0.50]{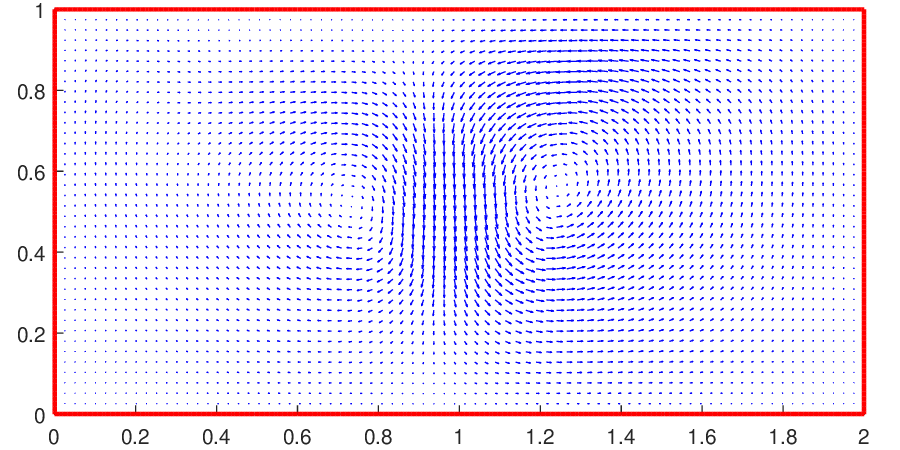}
	\includegraphics[scale=0.50]{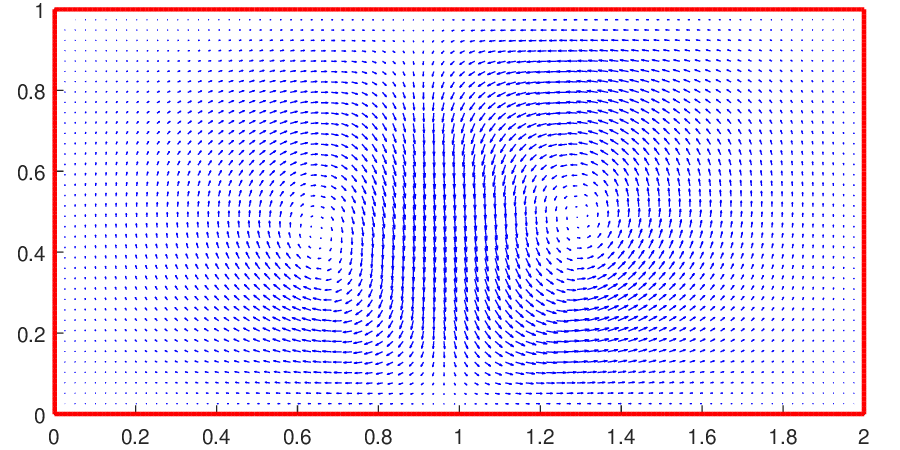}
	\includegraphics[scale=0.50]{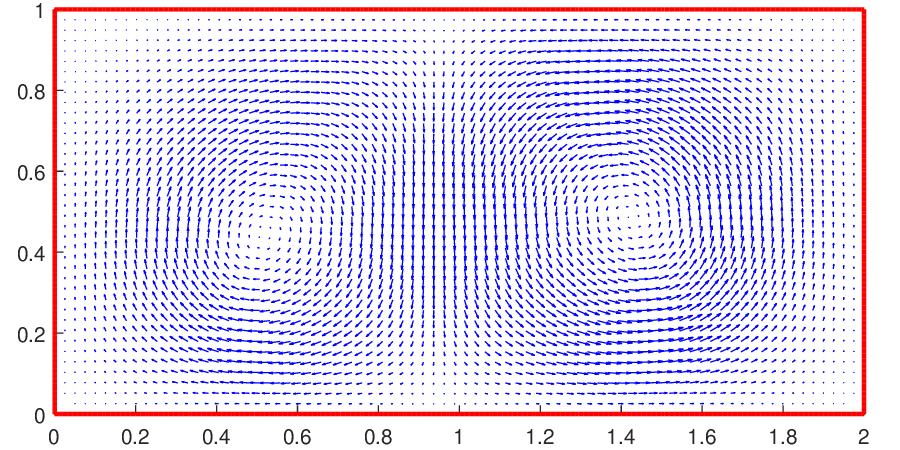}

\columnbreak

	\includegraphics[scale=0.50]{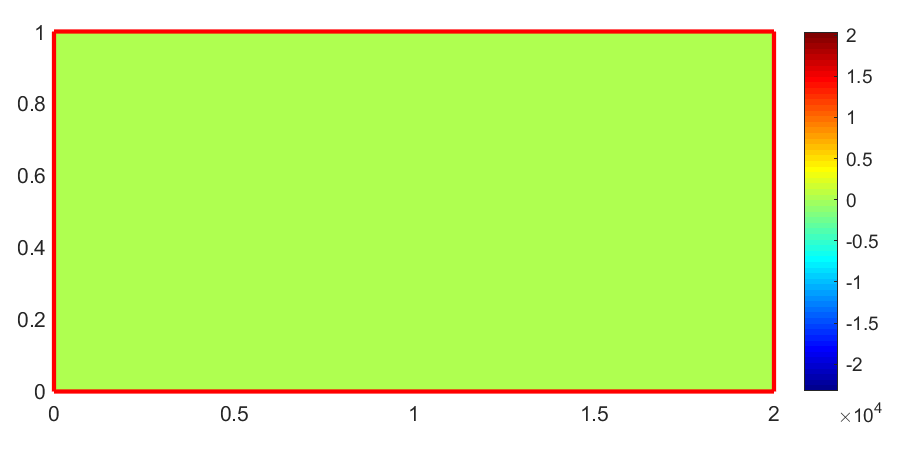}
	\includegraphics[scale=0.50]{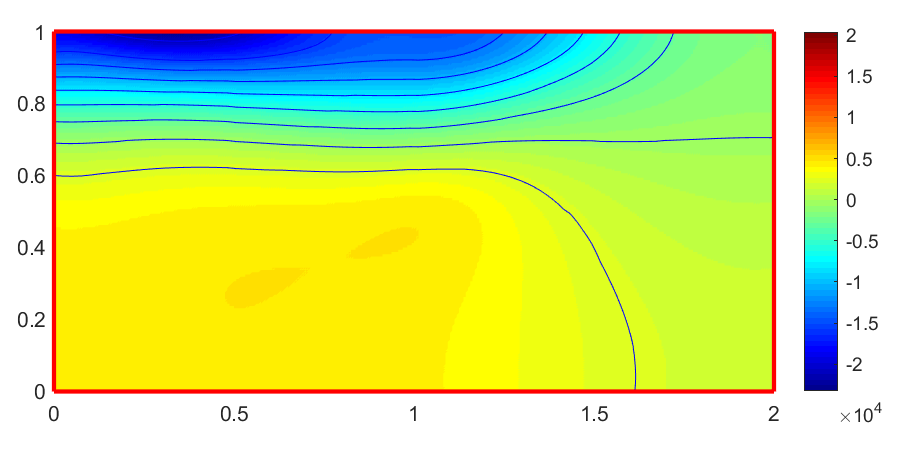}
	\includegraphics[scale=0.50]{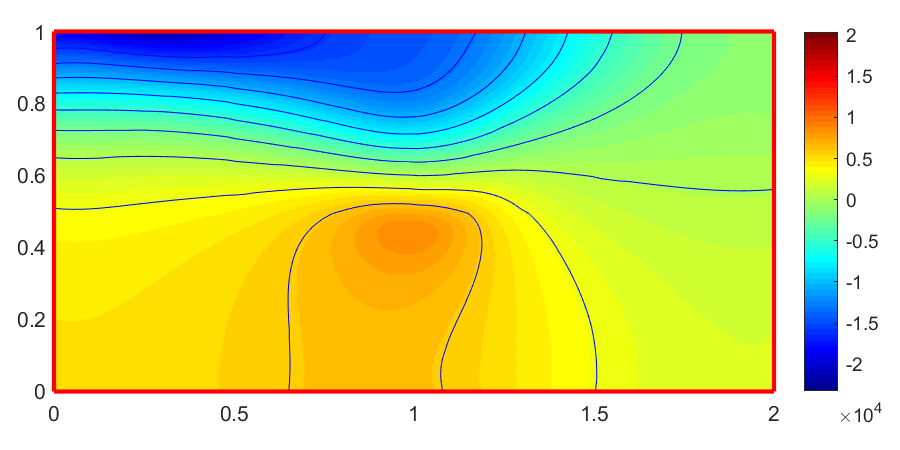}
	\includegraphics[scale=0.50]{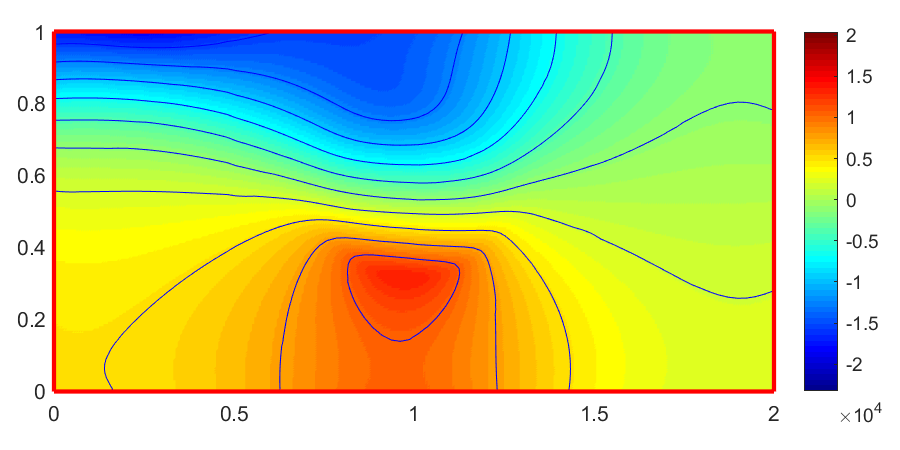}
	\includegraphics[scale=0.50]{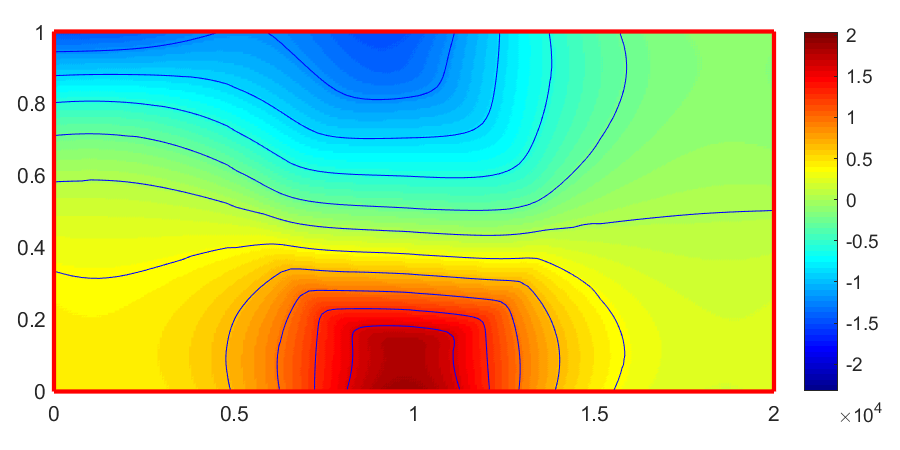}
	\includegraphics[scale=0.50]{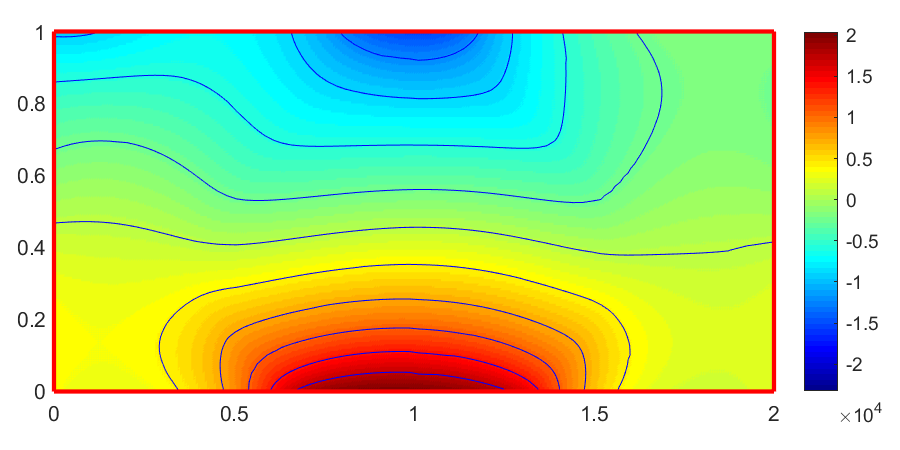}

\end{multicols}
\caption{The velocity vector (left) and the pressure contour (right) of the fluid at different time levels $t = 0, 0.0002, 0.0005, 0.001, 0.002$ and $0.005$ (top to down). }
\label{fig:vecp}
\end{figure}

\begin{figure}
\centering
\begin{multicols}{2}  
	\includegraphics[scale=0.50]{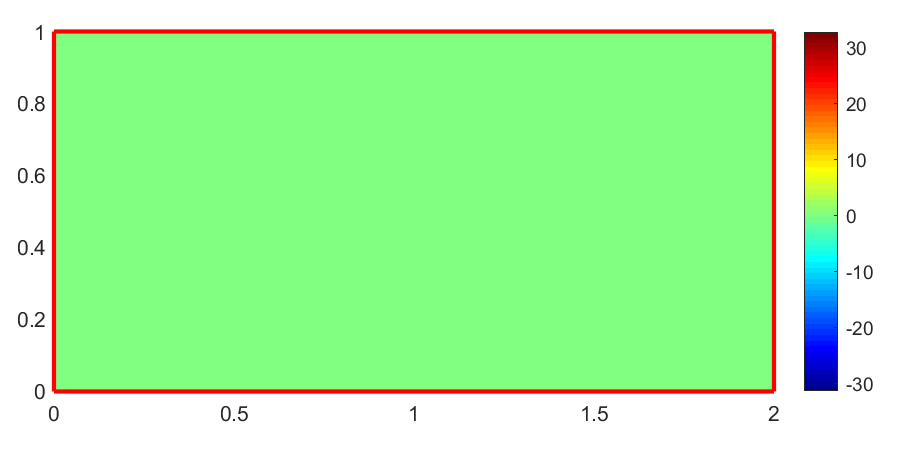}
	\includegraphics[scale=0.50]{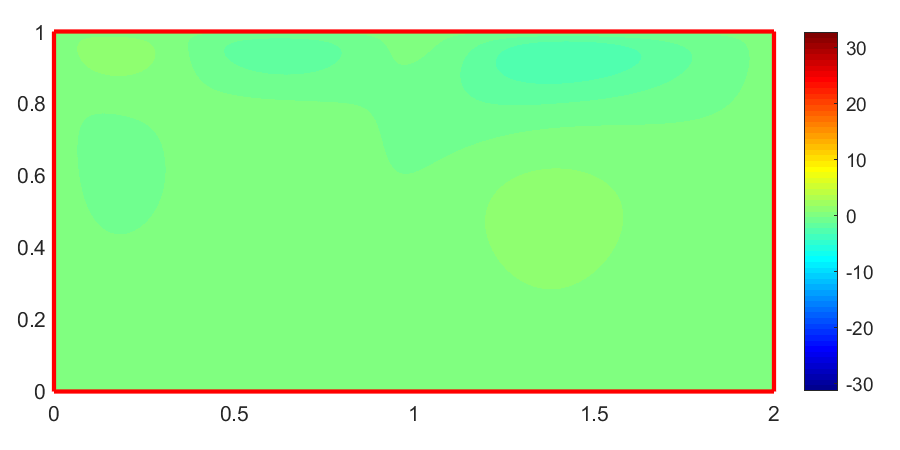}
	\includegraphics[scale=0.50]{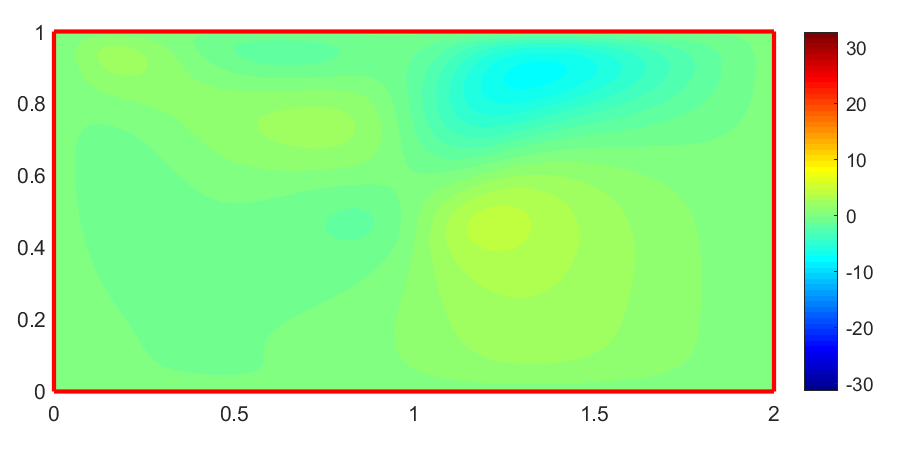}
	\includegraphics[scale=0.50]{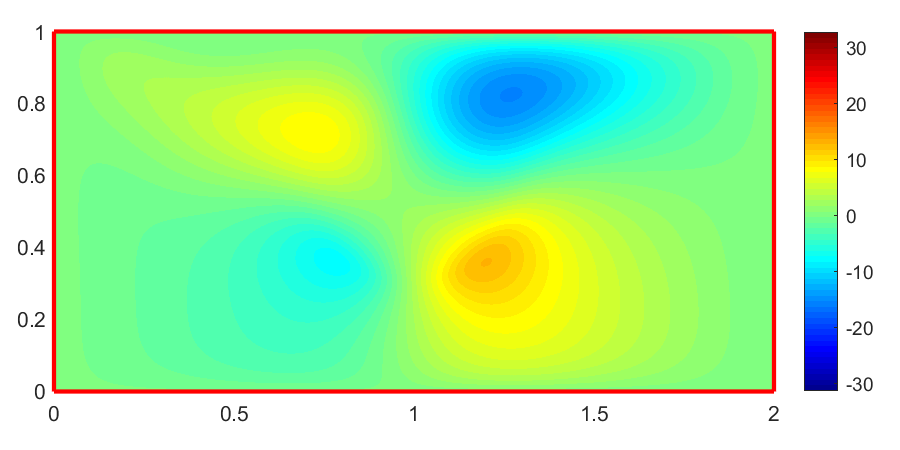}
	\includegraphics[scale=0.50]{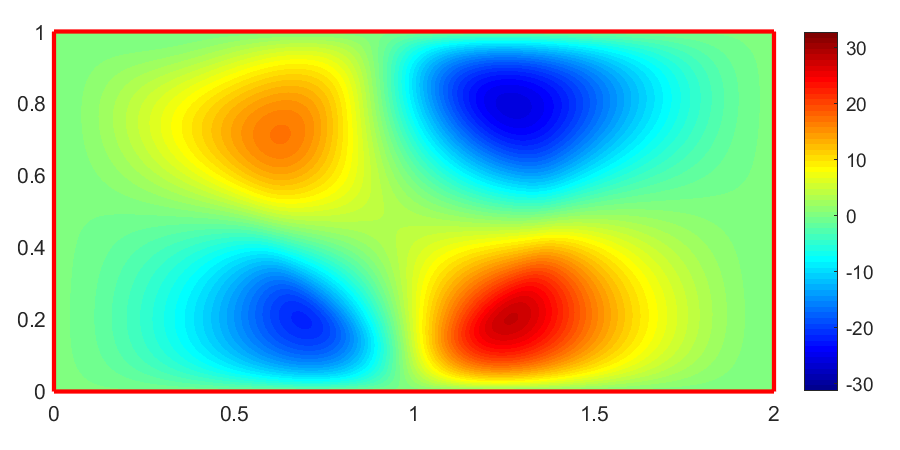}
	\includegraphics[scale=0.50]{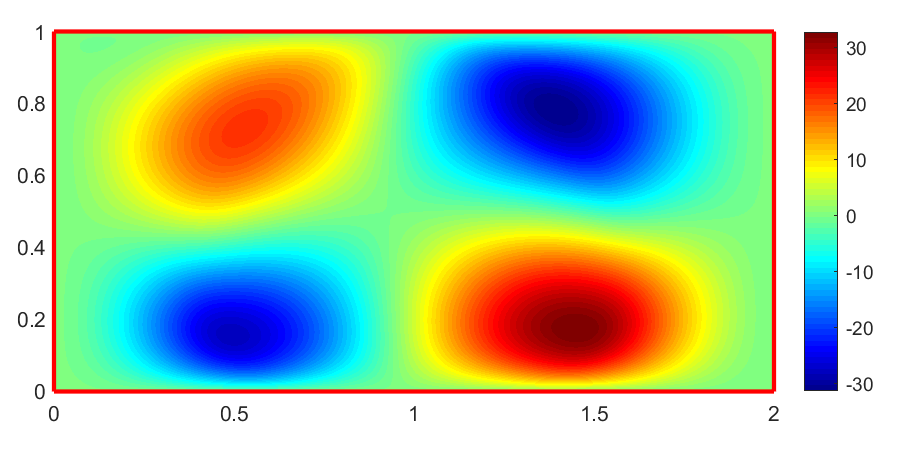}

\columnbreak

	\includegraphics[scale=0.50]{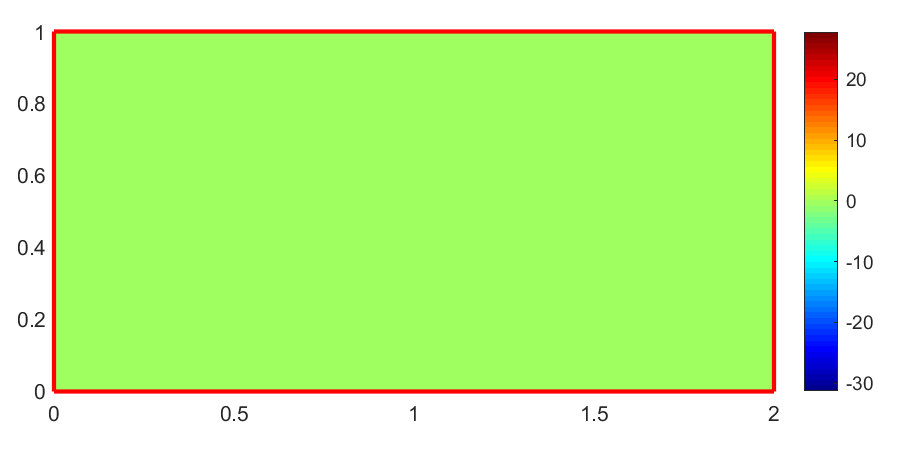}
	\includegraphics[scale=0.50]{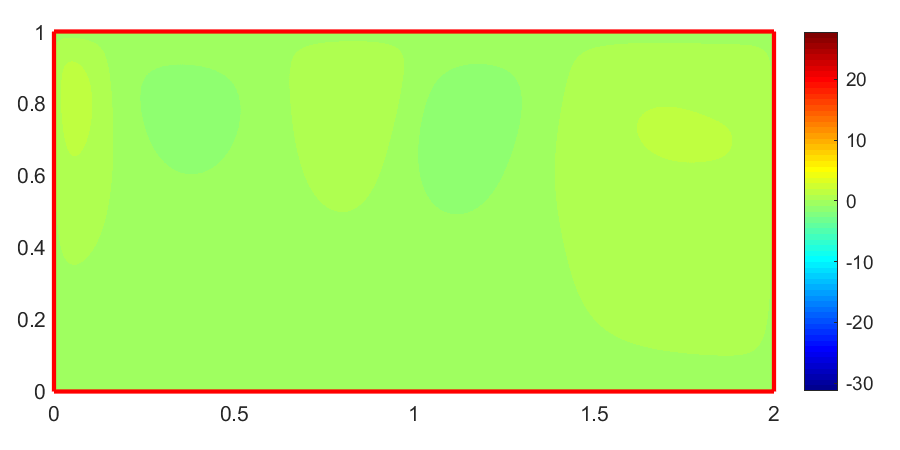}
	\includegraphics[scale=0.50]{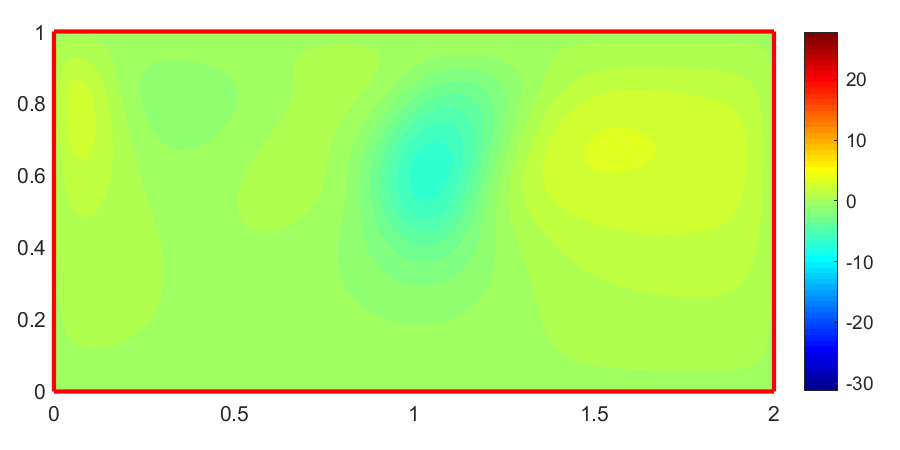}
	\includegraphics[scale=0.50]{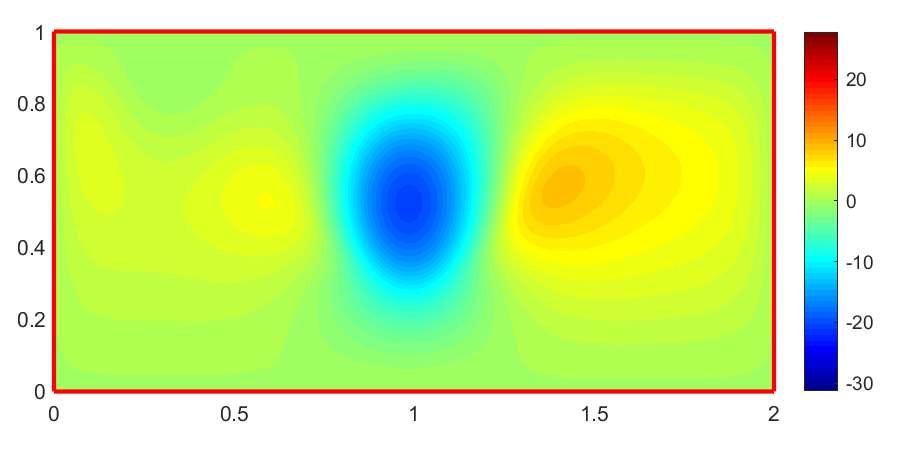}
	\includegraphics[scale=0.50]{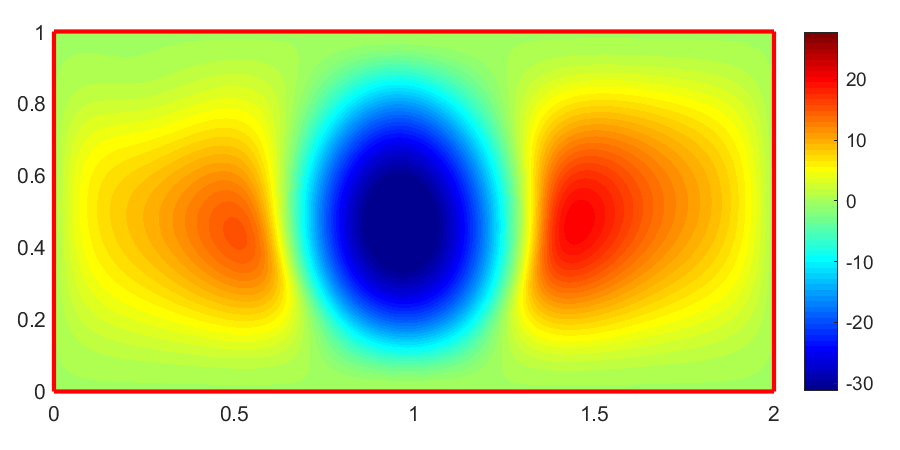}
	\includegraphics[scale=0.50]{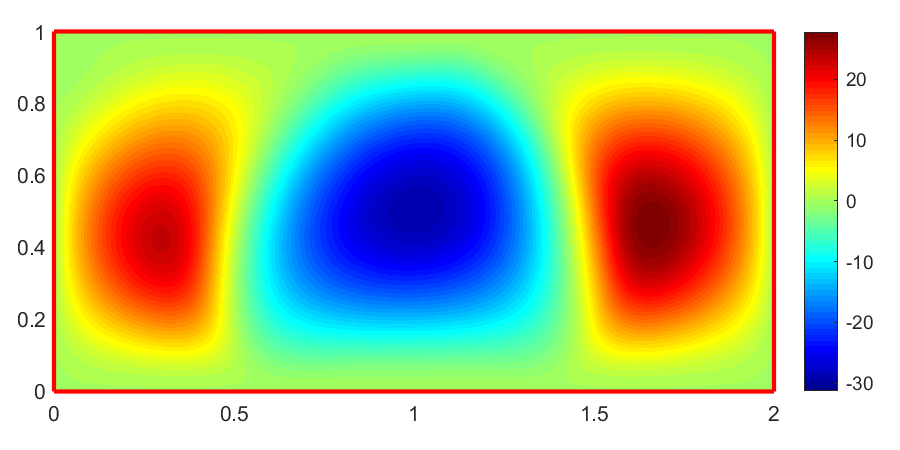}

\end{multicols}
\caption{The evolution of the velocity components (left: primary component, right: secondary component) of the fluid at different time levels $t = 0, 0.0002, 0.0005, 0.001, 0.002$ and $0.005$ (top to down). }
\label{fig:uu}
\end{figure}

For the experiment, we set $\nu=10, \mu=4, \kappa=1, \gamma=6, \beta=8$ and $\Phi(x,y)=-1000y$. We choose mesh size $h=\frac{1}{160}$ and time step $\Delta t = 10^{-5}$. We present the experiment result at time $t = 0, 0.0002, 0.0005, 0.001, 0.002$ and $0.005$. For the simulation, we discretize the domain based on $(\bm{P}_1,P_0,P_1,P_1)$ pairs. Figure \ref{fig:etac} illustrates the progression of cell density and concentration, while Figure \ref{fig:vecp} showcases the velocity vector field and the pressure contour. Additionally, Figure \ref{fig:uu} depicts the evolution of the velocity components.

As stated in the introduction, the system under consideration pertains to the migration of the cells along a greater concentration. From Figure \ref{fig:etac}, we observe that initially, the concentration is higher in the center of the domain, and the cells form two distinct clusters in the upper part of the domain. As time passes, the cells begin to orient their movement toward higher concentrations. It is observed that the clusters of organisms form a bridge between them and accumulate in the center of the domain. This is due to the effect of the cross-diffusion term (chemotaxis term), which is prominent in the early stages of movement. The chemical concentration decays with time, and the cross-diffusion term loses strength. As a result, cell self-diffusion becomes dominant, and cells begin to spread uniformly over the region. Finally, the cells migrate downwards of the domain due to the external forces acting on the fluid equation ($\nabla \Phi = (0,-1000)$). It is also observed that, as time progresses, there are some changes in velocity and pressure due to the movement of the cells as depicted in Figures \ref{fig:vecp} and \ref{fig:uu}.

% Eta and c

% Vector and p

% u1 and u2

\begin{figure}
	\centering
	\includegraphics[scale=0.55]{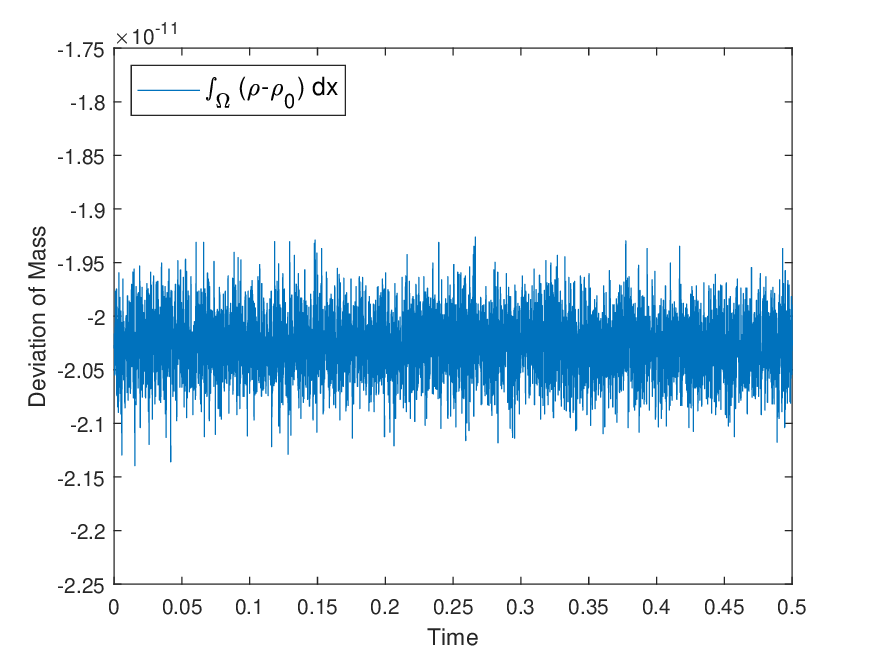}
	\includegraphics[scale=0.55]{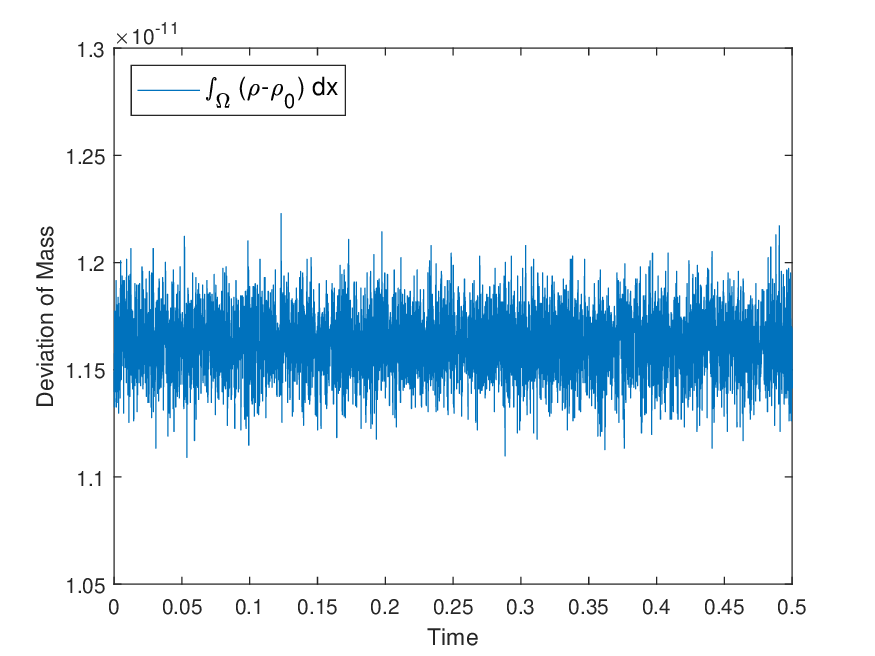}
	\caption{Deviation of mass of discrete cell density for linear(left) and quadratic(right) elements.}
	\label{fig:mass}
\end{figure}

Moreover, in Figure \ref{fig:mass}, we present the deviation of discrete mass for cell density which prove the global mass conservation property (Lemma \ref{lem:mass}) .
All the numerical simulations have been computed using the software FreeFem++ \cite{Hec12} and graphs are plotted in Matlab.

\section{Conclusion}
In this work, we have analyzed a fully discrete discontinuous Galerkin method for the coupled Chemotaxis (Keller-Segel)-Navier-Stokes system. We have obtained optimal error estimates in $L^2$-norm with $\mathcal{O}(h^{k+1}+\Delta t)$ and in $H^1$-norm with $\mathcal{O}(h^k+\Delta t)$ for the fully discrete cell density, concentration and fluid velocity and in $L^2$-norm with $\mathcal{O}(h^k+\Delta t)$ for the fluid pressure with the help of introducing a new projection for the chemotaxis term. The numerical simulations are performed which verify our theoretical findings.

\section*{Acknowledgments}
The first author would like to express his gratitude to the Indian Institute of Technology Bombay (IITB), for the financial support (Institute post-doctoral fellowship No. HR-1(HRM-1)/Rec/34/2022).

\bibliography{ref_CNSE}
\bibliographystyle{plain}

\end{document}